\documentclass[a4paper,reqno]{paper}
\usepackage[utf8]{inputenc}
\usepackage[T1]{fontenc}
\usepackage{amsfonts,enumerate}
\usepackage[dvipsnames]{xcolor}
\usepackage{amsmath}
\usepackage{amssymb}
\usepackage{amsthm}
\usepackage{graphics}
 \usepackage{array,multirow,graphicx}
\usepackage[backref = page]{hyperref}
\hypersetup{
colorlinks =  true, 
linkcolor = purple!80, 
citecolor =teal
}
\usepackage[capitalize]{cleveref}
\usepackage{tikz}
\usepackage{geometry}
\usepackage{amsfonts}
\usepackage{graphicx}
\usepackage{genyoungtabtikz}
\usepackage{dsfont}
\usepackage{mathrsfs}
\usepackage{stmaryrd}
\usepackage{mynewcommands}
\usepackage{booktabs}
\usepackage{ytableau}
\usepackage{float}
\usepackage{titlesec}

\Crefname{equation}{}{}

\titleformat{\subsection}
{\normalfont\fontsize{11}{17}\sffamily\bfseries}
{\thesubsection}
{1em}
{}

\let\svthefootnote\thefootnote
\newcommand\freefootnote[1]{
  \let\thefootnote\relax
  \footnotetext{#1}
  \let\thefootnote\svthefootnote
}

\begin{document}

\sloppy

\title{Atomic length on Weyl groups}
\author{
Nathan Chapelier-Laget\thanks{Institut Denis Poisson, Universit\'e de Tours, CNRS, Parc de Grandmont, 37200 Tours, France.
Email address: {\tt nathan.chapelier@gmail.com}
},
Thomas Gerber\thanks{\'Ecole Polytechnique F\'ed\'erale de Lausanne, 1015 Lausanne, Switzerland. Supported by the \textit{Ambizione} grant PZ00P2\_180120 of the Swiss National Science Foundation.
Email address: {\tt thomas.gerber@epfl.ch}
}
}
\maketitle

\hrule 

\begin{abstract} 
We define a new statistic on Weyl groups
called the \textit{atomic length} and investigate its combinatorial and representation-theoretic properties.
In finite types, we show a number of properties of the atomic length which 
are reminiscent of the properties of the usual length.
Moreover, we prove that, with the exception of rank two, 
this statistic describes an interval.
In affine types, our results shed some light on classical enumeration problems,
such as the celebrated Granville-Ono theorem on the existence of core partitions,
by relating the atomic length to the theory of crystals.
\end{abstract}

\hrule 

\setcounter{tocdepth}{1}

\section*{Introduction}
\addcontentsline{toc}{section}{\protect\numberline{}Introduction}

Let $W$ be a Coxeter group and $S$ be a finite set of generators called the \textit{simple reflections}.
The length function $\ell : W\to \N$ 
is a fundamental tool for studying the combinatorial and algebraic properties of $W$ and of related structures.
The number $\ell(w)$ is defined as the minimal number of 
simple reflections necessary to decompose $w$ 
(in particular, the length function depends on $S$).
Expressions of $w$ as a product of $\ell(w)$ generators are called \textit{reduced},
and Matsumoto's theorem ensures that any two reduced expressions of $w$ are related by a sequence of \textit{braid moves}.
It is readily shown that $\ell(w)$ equals the number of inversions of $w$,
that is, the number of positive roots of the corresponding root system that are sent to a negative root by $w^{-1}$.

\medskip

The length function naturally induces a partial order on $W$, which contains a famous ordering on $W$: the weak (Bruhat) order.
The weak order is ubiquitous in the theory of Coxeter groups and their deformations, the \textit{(Iwahori-)Hecke algebras}  \cite{GeckPfeiffer2000}.
Moreover,  it is contained in the (strong) Bruhat order, which plays a crucial role in Kazhdan-Lusztig theory 
and related topics \cite{lusztig2003hecke}, allowing for example to compare the (Zariski) adherence of cells in Schubert varieties.
There are also variants of the length function that also play an important role in representation theory  and combinatorics.

The \textit{twisted length} function is a central tool in the theory of certain partial pre-orders on $W$
defined with respect to a subset of reflections, and behaving similarly to the Bruhat order.
These partial pre-orders have been used to prove the existence of new Kazhdan-Lusztig polynomials \cite{Shelling2}. 
The twisted length function also helps determine which subsets of reflections induce a partial order \cite{Shelling2, edgar2007sets}.

In another direction, there is a notion of \textit{absolute length} or \textit{reflection length},
which is defined as the minimal number of reflections (not necessarily simple) required to decompose an element of $W$.
It also has many interesting properties, but is less well-understood than the usual length \cite{Carter1972}, \cite{Dyer2001}, \cite{LMPS2019}.

\medskip

In this paper, we introduce a different variant of the length function on Weyl groups denoted $\sL$, which we call the \textit{atomic length}.
More precisely, $\sL(w)$ is defined similarly to $\ell(w)$ by considering the inversion set of $w$,
but by counting each inversion of $w$ not just once, 
but as many times as its height (that is, the number of simple roots needed to decompose it).
In type $A$ (that is, when the Weyl group is the symmetric group), 
the atomic length turns out to coincide with half the \textit{entropy} of permutations, a well-studied statistic
(see \Cref{sec_entropy}).
We quickly realise that the atomic length is just a special case of 
a more general statistic depending on a parameter $\la$ which is a \textit{dominant weight}, 
denoted $\sL_\la$.
More precisely, specialising $\sL_\la$ at $\la=\rho$, the half-sum of positive roots, recovers the statistic $\sL$.
In fact, we will see that $\sL_\la$ can be defined
in the more general context of Weyl groups associated to Kac-Moody algebras \cite{Kac1984}, and we will focus on the affine case
(which contains all the finite types by restricting to an appropriate sub-root system).

\medskip

We will see that $\sL_\la$ enjoys a number of properties echoing properties of the usual length function $\ell$.
In particular, in finite types, the longest element $w_0\in W$ will play a crucial role.
We will prove that $w_0$ realises the largest value of the atomic length on $W$ (for any dominant weight $\la$),
and the main theorem of this paper states that 
\begin{center}
$\sL:W\to \lbra 0, \sL(w_0) \rbra$ is surjective,
 except for rank $2$ root systems.
\end{center}
The relevance of this result is two-fold.
First, our approach illuminates classic results about the entropy of permutations
by using systematic root systems combinatorics.
In particular, we recover the surjectivity of the half-entropy
on symmetric groups $S_n$ for $n\geq 4$ \cite{SackUlfarsson2011} as a special case.
Second, the affine variant of this problem leads us to a famous result by Granville and Ono \cite{GO1996},
ensuring that there exists an $(n+1)$-core partition of every size as soon as $n\geq 3$.
More precisely, take $W$ to be the affine Weyl group of type $A_n^{(1)}$.
Then, using a convenient interpretation of the statistic $\sL_\la$ in the context of crystals,
we are able to rephrase the main theorem of \cite{GO1996} as follows: 
\begin{center}
$\sL_{\La_0}: W \to \N$ is surjective if and only if $n\geq 3$.
\end{center}
Therefore, we ask the general question: when is the image of $\sL_\la$ an integer interval?
Like Granville and Ono's applications in block theory
for symmetric and alternating groups, 
this would enable us to understand defect zero blocks for more complicated structures such as (cyclotomic) Hecke algebras \cite{Fayers2019, JaconLecouvey2021}.
Also, we observe in \Cref{cores_paths} that the specialisation 
of $\sL_\la$ at $\la=\La_0$ in affine type $A$
already turned up in other interesting number-theoretic and probabilistic contexts \cite{GKS1990, BrunatNath2022, ThielWilliams2017, STW2021}.
We believe that further investigation of the properties of the atomic length in other types will reveal more exciting applications.

\medskip

The paper is structured as follows.
We motivate the study of the atomic length by reviewing a number of 
eclectic results about the entropy of permutations in \Cref{sec_entropy}.
In \Cref{sec_Cox}, we introduce the necessary background on Coxeter groups, root systems and reflection subgroups,
and we also recall some important features in the Weyl group case.
In \Cref{sec_affine_WG}, we use the setting of Kac-Moody algebras to introduce Weyl groups of affine type in full generality,
and we recall how to recover the untwisted affine Weyl groups in their classical geometric construction.
This formalism enables us to introduce the notion of atomic length $\sL$ in \Cref{sec_AL}, firstly using inversion sets, 
and secondly using a more general formula involving the parameter $\la$, giving rise to the statistic $\sL_\la$.
We prove several important properties of $\sL_\la$, reminiscent of classical properties of the usual length function.
\Cref{sec_susanfe} is devoted to the study of special reflections of $W$ which we call \textit{Susanfe}.
In classical types, we introduce particular Susanfe elements, which are subsequently used in \Cref{sec_surj_AL}.
In fact, we give two independent proofs of \Cref{thm_surj}
which states that $\sL$ surjects onto the expected integer interval (except in rank $2$), 
one of which extensively using the properties of $w_0$, and the other based on Susanfe theory.
In \Cref{sec_cryst}, we explain how the atomic length $\sL_\la$ can be interpreted in the theory of crystals 
for Kac-Moody algebras representations, namely by simply looking at the depth  of certain $W$-orbit elements
in the corresponding crystal graph.
This enables us to rephrase the Granville-Ono theorem
in terms of atomic length, and thereby motivates the study of the affine atomic length, which we initiate in \Cref{sec_affine_AL}.


\section{Entropy of permutations}\label{sec_entropy}

In this section, we survey some results
about several important statistics on the symmetric group,
including the notion of entropy.
This serves as a motivation for studying the atomic length defined in \Cref{sec_AL}, as it will turn out to be simply half the entropy in type $A$.
We denote $\N$ the set of nonnegative integers.

\subsection{Entropy and inversion sum}\label{subsec_entropy_inv}

Let $S_n$ be the symmetric group and let ${{w}} \in S_n$. Denote 
$$
\cos({{w}}) = \sum_{k=1}^n k \ {{w}}(k).
$$
The map $\cos: S_n\to\N$ is called the \textit{cosine}.
In \cite[Theorem 2.2]{SackUlfarsson2011}, it was proved that
all nonnegative integers with the exception of $16$ can be expressed as the cosine of some permutation.
Recently, the cosine map has been studied from a combinatorial point of view in \cite[Section 5.2]{elder2022homomesies} where the authors show that it is a $\frac{n(n+1)^2}{4}$-homomesy.
In order to prove \cite[Theorem 2.2]{SackUlfarsson2011},  
the authors introduce two statistics on $S_n$,
based on the notion of inversion and non-inversion of a permutation. 
Let ${{w}} \in S_n$ and write ${{w}} = {{w}}_1 {{w}}_2\dots {{w}}_n$ in one-line notation,
that is $ {{w}}(k) = {{w}}_k$. 
An inversion of ${{w}}$ is a pair $(i,j)$ such that $i < j$ and ${{w}}_i > {{w}}_j$,
and a non-inversion of ${{w}}$ is a pair $(i,j)$ such that $i < j$ and ${{w}}_i < {{w}}_j$. 
The set of all inversions (respectively non-inversions) of $w$ is denoted by 
$\boldsymbol{N}({{w}})$ (respectively $\boldsymbol{N'}({{w}})$).  
We then set
$$
\mathsf{invsum}({{w}}) = \sum\limits_{(i,j) \in \boldsymbol{N}({{w}})}(j-i)
\mand
\mathsf{ninvsum}({{w}}) = \sum\limits_{(i,j) \in \boldsymbol{N'}({{w}})}(j-i).
$$
We will see that $\mathsf{invsum}$
coincides with the type $A_{n-1}$ atomic length of \Cref{sec_AL}, see \Cref{rem_def_AL} (2).
Furthermore, we have the following formula  \cite[Eq.(1)]{SackUlfarsson2011}  that connects these two statistics
\begin{align}\label{formula 2 stats}
\mathsf{invsum}({{w}}) + \mathsf{ninvsum}({{w}}) = \binom{n+1}{3}.
\end{align}

Let $w_0$ be the permutation $n(n-1)(n-2)\cdots 1$\footnote{
The permutation $w_0$ is the \textit{longest element} of $S_n$ in the Coxeter, see \Cref{sec_Cox}, which justifies the notation.
}.
An important tool for proving \cite[Theorem 2.2]{SackUlfarsson2011} is the following theorem.

\begin{Th}[{\cite[Theorem 2.5]{SackUlfarsson2011}}]\label{Theo Sack 2.5}
For ${{w}} \in S_n$ we have
$$
 \cos({{w}}) = \cos(w_0 ) +  \mathsf{ninvsum}({{w}}).
 $$
\end{Th}

We are ready to give the central definition of this section.

\begin{Def}\label{def_entropy}
Let $w\in S_n$. The entropy of ${{w}}$ is the nonnegative integer
$
E({{w}}) = \sum\limits_{i=1}^n(i-{{w}}(i))^2.
$
\end{Def}

\begin{Rem}\label{rem_entropy}\
\begin{enumerate}
\item 
The entropy is a particular case of the metric $S$ introduced in \cite{diaconis1977spearman}, 
showing that it is the ``Spearman’s rho'' of a permutation and the identity permutation.
\item 
The entropy is related to another important notion: \textit{bigrassmannian} permutations.  
A bigrassmannian permutation is a permutation that has 
only one left descent and only one right descent (in the Coxeter sense).
The study of these permutations goes way back, beginning with the work of Lascoux and Schützenberger \cite{lascoux1996treillis} where they show that ${{w}} \in S_n$ is bigrassmannian  if and only if it is 
join-irreducible for the Bruhat order on $S_n$ \cite[Theorem 4.4]{lascoux1996treillis}.  
The notion of join-irreducible elements in Weyl groups appears more generally in quiver-representation-theoretic context in \cite{iyama2018lattice}.

Precise computations and expressions of bigrassmannian permutations can also be found in \cite{geck1997bases}.  It turns out that various aspects of bigrassmannian permutations are useful in a broad range of subjects and therefore were strongly studied in the past three decades \cite{eriksson1996combinatorics,  lascoux1996treillis,  geck1997bases,  reading2002order,  kobayashi2010bijection,  kobayashi2011enumeration,  reiner2011presenting, engbers2015comparability,  iyama2018lattice}.  
For example, in \cite{reiner2011presenting},  Reiner uses this notion to describe the cohomology of Schubert varieties. 
In another direction,  bigrassmannian elements are used to determine the socle of the cokernel of an inclusion of Verma modules in type $A$ \cite[Section 2.2, Corollary 5]{ko2021bigrassmannian}. 
Let $\leq_B$ denote the (strong) Bruhat order on $S_n$.
In \cite[Theorem]{kobayashi2011enumeration} Kobayashi shows that the cardinality of the set $\{v \leq_B {{w}}~|~ v~\text{is bigrassmannian}\}$ is equal to $E({{w}})/2$.
\end{enumerate}
\end{Rem}

The entropy is related to the inversion sum by the following formula, announced previously.

\begin{Prop}\label{entropy_vs_invsum}
Let ${{w}} \in S_n$. We have
$$
\displaystyle\frac{E({{w}})}{2} = \mathsf{invsum}({{w}}).
$$
\end{Prop}

\begin{proof}
Let us denote $a_n = \sum_{i=1}^n i^2$.  On the one hand we have
$$
E({{w}}) =  \sum\limits_{i=1}^n(i-{{w}}(i))^2 =  
\sum\limits_{i=1}^n (i^2-2i{{w}}(i)+{{w}}(i)^2) = 
\sum\limits_{i=1}^n i^2 -2\sum\limits_{i=1}^n i {{w}}(i) + \sum\limits_{i=1}^n {{w}}(i)^2 = 
2a_n -2\cos({{w}}).
$$
Hence $E({{w}})/2  = a_n - \cos({{w}})$.  
Moreover by \Cref{Theo Sack 2.5}, 
we know that $ \cos({{w}}) = \cos(w_0 ) +  \mathsf{ninvsum}({{w}})$,
and by \Cref{formula 2 stats}, we know that 
$\mathsf{invsum}({{w}}) + \mathsf{ninvsum}({{w}}) = \binom{n+1}{3}$.  
By definition of the cosine we have
$$
\cos(w_0 ) = \sum\limits_{i=1}^n i \ w_0 (i) = \sum\limits_{i=1}^n i(n+1-i) = (n+1) \sum\limits_{i=1}^n i -  a_n.
$$
Therefore we obtain
\begin{align*}
a_n - \cos({{w}}) & = a_n - \cos( w_0 ) - \mathsf{ninvsum}({{w}}) \\
									 & = a_n - \left[(n+1) \sum\limits_{i=1}^n i -  a_n \right]- \left[ \binom{n+1}{3}-\mathsf{invsum}({{w}})  \right] \\
									 & = 2a_n - (n+1)\sum\limits_{i=1}^n i -  \binom{n+1}{3} + \mathsf{invsum}({{w}}) \\
									 & = \frac{n(n+1)(2n+1)}{3} - \frac{n(n+1)^2}{2} - \frac{(n-1)n(n+1)}{6} + \mathsf{invsum}({{w}}) \\
									& =  \mathsf{invsum}({{w}}).
\end{align*}
\end{proof}

\subsection{Entropy and the permutohedron}

A geometrical interpretation of the entropy is given in terms of the $n$-dimensional permutohedron $\mathcal{P}_n$ 
(called the ``Voronoi cell'' of a certain lattice in \cite[Theorem 7 page 474]{conway2013sphere}).
Let us recall the construction of the permutohedron in type $A_{n-1}$.

\medskip 

Fix $V = \mathbb{R}^n$ with canonical basis $\{e_i~|~i=1,\dots, n\}$.  
The symmetric group $S_n$ acts transitively
on $V$ by permuting the coordinates,
denote $S_n\times V\to V, (w,x)\mapsto w(x)$ this action.
The reflection hyperplanes of $S_n$ are $H_{ij} = \{x \in V~|~x_i-x_j = 0\}$ with $1 \leq i < j \leq n$. 
Choose now a generic point $x \in V$ not located on any reflection hyperplane of $S_n$.  
Then the permutohedron 
is defined as the convex hull of the orbit under $S_n$ of the point $x$, that is
$$
\mathcal{P}_n(x) = \mathrm{conv}\{ w(x)~|~ w \in S_n \}.
$$
We see then that there are as many realizations of $\mathcal{P}_n$ as points $x$ in $V$ not located on the reflection hyperplanes.  Some of them are move convenient and sometimes it can be useful to choose one over another.  We say that a point $x \in V$ is \textit{adequate} if it is not located on any reflection hyperplane of $S_n$
and if all its coordinates $x_i$ are in $[1, n]\cap\Z$.

\begin{Rem}
The permutohedron can be defined for any finite Coxeter group, see \cite[Section 1.5]{pilaud2020facial} for a good reference on the subject.
\end{Rem}

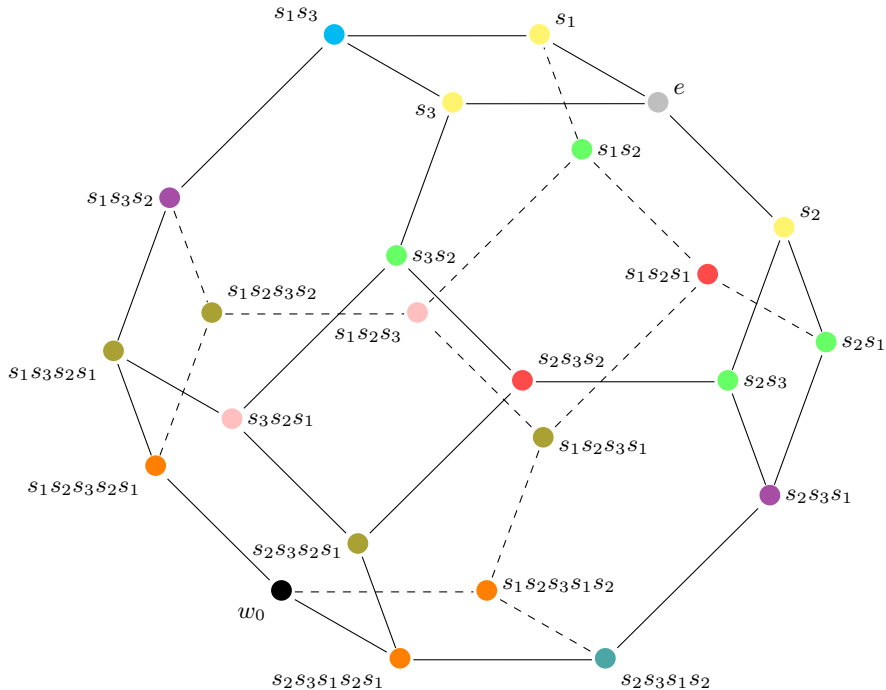
\begin{figure}[h!]
\centering
\begin{tikzpicture}[scale=1.8, every node/.style={inner sep=0,outer sep=0}]

\node[anchor=mid, scale=2](black) at (0,0) {$\color{black}\bullet$};
\path (black) ++(-30:1) node[scale=2] (or1) {$\color{orange}\bullet$};
\path (black) ++(0:1.5) node[scale=2] (or2) {$\color{orange}\bullet$};
\path (or1) ++(0:1.5) node[scale=2] (bl1) {$\color{teal!70}\bullet$};
\path (black) ++(135:1.3) node[scale=2] (or3) {$\color{orange}\bullet$};
\path (or3) ++(110:0.9) node[scale=2] (pu1) {$\color{olive!80}\bullet$};
\path (or1) ++(110:0.9) node[scale=2] (pu2) {$\color{olive!80}\bullet$};
\path (pu1) ++(-30:1) node[scale=2] (pink1) {$\color{pink}\bullet$};
\path (pu1) ++(70:1.2) node[scale=2] (br1) {$\color{violet!70}\bullet$};
\path (or3) ++(70:1.2) node[scale=2] (pu3) {$\color{olive!80}\bullet$};
\path (pink1) ++(45:1.7) node[scale=2] (gre1) {$\color{green!60}\bullet$};
\path (pu2) ++(45:1.7) node[scale=2] (gra1) {$\color{red!71}\bullet$};
\path (bl1) ++(45:1.7) node[scale=2] (br2) {$\color{violet!70}\bullet$};
\path (br2) ++(110:0.9) node[scale=2] (gre2) {$\color{green!60}\bullet$};
\path (gre1) ++(70:1.2) node[scale=2] (pi1) {$\color{yellow!70}\bullet$};
\path (pi1) ++(150:1) node[scale=2] (cy1) {$\color{cyan!80}\bullet$};
\path (pi1) ++(0:1.5) node[scale=2] (yel1) {$\color{lightgray}\bullet$};
\path (br2) ++(70:1.2) node[scale=2] (gre3) {$\color{green!60}\bullet$};
\path (gre2) ++(70:1.2) node[scale=2] (pi3) {$\color{yellow!70}\bullet$};
\path (yel1) ++(150:1) node[scale=2] (pi2) {$\color{yellow!70}\bullet$};
\path (gre3) ++(150:1) node[scale=2] (gra2) {$\color{red!71}\bullet$};
\path (gra2) ++(135:1.3) node[scale=2] (gre4) {$\color{green!60}\bullet$};
\path (gra2) ++(225:1.7) node[scale=2] (pu4) {$\color{olive!80}\bullet$};
\path (gre4) ++(225:1.7) node[scale=2] (pink2) {$\color{pink}\bullet$};

\node[anchor=north east, xshift=-2mm, yshift=-2mm] at (0,0) {\small $w_0$};
\node[anchor=south west, xshift=2mm, yshift=0mm] at (or2) {\small $s_1s_2s_3s_1s_2$};
\node[anchor=north east, xshift=-2mm, yshift=-2mm] at (or1) {\small $s_2s_3s_1s_2s_1$};
\node[anchor=north east, xshift=-2mm, yshift=-2mm] at (or3) {\small $s_1s_2s_3s_2s_1$};
\node[anchor=north east, xshift=-2mm, yshift=-2mm] at (pu1) {\small $s_1s_3s_2s_1$};
\node[anchor=north east, xshift=-2mm, yshift=-0mm] at (pu2) {\small $s_2s_3s_2s_1$};
\node[anchor=south west, xshift=2mm, yshift=-2mm] at (pu4) {\small $s_1s_2s_3s_1$};
\node[anchor=south west, xshift=2mm, yshift=2mm] at (pu3) {\small $s_1s_2s_3s_2$};
\node[anchor=north west, xshift=2mm, yshift=-2mm] at (bl1) {\small $s_2s_3s_1s_2$};
\node[anchor=north east, xshift=-2mm, yshift=-2mm] at (pink2) {\small $s_1s_2s_3$};
\node[anchor=west, xshift=2mm, yshift=-0mm] at (pink1) {\small $s_3s_2s_1$};
\node[anchor=east, xshift=-2mm, yshift=-0mm] at (br1) {\small $s_1s_3s_2$};
\node[anchor=west, xshift=2mm, yshift=-0mm] at (br2) {\small $s_2s_3s_1$};
\node[anchor=west, xshift=2mm, yshift=-0mm] at (gre1) {\small $s_3s_2$};
\node[anchor=west, xshift=2mm, yshift=-0mm] at (gre2) {\small $s_2s_3$};
\node[anchor=west, xshift=2mm, yshift=-0mm] at (gre3) {\small $s_2s_1$};
\node[anchor=west, xshift=2mm, yshift=-0mm] at (gre4) {\small $s_1s_2$};
\node[anchor=south east, xshift=-2mm, yshift=2mm] at (cy1) {\small $s_1s_3$};
\node[anchor=south west, xshift=2mm, yshift=2mm] at (gra1) {\small $s_2s_3s_2$};
\node[anchor=east, xshift=-2mm, yshift=-0mm] at (gra2) {\small $s_1s_2s_1$};
\node[anchor=east, xshift=-2mm, yshift=-1mm] at (pi1) {\small $s_3$};
\node[anchor=west, xshift=2mm, yshift=2mm] at (pi3) {\small $s_2$};
\node[anchor=west, xshift=2mm, yshift=2mm] at (pi2) {\small $s_1$};
\node[anchor=west, xshift=2mm, yshift=2mm] at (yel1) {\small $e$};

\draw[] (black) -- (or1);
\draw[] (black) -- (or3);
\draw[dashed] (black) -- (or2);
\draw[] (or3) -- (pu1);
\draw[] (or1) -- (pu2);
\draw[] (pu1) -- (pink1);
\draw[] (pu2) -- (pink1);
\draw[] (pu1) -- (br1);
\draw[dashed] (or3) -- (pu3);
\draw[dashed] (br1) -- (pu3);
\draw (pink1) -- (gre1);
\draw (pu2) -- (gra1);
\draw (gre1) -- (gra1);
\draw (gre2) -- (gra1);
\draw (gre2) -- (pi3);
\draw (pi3) -- (yel1);
\draw (pi2) -- (yel1);
\draw (pi1) -- (yel1);
\draw (pi1) -- (cy1);
\draw (pi2) -- (cy1);
\draw (br1) -- (cy1);
\draw (gre1) -- (pi1);
\draw (pi3) -- (gre3);
\draw (br2) -- (gre3);
\draw (br2) -- (gre2);
\draw (br2) -- (bl1);
\draw[dashed] (or2) -- (bl1);
\draw (or1) -- (bl1);
\draw[dashed](pu3)--(pink2);
\draw[dashed](pu4)--(pink2);
\draw[dashed](gre4)--(pink2);
\draw[dashed](gre4)--(gra2);
\draw[dashed](pu4)--(gra2);
\draw[dashed](gre3)--(gra2);
\draw[dashed](gre4)--(pi2);
\draw[dashed](pu4)--(or2);

\end{tikzpicture}
\caption{The permutohedron for $S_4$.  The vertices are labelled by reduced expressions of the group elements (i.e., in terms of the Coxeter generators).
Vertices of the same color are at the same distance to $e$,
therefore have the same entropy.
The values of the half-entropy (which will coincide with the atomic length) 
for each color is given by
\newline
\centerline{
${\color{lightgray}\bullet}:0$ \quad
${\color{yellow!70}\bullet}:1$ \quad
${\color{cyan!80}\bullet}:2 $  \quad
${\color{green!60}\bullet}:3 $  \quad
${\color{red!71}\bullet}:4 $  \quad
${\color{violet!70}\bullet}:5 $  \quad
${\color{pink}\bullet}:6 $  \quad
${\color{olive!80}\bullet}:7 $  \quad
${\color{teal!70}\bullet}:8 $  \quad
${\color{orange}\bullet}:9 $  \quad
${\color{black}\bullet}:10 $
}
}
\label{permutohedron}
\end{figure}

We can give the expected formula relating the entropy and the permutohedron.
The following proposition states that the entropy of a permutation $w$ is the square of the distance
between two points of $\mathcal{P}_n(x)$ permuted by $w$.
This is illustrated in \Cref{permutohedron}.

\begin{Prop}
Let $x\in V$ be adequate, and let ${{w}} \in S_n$.
Then 
$$
E({{w}})= | w(x) - x|^2.
$$
\end{Prop}

\begin{proof}
Write $x = (x_1,x_2,\dots, x_n)$.  By definition $w(x) = ({{w}}(x_1), {{w}}(x_2),\dots, {{w}}(x_n))$.  Moreover we have the following formula $| w(x) - x|^2 = \sum_{i=1}^n ({{w}}(x_i) - x_i)^2$.  
Since $x$ is adequate, we can make the change of variables $x_i \mapsto i$ and we obtain 
$| w(x) - x|^2 = \sum_{i=1}^n ({{w}}(i) - i)^2 = E({{w}})$, which ends the proof.
\end{proof}

\section{Coxeter groups and Weyl groups}
\label{sec_Cox}

\subsection{Inversion sets}

Let $(W,S)$ be a Coxeter system with $S$ finite.
Consider a  corresponding root system\footnote{The geometric representation of $(W,S)$ gives a standard way to define root systems for any Coxeter system,  see for example \cite[Section 2.3]{SRLE}.} $\Phi$
and denote $\Delta=\{\al_1,\ldots, \al_n\}$ the set of simple roots, so that $S=\{s_1,\ldots, s_n\}$
where $s_i=s_{\al_i}$ is the \textit{simple reflection} associated to the simple root $\al_i$. 
Throughout this paper, we shall always assume that $\Phi$ is irreducible.
Let $\Phi^+ \subseteq \Phi$ 
be the associated set of positive roots, so that
$\Phi^- = -\Phi^+$ is the corresponding set of negative roots.
In particular, $\Phi = \Phi^+ \sqcup \Phi^-$.  
We denote by $T = \bigcup_{w \in W}wSw^{-1}$ the set of all reflections of $W$.  
The set $\Phi^+$ is in bijection with $T$, as any reflection 
$t$ writes $t=s_{\al}$ for a certain $\al \in \Phi^+$.

\medskip

Let ${\ell : W \to \N}$ be the length function, that is, 
$\ell(w)$ is the smallest number $r$ such that there exists an expression $w = s_{i_1}\dots s_{i_r}$ with ${s_{i_k} \in S}$. By convention, $\ell(e) = 0$. An expression of $w \in W$ is called a reduced expression if it is a product of $\ell(w)$ generators.

\medskip

Let $w \in W$. The inversion set of $w$ is
\begin{align*}
 N(w) & = \{ \al \in \Phi^+\mid w^{-1}(\al) \in \Phi^- \},
\end{align*}
and an alternative description of the inversion set is given by
\begin{align*}
N(w) = \{ \al \in \Phi^+\mid \ell(s_{\al}w) < \ell(w) \}.
\end{align*}

It is well-known that 
$\ell(w) = |N(w)|$, 
and $w = w'$ if and only if $N(w) = N(w')$. 
The following result, 
which can be found for example in \cite[Proposition 2.1]{ISWO},
enables us to construct inversion sets starting from reduced expressions.
We will use it several times in the rest of the paper.

\begin{Prop}\label{Decomposition_N}
Let $w \in W$ 
and write $w=uv$ with $u, v\in W$ veriyfing $\ell(w) = \ell(u)+\ell(w)$.
Then $N(w) = N(u) \sqcup u(N(v))$. 
In particular if $w = s_1s_2\dots s_{n-1}s_n$ is a reduced expression then 
$N(w) = \{ \al_1, s_1(\al_2), s_1s_2(\al_3),\dots, s_1s_2\dots s_{n-1}(\al_n) \}.$
\end{Prop}

In fact, inversion sets are well-behaved with respect to the 
right weak order, defined as follows. For $ w, w'\in W$, we write $w \leq w'$ if there exist reduced expressions
$\underline{w}, \underline{w'}$ of $w$ and $w'$ respectively such that $\underline{w}$ is a prefix of $\underline{w'}$
We have the following characterisation proved in \cite[Corollary 2.10]{SRLE}.

\begin{Prop}\label{inversions_weak_order_1}
For all $w,w'\in W$, we have
 $N(w)\subseteq N(w') \Leftrightarrow w\leq w'$.
\end{Prop}

In fact, one direction of the above equivalence can be refined as follows, see \cite[Proposition 1.1]{ISWO}.

\begin{Prop}\label{inversions_weak_order_2}
The map $N$ is a poset monomorphism from $(W, \leq)$ to ${(\mathcal{P}(\Phi^+), \subseteq)}$. 
\end{Prop}

We now give some properties that will be needed in the following sections.
Recall that if $A \subseteq T$ is a subset of reflections, then the corresponding reflection subgroup of $W$ is the subgroup 
$W_A$ generated by the reflections in $A$, that is $$
W_A =  \langle s_{\al} \mid s_{\al} \in A \rangle .
$$
From \cite{Reflection_subgroups} we know that the set $\Phi_{A} = \{ \al \in \Phi\mid s_{\al} \in W_A\}$ is a root system of $W_A$, 
with simple roots $\Delta_{A} = \{ \al \in \Phi^+\mid N(s_{\al}) \cap \Phi_{A}=\{\al\}\}$.
Moreover, we have $\Phi_{A}^+ = \Phi_{A} \cap \Phi^+$. Following \cite{PC}, we denote ${S_A = \{s_{\al}\mid \al \in \Delta_{A} \}}$. Therefore, the pair $(W_A, S_A)$ is a Coxeter system, and we denote $\ell_A$ its length function.

 Let us write 
 \begin{align}\label{A-coset}
 {}^AW & = \{w \in W\mid \ell(s_{\al}w) > \ell(w)~~\forall \al \in \Phi_A \}   = \{w \in W\mid N(w) \cap  \Phi_A = \emptyset \}.
 \end{align}
Note finally that there is an easier characterization of the previous set thanks to the functoriality of the Bruhat graph (see \cite[Section 2.5]{SRLE} for more details):
 \begin{align}\label{functoriality}
 {}^AW  = \{w \in W\mid \ell(s_{\al}w) > \ell(w)~~\forall \al \in \Delta_A \}.
 \end{align}
It is also known that for any $w \in W$, there exists a unique $w_A \in W_A$ and a unique ${}^Aw \in {}^AW$ such that $w = w_A{}^Aw$. The decomposition $w = w_A{}^Aw$ is called the $A$-decomposition of $w$. 

\begin{Rem}\label{rem_A_dec}\
\begin{enumerate}
 \item The $A$-decomposition defined above is usually called ``left $A$-decomposition'',
and accordingly, there is a notion of right $A$-decomposition that uses an analogous subset $W^A$. 
In this article we will only use the left $A$-decomposition, which is why we use this simpler terminology.
\item 
In the particular case  $A \subseteq S$, 
we recover the parabolic decomposition, see \cite{Hum}.
In this case, we have the formula 
$\ell(w) = \ell(w_A)+ \ell({}^Aw)$ but this does not holds in general if $A \not\subseteq S$, see \cite{BG}. 

\end{enumerate}

\end{Rem}

If $W_A$ is a reflection subgroup of $W$ and $w\in W$, the group $wW_Aw^{-1}$ is also a reflection subgroup of $W$,
generated by the reflections $wtw^{-1}, t\in A$.
For simplicity, we denote $B= \{ wtw^{-1} \,;\, t\in A\}$, so that $\Phi_B$ denotes the root system corresponding to $W_B=wW_Aw^{-1}$.
Recall that we have the following proposition.

\begin{Prop}\label{root sytem reflection}
Let $A\subseteq T$ and $w\in W$, and set $B=wAw^{-1}$.
We have $\Phi_B 
= w(\Phi_A)$.
\end{Prop}

\begin{Def}
Let $A \subseteq T$ and $w \in W$.  
Define $W_B = wW_Aw^{-1}$.
We say that $w$ is \textit{$A$-utopic} if the following map is a bijection
$$
\begin{array}{ccccc}
\Gamma_w & : & W_B & \longrightarrow & W_A \\
                    &   &   x      & \longmapsto      & (wx)_A.
\end{array}
$$
The set of $A$-utopic elements is denoted by $\mathcal{U}(A)$.
\end{Def}

\begin{Prop}\label{Proposition_Utopic}
Let $I \subseteq S$. Then $T \subseteq  \mathcal{U}(I)$.
\end{Prop}
 
 \begin{proof}
 Let $t \in T$ and $W_B = tW_It$.  Let us show that the map $\Gamma_t$  is a bijection. To do so, since $W_B$ and $W_I$ have the same number of elements it is enough to show that $\Gamma_t$ is injective.  Let $x$ and $y$ be two elements of $W_B$. Therefore, there exist $u, v \in W_I$ such that $x = tut$ and $y = tvt$. Assume that $\Gamma_t(x) = \Gamma_t(y)$, that is $(tx)_I = (ty)_I$, which is equivalent to $(ut)_I = (vt)_I$ and then $(ut_I {}^It)_I = (vt_I {}^It)_I$.  But since $u$ and $v$ belong to $W_I,$ it follows that $(ut_I {}^It)_I = ut_I$ and $(vt_I {}^It)_I = vt_I$. Thus $ut_I = vt_I$ and it follows that $u=v$.  Hence $x=y$, which ends the proof.
 \end{proof}

For $w \in W_A$ we denote by $N_{A}(w)$ its inversion set, 
which is a subset of $\Phi_A^+$.  
The relation between $N$ and $N_{A}$ is discussed in Lemma \ref{lemma Shi coeff dans J}.
It turns out that inversion sets behave nicely when restricted to reflection subgroups, as expressed in the following proposition,
found in \cite[Proposition 2.16]{SRLE}.
 
\begin{Prop}\label{propDy}
 Let $W_A$ be a reflection subgroup of $W$ and let $w \in W$. 
 Write $w=w_A{}^Aw$ be the $A$-decomposition of $w$. Then $N(w) \cap \Phi_A = N_{A}(w_A)$.
 \end{Prop}

 \subsection{The case of finite Weyl groups}
 \label{Weyl_groups}
 
 Let $V$ be a Euclidean space with inner product $(\cdot\mid \cdot)$. 
 Let $\Phi$ be an irreducible crystallographic root system in $V$, that is 
 $2\frac{(\al \mid \beta)}{(\al \mid \al)} \in \mathbb{Z}$ for any $\al, \beta \in \Phi$.
 We denote again $\De=\{\al_1,\ldots, \al_n\}$ a simple system and $\Phi^+$ the corresponding positive roots.
 The corresponding reflection group is called a \textit{Weyl group}.
 These are particularly important as they are attached
 to simple Lie algebras over $\C$, for which we have a classification.
 This is achieved via the different \textit{Dynkin} types, namely
 $A_n$ for $n\geq 1$, $B_n, C_n$ for $n\geq 2$, $D_n$ for $n\geq 4$ (the classical types),
 $E_6,E_7, E_8. F_4$ and $G_2$ (the exceptional types).
 We refer to the book \cite{BOURB} for details, and we will use its conventions in this paper.
 We will sometimes use the convenient notation $W(X_n)$
 for the Weyl group associated to a root system of type $X_n$.
 In this paper, we will mostly focus on Weyl groups rather than general Coxeter groups, 
 since we will make use of several representation-theoretic properties.
 
 \medskip
 
Let us fix some notation that will be useful in the rest of this paper.
First of all the \textit{height} of a root $\al\in\Phi$ is the number of simple roots appearing in the decomposition of $\al$, that is,
if $\al=\sum_{i=1}^n a_i\al_i$ with $a_i\in\Z$, then
$$\h(\al) = \sum_{i=1}^n a_i.$$
We denote by $e_i$  the $i$-th canonical vector of its corresponding ambient space $V$. 
We will use the shorthand notation
$e_{ij} = e_i - e_j$  and $e^{ij} = e_i + e_j$.
In clasical types, the roots are either $e_{ij}$, $e^{ij}$, $e_i$ or $2e_i$ and we refer the reader to read \cite[Planches]{BOURB} for more details. 
Sometimes we will use a comma between the labels $i$ and $j$ to prevent any ambiguity, for instance we will write $e^{i,n+1}$ rather than $e^{i(n+1)}$.

\begin{Rem}\label{rem_inv_type_A}
In the case where $W=W(A_n)$, that is $W=S_{n+1}$,
the inversion set $N(w)$ is in bijection with the set $\boldsymbol{N}({{w}})$
introduced in \Cref{sec_entropy} via the map $e_{ij}\mapsto (i,j)$.
In this case, the simple roots are $e_{i,i+1}$ for $1\leq i\leq n$,
and one checks that $\h(e_{ij})=j-i$, which is the quantity that appears in the definition of
the statistic $\mathsf{invsum}$.
\end{Rem}

Finally, we will use two important features of Weyl groups.
First, there exists a unique element $w_0$ of maximal length,
see for instance \cite[Chapter VI, \S 1, Corollary 3]{BOURB}.
Second, there exists a unique root $\tal$
verifying, for all $\be\in\Phi$, $\tal-\beta=\sum_{i=1}^n a_i \al_i$ for some nonnegative integers $a_1,\ldots, a_n$, see \cite[Chapter VI, \S 1, Proposition 25]{BOURB}. It follows from this property that $\tal\in\Phi^+$ and that $\tal$ is in fact
the unique element with maximal height. Thus, $\tal$ is called  the \textit{highest root} of $\Phi$.


\section{Affine Weyl groups and Shi coefficients}
\label{sec_affine_WG}

In this section, we recall the construction of affine Weyl groups
via the theory of Kac-Moody algebras, following \cite{Kac1984}.
The Weyl groups appearing in this context are either of untwisted or twisted type.
In the untwisted cases, we will recall the classical geometric constructions in \Cref{subsec_shi}.

\subsection{Kac-Moody algebras and affine Weyl groups}

Let $A = (a_{ij})_{0\leq i, j\leq n}$
be a generalised (symmetrisable) Cartan matrix such that
the corresponding Kac-Moody algebra $\mathfrak{g}$ is of affine type.
In particular, the rank of $A$ is $n$.
Let $\mathfrak{h}^\ast$
be a real vector space of dimension $n+2$
and $\Delta=\{ \al_0,\ldots, \al_n\}\subseteq  \mathfrak{h}^\ast$,
$\Delta^\vee=\{ \al_0^\vee,\ldots, \al_n^\vee\}\subseteq  \mathfrak{{h}}$
be a realisation of $\mathfrak{g}$ (see \cite{Kac1984} for details), so that
$\langle \al_j, \al_i^\vee \rangle = a_{ij},$
where $\langle . , .\rangle$ denotes the natural pairing between $\mathfrak{h}^\ast$ and $\mathfrak{{h}}$.  Write 
$V_0=\bigoplus_{i=1}^n \R\al_i$ and $V = V_0\oplus \R\al_0$ together with $V_0^*=\bigoplus_{i=1}^n \R\al_i^{\vee}$ and $V^* = V_0^*\oplus \R\al_0^{\vee}$.  
For $0\leq i\leq n $, the reflections $s_i : \mathfrak{h}^\ast \to \mathfrak{h}^\ast$ given by the formula
$$s_i(x) = x -\langle x, \al_i^\vee \rangle \al_i$$
generate a subgroup $W$ of $\mathrm{GL}(\mathfrak{h}^\ast)$ called the \emph{Weyl group}  of $\mathfrak{g}$.  The subgroup of $W$ generated by $s_1,\dots, s_n$ is called the \emph{finite Weyl group} of $\mathfrak{g}$ and we denote it by $W_0$.  
Finally, we write $\De_0 =\{\al_1,\ldots, \al_n\}$ and $\De_0^{\vee} =\{\al_1^{\vee},\ldots, \al_n^{\vee}\}$.

\medskip

The Weyl group $W$ acts naturally on $\mathfrak{h}^\ast$ by isometry and induces an action on $\mathfrak{{h}}$ via the contragredient representation,  defined by $w(f) = f \circ w^{-1}$ f
or any $f \in \mathfrak{{h}}$.  In particular, one can show \cite[Proposition 16.14]{Carter2005} that this action is determined by the formulas $s_i(f) = f - \langle f,  \alpha_i \rangle \alpha_i^{\vee}$
for all $0\leq i \leq n$. Moreover, the finite Weyl group $W_0$ acts by restriction on $V_0$ and $V_0^*$.
We denote by $\Phi_0$ (respectively  $\Phi_0^{\vee}$) the crystallographic root system (respectively coroot system) associated to $W_0$, that is $\Phi_0 = W_0(\Delta_0)$ (respectively  $\Phi_0^{\vee} = W_0(\Delta_0^{\vee})$.  For $\alpha=w(\alpha_i)\in \Phi_0$, where $\alpha_i \in \Delta_0$,
we denote $\alpha^{\vee} = w(\alpha_i^{\vee})$ the corresponding coroot. 
The definition of $\alpha^{\vee}$ is independent of the choice of $w$.
Denote $\Phi = W(\Delta)$. Then $\Phi$ is a root system associated to $W$ with simple system $\Delta$.
The set $Q =  \Z \De_0$ is called the \textit{root lattice} and $Q^{\vee} = \Z \De_0^{\vee}$ is called the \emph{coroot lattice}.
We can extend the notion of height to any element of $Q$ by setting,
for $\be = \sum_{i=1}^n b_i\al_i\in Q$,
$\h(\be)=\sum_{i=1}^n b_i$.

\medskip

We now introduce the following two important elements:
$$\delta = \sum_{i=0}^n a_i\al_i\in  V
\mand
c = \sum_{i=0}^n a_i^\vee \al_i^\vee \in V^\ast
$$
where 
$a_0, \ldots, a_n$ and 
$a_0^\vee, \ldots, a_n^\vee$ are defined in \cite[Theorem 4.8]{Kac1984} and are determined by the Dynkin type.
Note that we always have $a_0^\vee=1$, regardless of the type, and $a_0 = 1$ except in type $A_{2n}^{(2)}$, in which case $a_0 = 2$.

\medskip

The space $V$ is equipped with a symmetric bilinear form $(.\mid.)$
which satisfies:
$$
\left( \al_i \mid \al_j  \right) = a_i^\vee a_i^{-1} a_{ij} \text{ for } i,j\neq 0
\ , \quad 
\left( \delta \mid \al_i  \right) = 0 \text{ for } i \neq 0
\mand
\left( \delta \mid \delta  \right) = 0.
$$
The space $V$ endowed with $(\cdot\mid\cdot)$ is a quadratic space 
with isotropic cone  $\mathbb R\delta$.
We extend $(.\mid .)$ to $	\mathfrak{h}^\ast$ by introducing an isotropic element 
$\La_0 \in  \mathfrak{h}^\ast\setminus V$ verifying
$
\left( \La_0 \mid \al_i  \right) = 0 \text{ for } i\neq 0$
and
$\left(\La_0 \mid \delta \right) = 1.$
In particular, $\left(\La_0 \mid \al_0\right) = a_0^{-1}$.  For $x\in \mathfrak{h}^\ast$, we denote $|x| = (x\mid x)^{1/2}$.
The bilinear form  enables us to express easily the coroots $\alpha^{\vee}$ for any $\alpha \in \Phi_0$ via the formula  \cite[Proposition 5.1 (d)]{Kac1984}
\begin{equation}\label{coroot in terms of root}
\al^\vee = 2\frac{(\al\mid .)}{(\al\mid\al)},
\end{equation}
and in fact, for $0\leq i\leq n$, we have
\begin{equation}\label{alpha_vs_alpha_check}
\al_i^\vee = \frac{a_i}{a_i^\vee}(\al_i\mid .).
\end{equation}

Consider now the element defined by
$$
\theta = \delta - a_0\alpha_0 =  \sum_{i=1}^n a_i \al_i \in V_0.
$$
Depending on the type, $\theta$ is either the highest root of $\Phi$ or the the highest short root of $\Phi$, and it satisfies the equality $(\theta \mid \theta) = 2a_0$, see \cite[Proposition 17.18]{Carter2005} for more details on $\theta$. 
Moreover, since $\theta \in V_0$ and $\theta^\vee \in V_0^*$, the sets $\{\delta, \alpha_1, \dots, \alpha_n\}$ and $\{c, \alpha_1^\vee, \dots, \alpha_n^\vee\}$  form a basis of $V$ and $V^*$ respectively.  In particular one has
$
V = V_0 \oplus \mathbb{R}\delta$ and  $V^* = V_0^* \oplus \mathbb{R}c.
$
Since $ \Lambda_0 \in \mathfrak{h}^\ast\setminus V$, it follows that
$$
\mathfrak{h}^\ast = V_0 \oplus \mathbb{R}\delta \oplus \mathbb{R}\Lambda_0.
$$

Let $\Omega$ be the $W_0$-orbit of $\theta^\vee$,
and $M$ be the preimage of the lattice $\Z\Omega$ under the isomorphism $V_0\to V_0^\ast$ induced
by the scalar product $(\cdot\mid\cdot)$, that is,
$$M=\{ \be \in V_0 \mid (\be\mid .)\in\Z\Omega\}.$$

For $x \in V_0$ we denote by $t_x : \mathfrak{h}^\ast \rightarrow \mathfrak{h}^\ast$ the linear map defined by:
\begin{equation}\label{translations M}
t_x(v) = v + \langle v,c \rangle x - \left( (v \mid x) + \displaystyle\frac{1}{2}|x|^2 \langle v, c \rangle\right)\delta.
\end{equation}

Let $T(M)$ be the subgroup of $GL(\mathfrak{h}^\ast)$ generated by the $t_{\beta}$ for $\beta \in M$. This group is called the group of translations of $M$ and it acts faithfully on $\mathfrak{h}^\ast$ by Formula (\ref{translations M}).  Moreover,  for all $x,y \in V_0$,  all $v \in \mathfrak{h}^\ast$ and all $w \in W$,  one has $t_xt_y(v) = t_{x+y}(v)$ and $wt_xw^{-1} = t_{w(x)}$.
Then one can express $W$ as follows (\cite[Proposition 17.22]{Carter2005}):
\begin{equation}\label{Weyl group Kac semi-direct}
W = T(M) \rtimes W_0.
\end{equation}

Hence, any element $w \in W$ decomposes uniquely as $w=t_{\beta}\overline{w}$ for $\beta \in M$ and $\overline{w} \in W_0$.
We recall the classification of affine Weyl groups that arise this way in \Cref{aff_weyl_gps}.
\begin{figure}[!h]
$$
\begin{array}{@{}l@{\hskip 10pt} @{}l@{\hskip 20pt} @{}l@{\hskip 20pt} @{}l@{\hskip 20pt} @{}l@{} @{}l@{}}
\cline{2-5}
&
\text{Type}
& 
\text{Alt. notation}
&
\text{Type of } 
W_0
&
\text{Lattice } M
& \phantom{\parbox[t]{8mm}{\multirow{8}{*}{\rotatebox[origin=c]{90}{untwisted types} 
$\left\{\rule{0mm}{3.2cm}\right.$}}}
\\
\cline{2-5}
\parbox[t]{8mm}{\multirow{7}{*}{\rotatebox[origin=c]{90}{untwisted types} 
$\left\{\rule{0mm}{2.8cm}\right.$}} 
&
A_n^{(1)} \ n\geq 1
&
\widetilde{A_n}
&
A_n
&
\Z\al_1+\cdots+\Z\al_{n-1}+\Z\al_n
&\\ &
B_n^{(1)}  \ n\geq 3
&
\widetilde{B_n}
&
B_n
&
\Z\al_1+\cdots+\Z\al_{n-1}+2\Z\al_n
&\\ &
C_{n}^{(1)}  \ n\geq 2
&
\widetilde{C_n}
&
C_n
&
2\Z\al_1+\cdots+2\Z\al_{n-1}+\Z\al_n
&\\ &
D_{n}^{(1)}   \ n\geq 4
&
\widetilde{D_n}  
&
D_n
&
\Z\al_1+\cdots+\Z\al_{n-1}+\Z\al_n
&\\ &
E_{n}^{(1)}\  n=6,7,8
&
\widetilde{E_n}
&
E_n
&
\Z\al_1+\cdots+\Z\al_n
&\\ &
F_{4}^{(1)}
&
\widetilde{F_4}
&
F_4
&
\Z\al_1+\Z\al_2+2\Z\al_3+2\Z\al_4
&\\ &
G_{2}^{(1)}
&
\widetilde{G_2}
&
G_2
&
\Z\al_1+3\Z\al_2
&\\
\parbox[t]{8mm}{\multirow{6}{*}{\rotatebox[origin=c]{90}{twisted types}
$\left\{\rule{0cm}{2.3cm}\right.$}}
&
A_2^{(2)}
&
\widetilde{A_1}'
&
A_1
&
\frac{1}{2}\Z\al_1
&\\ &
A_{2n-1}^{(2)}  \ n\geq 2
&
\widetilde{B_n}^{\mathrm{t}}
&
C_n
&
\Z\al_1+\cdots+\Z\al_{n-1}+\Z\al_n
&\\ &
D_{n+1}^{(2)}  \ n\geq 2
&
\widetilde{C_n}^{\mathrm{t}}
&
B_n
&
\Z\al_1+\cdots+\Z\al_{n-1}+\Z\al_n
&\\ &
A_{2n}^{(2)}  \ n\geq 2
&
\widetilde{C_n}'
&
C_n
&
\Z\al_1+\cdots+\Z\al_{n-1}+\frac{1}{2}\Z\al_n
&\\ &
E_{6}^{(2)}
&
\widetilde{F_4}^{\mathrm{t}}
&
F_4
&
\Z\al_1+\Z\al_2+\Z\al_3+\Z\al_4
&\\ &
D_{4}^{(3)}
&
\widetilde{G_2}^{\mathrm{t}}
&
G_2
&
\Z\al_1+\Z\al_2
&\\
\cline{2-5}
\end{array}
$$
\caption{The affine Weyl groups are of the form $T(M)\rtimes W_0$
where $W_0$ and $M$ are classified in the above table,
which can be recovered from
\cite[Formula 6.5.8]{Kac1984} 
or \cite[Proposition 17.23]{Carter2005}.
The terminology for the Dynkin type (first column)
is taken from Kac' book \cite[Chapter 4]{Kac1984}.
In the twisted types, the superscript coincides with the ratio of the squared lengths of the long and short roots (except in type $A_{2}^{(2)}$ where this ratio is $4$, and in type 
$A_{2n}^{(2)}$ where there are three different root lengths, with consecutive ratios $2$).
The alternative terminology (second column) is taken from \cite{Carter2005}.
It has the advantage of making the untwisted counterpart of each twisted type with a superscript ``$\mathrm{t}$'' appear.}
\label{aff_weyl_gps}
\end{figure}

\begin{Rem}\label{rem_aff_WG}\
\begin{enumerate}
\item In untwisted Dynkin types, we always have $M\simeq Q^\vee$ (the coroot lattice),
which can be checked from \Cref{aff_weyl_gps}.
This implies that the group $W$ is isomorphic (as a Coxeter group)
to the affine Weyl group $W_a$ constructed geometrically in \cite[Chapter VI, \S 2]{BOURB}.
We will detail this contruction in \Cref{subsec_shi} below.
\item The terminology ``group of translations'' for $T(M)$
is justified by \cite[Proposition 17.24]{Carter2005} or \cite[Formula 6.6.3]{Kac1984}.
\end{enumerate}
\end{Rem}

\subsection{Weight lattice and classic formulas}

Before proceeding  to the further study of $W$,
let us briefly recall some classical constructions related to 
the (finite) Weyl group $W_0$.
The \textit{fundamental weights} are defined as
the elements $\om_i\in V_0$, $1\leq i\leq n$ verifying
$$
\left\langle \om_i, \al_j^\vee \right\rangle
=\delta_{ij}.
$$
The set $P_0= \bigoplus_{i=1}^n \Z \om_i$ is the corresponding \textit{weight lattice}, 
and  $P_0^+= \bigoplus_{i=1}^n \N \om_i$ the set of integral dominant weights.
We set
$$ \overline{\rho} = \sum_{i=1}^n \om_i
=\frac{1}{2}\sum_{\al\in\Phi^+} \al.$$
In turn, there is a coweight lattice $P_0^\vee$
with basis consisting of the fundamental coweights $\om_i^\vee\in V_0^\ast$, $i=1, \ldots, n$,
defined by
$$
\left\langle \al_j, \om_i^\vee \right\rangle
=\delta_{ij}.
$$
We define similarly
$$\overline{\rho}^\vee = \sum_{i=1}^n \om_i^\vee= \frac{1}{2}\sum_{\al\in\Phi^+} \al^\vee.$$
We see that
$\left\langle\al_i, \overline{\rho}^\vee\right\rangle = 1$ 
for all $1\leq i\neq j\leq n$,
which implies in particular that, for all $\be\in Q$, 
\begin{equation}\label{inner_prod_height}
\left\langle \be, \overline{\rho}^\vee\right\rangle = \h(\be).
\end{equation}

\medskip

Let us come back to the general affine setting.

\medskip

Similarly to the finite setting,
we let $\La_i, 1\leq i\leq n$ be the remaining (affine) 
fundamental weights, defined by  $\La_i =a_i^\vee\La_0+\omega_i$
for $1\leq i\leq n$, so that
$\langle \La_i, \al_j^\vee\rangle = \delta_{ij}$.
We define  similarly $\La_i^\vee\in \mathfrak{{h}}$,  $1\leq i \leq n$, in particular
$\langle \al_j , \Lambda_i^\vee\rangle =\delta_{ij}$,
and we finally set $\La_0^\vee = a_0 (\La_0\mid \cdot) \in \mathfrak{{h}}$.
In particular, we have
$\langle \al_i , \La_0^\vee\rangle = \delta_{0i}$
for $0\leq i\leq n$ and $\langle \La_0, \La_0^\vee\rangle = 0$.
We can now  set
$$\rho = 
\La_0+\ldots+\La_{n} = h^{\vee}\La_0 + \overline{\rho}\in \mathfrak{h}^\ast
\mand 
\rho^\vee = 
\La_0^\vee+\ldots+\La_{n}^\vee = h \La_0^\vee + \overline{\rho}^\vee\in \mathfrak{{h}}
$$
where $h =\sum_{i=0}^n a_i$ is the 
\textit{Coxeter number} and 
$h^\vee =\sum_{i=0}^n a_i^\vee$ is the Coxeter number corresponding to the transpose of $A$.

The affine \textit{weight lattice} is
$\ds
 P = \bigoplus\limits_{i=1}^{n} \mathbb{Z}\omega_i \oplus \mathbb{Z}\delta \oplus \mathbb{Z}\Lambda_0,
$
and the lattice of affine \textit{dominant} weights is given by 
$
 P^+ 
= 
\bigoplus\limits_{i=0}^{n} \N\La_i \oplus \mathbb{Z}\delta.
$
Therefore, any dominant weight $\la\in P^+$ writes
$\la=\sum_{i=0}^n m_i\La_i +z\delta$ for some $m_i\in\N, z\in\Z$. Alternatively, $\la$ decomposes as
\begin{equation}\label{dom_weight_aff}
\la= \overline{\la} + \ell\La_0 +  z\delta
\end{equation}
where $\overline{\la}=m_1\om_1+\cdots+m_n\om_n\in \overline{P}^+$ and the \textit{level} of $\la$ 
(by using $a_0^\vee=1$) is the nonnegative integer
$$
\ell = m_0 + \sum_{i=1}^n a_i^\vee m_i
=\sum_{i=0}^n a_i^\vee m_i =  \langle \la, c \rangle.
$$

We record below the following important formulas of the bilinear form 
on $\mathfrak{h}^\ast$ and the pairing.
\begin{equation}\label{formulas}
\begin{array}{ll}
 \left\{
 \begin{array}{ll}
 \left( \alpha_i \mid \alpha_j  \right) = \displaystyle\frac{a_i^{\vee}}{a_i}a_{ij}&
0\leq i, j\leq n
\\
 \left( \alpha_i \mid \om_j  \right) = \displaystyle\frac{a_i^{\vee}}{a_i}\delta_{ij}&
0\leq i, j\leq n
\\
\left( \alpha_i \mid \delta  \right) = 0 & 
0\leq i \leq n
\\
\left( \alpha_i \mid \Lambda_0  \right) = 0 & 
1\leq i\leq n
\\
\left( \alpha_0 \mid \Lambda_0  \right) =  a_0^{-1} \\
\left( \delta \mid \delta  \right) = 0 \\
\left( \Lambda_0 \mid \Lambda_0  \right) = 0 \\
\left( \Lambda_0 \mid \delta  \right) = 1\\
 \end{array}
 \right.
 &
 \quad
 \left\{
 \begin{array}{ll}
\langle  \alpha_j ,\alpha_i^{\vee}\rangle =  a_{ij} &  
0\leq i, j\leq n
\\
\langle   \delta, \alpha_i^{\vee} \rangle = 
0 &  
0\leq i \leq n
\\
\langle \Lambda_0, \alpha_0^{\vee} \rangle = 1 \\
\langle \Lambda_0, \alpha_i^{\vee} \rangle = 0  & 
1\leq i\leq n
\\
\langle \alpha_i, \rho^{\vee} \rangle = 1  & 
0\leq i\leq n
\\
\langle \delta, \rho^{\vee} \rangle = h.
& 
\\
 \end{array}
 \right.
 \end{array}
\end{equation}

\subsection{Untwisted affine Weyl groups and Shi coefficients}
\label{subsec_shi}

We already mentioned in \Cref{rem_aff_WG}(1)
that the above construction recovers the construction of
the affine Weyl groups of untwisted Dynkin types.
Let us recall the direct geometric construction of $W_a$ in this case.
Let $k \in \mathbb{Z}$ and $\al \in \Phi_0$. We define the affine reflection $s_{\al,k} \in V_0 \rtimes \mathrm{GL}(V_0)$ (the affine group of $V_0$) by
$$
s_{\al,k}(x)=x- \left( \left( \al\mid x \right)-k\right)\frac{2 \al}{\left( \al\mid\al \right) }.
$$
Note that for all $\al\in\Phi_0$, $s_{\al,0} = s_\al$.
The group generated by all the affine reflections $s_{\al,k}$ with $\al \in \Phi_0$ and $k \in \mathbb{Z}$ is called the \textit{affine Weyl group} associated to $\Phi_0$ and is denoted by $W_a$.
The group $W_a$ is a Coxeter group and we also denote by $\ell$ its  length function. 

\medskip

We set the hyperplanes
$$H_{\al,k} = \{ x \in V_0\mid ( \al \mid  x ) = k\},$$ 
and the strips
$$H_{\al,k}^1  = \{x \in V_0\mid k < ( \al \mid  x ) < k+1 \}.$$
An \textit{alcove} of $V_0$ is  a connected component of
$$
 V_0 ~\backslash \bigcup\limits_{\begin{subarray}{c}
 ~ ~\al \in \Phi_0^{+} \\ 
  k \in \mathbb{Z}
\end{subarray}}
H_{\al,k}.
$$
We denote $A_e$ the \textit{fundamental alcove},
defined by $A_e = \bigcap_{\al \in \Phi_0^+} H_{\al,0}^1$. 
The fundamental chamber of $W_a$ is defined by
$$
C_0 = \{x \in V_0 \mid  \left( x\mid\beta\right) \geq 0 ~ \forall \beta \in \Delta_0\}.
$$
The group $W_a$ acts regularly on the set of alcoves. 
Therefore we have a bijective correspondence between 
$W_a$ and the set of alcoves, given by $w \mapsto A_w$ where $A_w = wA_e$.
Moreover, every alcove of $V_0$ can be written as an intersection of $|\Phi_0^+|$ strips, that is, for all $w\in W_a$, we have
$$
A_w = \bigcap\limits_{\al \in \Phi_0^+}H_{\al, k(w,\al)}^1.
$$
where $k(w,\al)\in\Z$ is called the \textit{Shi coefficient} of $w$ in position $\al.$ 
This integer indicates the number of hyperplanes $H_{\al,m}$ with $m \in \mathbb{Z}$ (also called the $\alpha$-hyperplanes) 
between the alcove $A_w$ and the fundamental alcove $A_e$. 
More precisely,  if $k(w,\al) \geq 0$ then the $\alpha$-hyperplanes between $A_e$ and $A_w$ are
$H_{\al,1}, H_{\al,2}, \dots, H_{\al,k(w,\al)}$. If  $k(w,\al) < 0$ then the $\alpha$-hyperplanes between $A_w$ and $A_e$ are
$H_{\al,0},H_{\al,-1}, \dots, H_{\al,k(w,\al)+1}$.

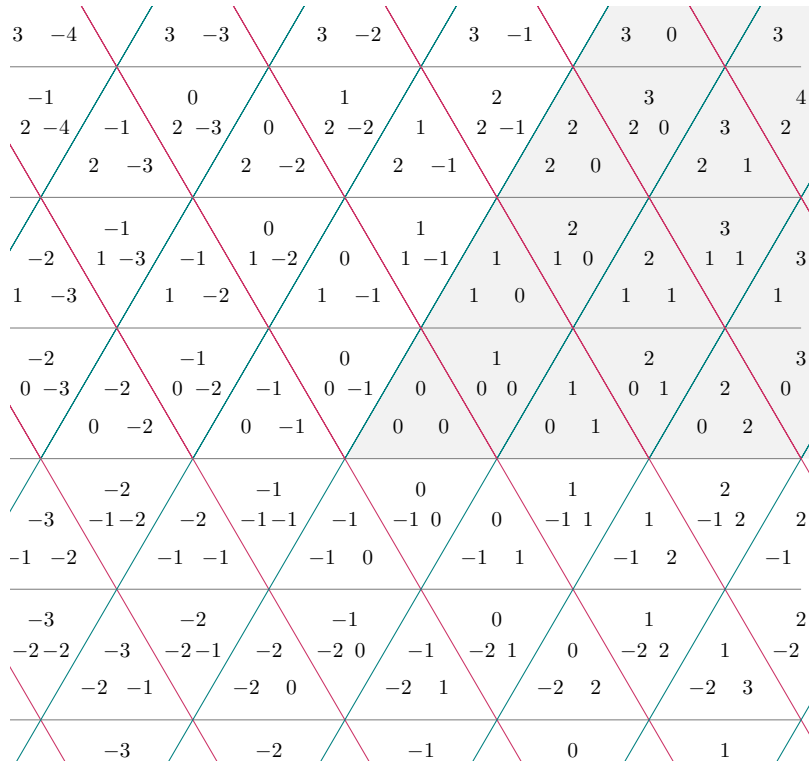
\begin{figure}[h!]
\centering
\begin{tikzpicture}[scale=2]
\clip (-2.2,-2) rectangle (3.05,3);
\path[fill=gray!10] (0,0) -- +(0:10) -- ++(60:10);
\foreach \i in {-4,...,4}
{
\draw[gray] (-3,{\i*sqrt(3)/2}) -- (3,{\i*sqrt(3)/2});
\draw[teal] (\i,0) -- +(60:4) -- ++(60:-4);
\draw[purple!80] (\i,0) -- +(120:4) -- ++(120:-4);
}

\foreach \j in {-4, ..., 4}
{
\foreach \i in {-4,...,4}
{
\node[anchor = mid, scale=0.8] at ( {\i+0.35+(\j*0.5)}, {0.2+(\j*sqrt(3)/2)} ) {$\j$};
\node[anchor = mid, scale=0.8] at ( {\i+0.65+(\j*0.5)}, {0.2+(\j*sqrt(3)/2)}) {$\i$};
\node[anchor= mid, scale=0.8] at ( {\i+0.5-(\j*0.5)}, {0.45+(\j*sqrt(3)/2)}) {$\i$};
\node[anchor= mid,scale=0.8] at ( {\i-0.1+(\j*0.5)}, {0.45+(\j*sqrt(3)/2)}) {$\j$};
\node[anchor= mid,scale=0.8] at ( {\i+1.1+(\j*0.5)}, {0.45+(\j*sqrt(3)/2)}) {$\i$};
\node[anchor= mid, scale=0.8] at ( {\i-(\j*0.5)}, {0.65+(\j*sqrt(3)/2)}) {$\i$};
}
}
\end{tikzpicture}
\caption{Alcoves in $W_a=W(A_2^{(1)})$.
In each alcove, we have written the Shi vector of
the corresponding Weyl group element,
where the bottom left (respectively bottom right, respectively top) position corresponds to the root $e_{12}$
(respectively $e_{23}$, respectively $e_{13}$).
The identity element corresponds to the alcove with $3$ zeros
and the elements of $W_0$ correspond to the alcoves with only zeros or $-1$'s, forming a hexagon.
The shaded region, consisting of only nonnegative integers, is the fundamental chamber $C_0$.}
\label{alcove A2}
\end{figure}

The vector $(k(w,\al))_{\al \in \Phi_0^+}$ is called the \textit{Shi vector} of $w$.  Shi vectors are usually arranged in a ``pyramidal'' shape as illustrated in \Cref{alcove A2}, \Cref{Shi_A}, 
\Cref{ShiB4} and \Cref{ShiC4}.
In \cite{JYS1}, Shi gave the following characterisation of the possible integers vectors $(k_\al)_{\al \in \Phi_0^+}$ 
that arise as Shi vectors of elements in $W_a$ (see Figure \ref{alcove A2} for an example).

\begin{Th}[{\cite[Theorem 5.2]{JYS1}}]\label{thJYS1}
Let $A = \bigcap\limits_{\al \in \Phi_0^+} H^1_{\al,k_{\al}}$ with $k_{\al} \in \mathbb{Z}$. 
Then $A$ is an alcove if and only if for all $\al$, $\beta \in \Phi_0^+$ satisfying  $\al + \beta \in \Phi_0^+$, we have the following inequality
\begin{equation*}\label{Shi ineq}
k_{\al} + k_{\beta} \leq k_{\al+\beta}  \leq k_{\al} + k_{\beta} + 1
\end{equation*}
\end{Th}

The following proposition is essential in view of proving \Cref{Proposition Fondamentale}.

\begin{Prop}[{\cite[Proposition 3.2]{NC1}}]\label{formula_Nathan}
Let $w \in W_a$ and let $t\in T$ be a reflection.
For all $\al \in \Phi_0$, we have  
$$k(tw,\al) = k(w,t(\al)) + k(t,\al).$$
\end{Prop}

We will describe in \Cref{lemma Shi coeff dans J} below the inversion set of affine Weyl group elements in terms of their Shi coefficients. 
This will be crucial for proving \Cref{Proposition Fondamentale} and, in turn, \Cref{thm_surj}.

\begin{Lem}[{\cite[Lemma 3.1]{JYS1}}]\label{Shi inversion}
Let $w \in W_0$ and $\al \in \Phi^+$.  Then $\al \in N(w)$ if and only if $k(w,\al) = -1$.
\end{Lem}

\medskip

Taking inversion sets is compatible with restricting to reflection subgroups as is shown in the following lemma\footnote{
Point (2) is well-known but
we find it useful to give the proof again.
}.

\bigskip ~

\begin{Lem}\label{lemma Shi coeff dans J}
Let $W$ be a finite Weyl group with simple system $S$. Let $w\in W$.
\begin{enumerate}
\item Let $A\subseteq T$ and consider the reflection subgroup $W_A\leq W$.
Then $N_{A}(w_A) = \{\al \in \Phi_A^+\mid k(w,\al) = -1\}$.  
\item Let $I\subseteq S$ and consider the standard parabolic subgroup $W_I\leq W$.
Then $ N_{I}(w_I) = N(w_I)$.   
\end{enumerate}
\end{Lem}

\begin{proof}\
\begin{enumerate}
\item By \Cref{propDy}, we have
$$N_A(w_A) = N(w) \cap \Phi_A = \{\al\in\Phi^+ \mid k(w, \al) =-1\} \cap \Phi_A 
= \{\al\in \Phi_A^+ \mid k(w, \al) =-1\}.$$
\item 
Let $w_I=s_{i_1}\dots s_{i_k}$ be a reduced expression of $w_I\in W_I$. 
Since $\ell_I(w_I) = \ell(w_I)$, this expression is also reduced in $W$. 
By \Cref{Decomposition_N} applied in both $W_I$ and in $W$, we get
$N_I(w_I) = \{\al_{i_1},  s_{i_1}(\al_{i_2}),s_{i_1}s_{i_2}(\al_{i_3}),\dots,s_{i_1} \dots s_{i_{k-1}}(\al_{i_k})\} = N(w_I)$.
\end{enumerate}
\end{proof}


\section{Atomic length in finite Weyl groups}\label{sec_AL}

We are ready to define the central notion of this paper,
namely the \textit{atomic length}
of a Weyl group element.
In this section, we only consider Weyl groups of finite type,
except in \Cref{def_lambda_AL},
which is given for Weyl groups of both finite and affine type Kac-Moody algebras (as introduced in \Cref{sec_affine_WG}).
Indeed, we will first investigate the properties of the atomic length in finite types, in both \Cref{sec_AL}
and \Cref{sec_surj_AL}.
In order to avoid cumbersome notation, we will simply 
use the notation $W$, $\rho$, ... for the finite Weyl group, 
the sum of the fundamental weights, and so on (instead of 
$W_0$, $\overline{\rho}$, ... as in \Cref{sec_affine_WG}).
The study of the affine atomic length will be  delayed to \Cref{sec_cryst} and \Cref{sec_affine_AL}.

\subsection{Inversion sets and atomic length}
\label{subsec_AL_inv}

We use the notation of the previous sections.
In particular, $\Phi^+$ denotes
the chosen set of positive roots,
and $N(w)$ denotes the set of inversions of $w\in W$.

\medskip

\begin{Def}\label{atomic_length}
Let $w\in W$. The \textit{atomic length} of $w$
is the nonnegative integer 
$$\sL(w) = \sum_{\al\in N(w)} \h(\al).$$
\end{Def}

\medskip

\begin{Rem}\label{rem_def_AL}\
\begin{enumerate}
 \item Recall from \Cref{sec_Cox}
that for all $w\in W$, we have $|N(w)|=\ell(w)$,  the length of $w$.
Therefore, the atomic length can be seen as a variant of the usual length function, where
each inversion $\al$ is counted with multiplicity $\h(\al)$,
that is, $\al$ is fully decomposed as a sum of simple roots,
hence the terminology ``atomic'' length. 
\item Assume that $W=W(A_n)$.
Recall the definition of $\mathsf{invsum}$ from \Cref{subsec_entropy_inv}.
We have already seen in \Cref{rem_inv_type_A} that $\h(e_{ij})=j-i$ for all $1\leq i<j\leq n+1$.
Therefore, for all $w\in W$, we have $\sL(w)=\mathsf{invsum}(w)$.
\end{enumerate}
\end{Rem}

\begin{Exa}\label{exa_atomic_length_S3} 
Let $W=W(A_2)$.
Denote $\al_i= e_{i,i+1}\in\De$, $i=1, 2$ the simple roots and $s_1, s_2\in W$ the corresponding simple reflections.
The values of both the usual and atomic length are recorded in the following table.
$$
\begin{array}{@{}l@{\hskip 20pt} @{}l@{\hskip 20pt} @{}l@{\hskip 20pt} @{}l@{}}
\toprule
w
& 
N(w)
&
\ell(w)
&
\sL(w)
\\
\midrule
1
&
\emptyset
&
0
&
0
\\
s_1
&
\{\al_1\}
&
1
&
1
\\
s_2
&
\{\al_2\}
&
1
&
1
\\
s_1s_2
&
\{\al_1, \al_1+\al_2\}
&
2
&
3
\\
s_2s_1
&
\{\al_2, \al_1+\al_2\}
&
2
&
3
\\
s_1s_2s_1
&
\{\al_1, \al_2, \al_1+\al_2\}
&
3
&
4
\\
\bottomrule
\end{array}
$$
\end{Exa}

The following lemma will enable us to extend the definition
of the atomic length.
Recall that we have defined $\rho$ in \Cref{sec_affine_WG}
to be the half-sum of the positive roots.

\begin{Lem}\label{lem_rho_inv}
For all $w\in W$, we have 
$$\rho- w(\rho) = \sum_{\al \in N(w)} \al.$$
\end{Lem}

\begin{proof} Fix $w\in W$.
First, we have
\begin{equation*}
\rho   = \ds \frac{1}{2}\sum_{\al\in\Phi^+} \al =
\ds
\frac{1}{2}\sum_{\substack{\al\in \Phi^+ \\ \al\in N(w)}}\al
+
\frac{1}{2}\sum_{\substack{\al\in \Phi^+ \\ \al\notin N(w)}}\al.
\end{equation*}
Second, we have
$$
\begin{array}{rcll}
w(\rho) 
& = &
\ds
\frac{1}{2}\sum_{\al\in\Phi^+} w(\al)
&
\\
& = &
\ds 
\frac{1}{2}\sum_{\substack{\al\in\Phi^+ \\ w(\al)\in \Phi^+}} w(\al)
+
\frac{1}{2}\sum_{\substack{\al\in\Phi^+ \\ w(\al)\notin \Phi^+}} w(\al)
&
\\
& = &
\ds
\frac{1}{2}\sum_{\substack{w^{-1}(\al)\in\Phi^+ \\ \al\in \Phi^+}} \al
+ 
\frac{1}{2}\sum_{\substack{w^{-1}(\al)\in\Phi^+ \\ \al\notin \Phi^+}} \al
&
\text{ by substituting $\al\leftarrow  w^{-1}(\al)$ in both sums}
\\
& = &
\ds
\frac{1}{2}\sum_{\substack{w^{-1}(\al)\in\Phi^+ \\ \al\in \Phi^+}} \al
- 
\frac{1}{2}\sum_{\substack{w^{-1}(\al)\notin\Phi^+ \\ \al\in \Phi^+}} \al
&
\text{ by substituting $\al\leftarrow  -\al$ in the second sum}
\\
& = &
\ds
\frac{1}{2}\sum_{\substack{\al\in \Phi^+\\ \al \notin N(w)}} \al
- 
\frac{1}{2}\sum_{\substack{\al\in \Phi^+\\ \al \in N(w)}} \al.
&
\end{array}
$$
Finally, taking the difference yields exactly the desired identity.
\end{proof}

For the following definition,
we allow $W$ to be a Weyl group of either a finite or affine
type Kac-Moody algebra
(see \Cref{sec_affine_WG}).
Moreover, let $\la\in P^+$, that is, 
$\la$ is a dominant weight.
We define a statistic $\sL_\la$ on $W$, which
will specialise at $\sL$ for $\la=\rho$ in finite types.
Recall the element $\rho^\vee$ introduced in \Cref{sec_affine_WG}, equal to the sum of the fundamental coweights.

\begin{Def}\label{def_lambda_AL}
Let $\la\in P^+$. The \textit{$\la$-atomic length} of $w$ is the nonnegative integer
$$\sL_\la(w) = \left\langle \la-w(\la), \rho^\vee \right\rangle.$$
\end{Def}

\begin{Rem}\label{AL_integer}\
\begin{enumerate}
 \item 
Let us quickly explain why $\sL_\la(w)$ is a nonnegative integer.
If $\mathfrak{g}$ denotes the corresponding 
Kac-Moody algebra, we can construct the irreducible $\mathfrak{g}$-module
$V(\la)$ with \textit{highest weight} $\la\in P^+$, see for instance \cite[Section 2.3]{HongKang2002}.
This module is a weight module, that is, it decomposes as the direct sum of its $\mu$-weight spaces
$V(\la)_\mu = \{ v \in V(\la) \mid xv = \langle \mu, x\rangle v \text{ for all } x\in V^\ast \},$
and where $\mu$ runs over $P$.
Those elements $\mu\in P$ such that $V(\la)_\mu \neq 0$ are called the \textit{weights} of $V(\la)$,
and we know that the Weyl group $W$ acts on the set of weights of $V(\la)$.
Therefore, for all $w\in W$, $w(\la)$ is a weight of $V(\la)$.
Since $\la$ is the highest weight of $V(\la)$, we can write
$\la-w(\la) = \sum_{i=1}^n a_i \al_i$ with $a_i\in\N$.
By \Cref{inner_prod_height}, 
$\sL_\la(w) =\sum_{i\in I} a_i\in\N$.
\item
One could decide to define the atomic length for $\la\in P$ using the same formula,
but the range of the map $\sL_\la$ would no longer be contained in $\N$.
For instance, take $\la=\om_2$ in type $B_2^{(1)}$ (it is a level zero weight but is not dominant).
Then the element $w = s_1 s_0 s_1 s_2 s_1 s_0 s_2$ verifies $\sL_\la(w) = -1$.
Note that this fails even though 
$\om_2$ is particularly nice in type $B_2$ as it is \textit{minuscule}, see \Cref{minusc}.
\end{enumerate}

\end{Rem}

The terminology of \Cref{def_lambda_AL} is justified by  
the following proposition, which we will use extensively in the rest of the article.

\begin{Prop}\label{prop_L_rho}
Assume that $W$ is of finite type.
Then $\sL = \sL_\rho$.
\end{Prop}

\begin{proof}
For all $w\in W$, we have
$$
\begin{array}{rcll}
\sL_\rho(w)
&
=
&
\ds 
\left\langle \rho- w(\rho),\rho^\vee\right\rangle
&
\\
&
=
&
\ds 
\left\langle \sum_{\al\in N(w)} \al,\rho^\vee\right\rangle
&
\text{ by \Cref{lem_rho_inv}}
\\
&
=
&
\ds 
\sum_{\al\in N(w)}\left\langle \al,\rho^\vee\right\rangle
&
\\ 
&
=
&
\ds 
\sum_{\al\in N(w)}\h(\al)
&
\text{ by (\ref{inner_prod_height})}
\\
& 
=
&
\sL(w).
&
\end{array}
$$
\end{proof}

\subsection{First properties of the atomic length}
\label{subsec_AL_prop}

We still assume that
$W$ is a finite Weyl group.
A first consequence of the characterisation of \Cref{prop_L_rho} is a remarkable symmetry
property for the atomic length in simply-laced Dynkin types.

\begin{Th}[Symmetry]\label{symmetry}
Assume that $W$ is of simply-laced finite type.
Then for all $w\in W$, $\sL(w)=\sL(w^{-1})$.
\end{Th}

\begin{proof} 
Since we are in type $A$, $D$ or $E$, all roots $\al$ verify $(\al\mid\al)=2$,
which implies that
$\ds\rho^\vee = \frac{1}{2}\sum_{\al\in\Phi^+} (\al\mid .).$
Therefore, for all $x\in V$,
$$\langle x,\rho^\vee\rangle 
= 
\left\langle x, \frac{1}{2}\sum_{\al\Phi^+} (\al\mid .) \right\rangle
=
\frac{1}{2}\sum_{\al\Phi^+} (\al\mid x)
=
\left( \frac{1}{2}\sum_{\al\Phi^+} \al \;\middle|\; x \right)
=
\left( \rho  \;\middle|\; x \right).
$$
Let $w\in W$.
We have
$$
\begin{array}{llll}
\sL(w) 
&
=
& 
\left\langle \rho - w(\rho) , \rho^\vee \right\rangle
&
\text{ by \Cref{prop_L_rho}} 
\\
&
=
&
\left( \rho \mid \rho - w(\rho) \right)
&
\text{ by the previous computation with $x=\rho-w(\rho)$}
\\
&
=
&
\left( \rho - w^{-1}(\rho) \mid \rho \right)
&
\text{ because $w$ is an isometry}
\\
&
=
&
\sL(w^{-1})
&
\text{ by \Cref{prop_L_rho}}.
\end{array}
$$
\end{proof}

\begin{Exa}\
\begin{enumerate}
\item 
Let $W$ be the Weyl group of type $D_4$, with simple roots $\al_i =  e_{i,i+1}$, $1\leq i\leq 3$ and $\al_4= e^{34}$,
and let $s_i$, $1\leq i\leq 4$ be the corresponding simple roots.
Take $w=s_4s_1s_2s_3s_1s_2s_1$ (reduced expression).
Then $$N(w)=\left\{ 
\al_1, 
\al_2, 
\al_3,
\al_1+\al_2,
\al_2+\al_3,
\al_1+\al_2+\al_3, 
\al_1+2\al_2+\al_3+\al_4
\right\}
$$
We have $w^{-1}=s_2s_3s_1s_2s_4s_3s_1$, so that
$$N(w^{-1})=\left\{ 
\al_1, 
\al_3, 
\al_4,
\al_2+\al_4,
\al_1+\al_2+\al_4, 
\al_2+\al_3+\al_4,
\al_1+\al_2+\al_3+\al_4
\right\},$$
and we see that $\sL(w)=\sL(w^{-1})=15.$
\item Let us give a counterexample in a non simply-laced Dynkin type.
Let $W$ be the Weyl group of type $G_2$, with simple roots $\al_1= e_1- e_2$
and $\al_2=-2 e_1+ e_2+ e_3$.
Let $w=s_2s_1$ so that $w^{-1}=s_1s_2$.
We have $N(w)=\{\al_1, \al_1+\al_2\}$ so $\sL(w) = 3$, 
but $N(w^{-1})=\{\al_1, 3\al_1+\al_2\}$ so $\sL(w^{-1}) = 5$.
\end{enumerate}
\end{Exa}

Obviously, for the usual length function,
regardless of the type,
we have $\ell(w) = \ell(w^{-1})$ for all $w\in W$.
The symmetry of the atomic length proved in \Cref{symmetry} can thus be seen as an analogue of this property of
the usual length.
In fact, we will now see that further properties of the usual
length have analogues for the atomic length.

\medskip

Recall that we denoted $w_0$
the longest element of $W$.
The map $\al_i \mapsto -w_0(\al_i)$
is an involution of $\De$, 
and we can write $\al_{\zeta(i)} = -w_0(\al_i)$,
where $\zeta : \{1,\ldots, n\} \to \{1,\ldots, n\}$ 
is the corresponding  Dynkin diagram automorphism.
The following result is an analogue of the well-known 
fact that $\ell(w_0 w) = \ell(w_0) - \ell(w)$.

\begin{Th}\label{AL_symmetry_w0}
Let $\la\in P^+$.
For all $w\in W$, we have $\sL_\la(w_0 w) = \sL_\la(w_0)-\sL_\la(w)$.
\end{Th}

\begin{proof}
Let us show that the quantities $\sL_\la(w)$ and $\sL_\la(w_0) - \sL_\la(w_0w)$ are equal.
On the one hand, we have 
$$\sL_\la(w) = \langle \la-w(\la), \rho^\vee\rangle
=\sum_{i\in I} a_i,$$
where $\la-w(\la) = \sum_{i\in I} a_i\al_i$ as we have already seen.
On the other hand,
$$
\begin{array}{rcl}
\ds\sL_\la(w_0) - \sL_\la(w_0w)
& = &
\ds
\langle \la -w_0(\la) - \la + w_0w(\la), \rho^\vee\rangle
\\
& = &
\ds
\langle -w_0( \la -  w(\la)), \rho^\vee\rangle
\\
&
=
&
\ds
\left\langle -w_0\left(\sum_{i\in I} a_i\al_i\right) , \rho^\vee\right\rangle
\\
&
=
&
\ds
\left\langle \sum_{i\in I} a_i (-w_0(\al_i)), \rho^\vee\right\rangle
\\
&
=
&
\ds
\left\langle \sum_{i\in I} a_i \al_{\zeta(i)}, \rho^\vee\right\rangle
\\
&
=
&
\ds\sum_{i\in I} a_i 
\left\langle \al_{\zeta(i}), \rho^\vee\right\rangle
\\
&
=
&
\ds
\sum_{i\in I} a_i.
\end{array}
$$
\end{proof}

Finally, we are able to give a short representation-theoretic proof of the fact
that the longest element $w_0$ is also the ``atomic-longest'' element of $W$.

\begin{Th}\label{atomic_longest_element}
Let $\la\in P^+$. Then $\sL_\la(w)\leq \sL_\la(w_0)$ for all $w\in W$.
\end{Th}

\begin{proof}
As in \Cref{AL_integer}, consider the highest weight module $V(\la)$ for the corresponding finite-dimensional simple Lie algebra.
Then it is well-known that $w_0(\la)$ is the lowest weight of $V(\la)$.
In other words, for all weight $\mu$ of $V(\la)$, we can write
$\mu- w_0(\la) = \sum_{i=1}^n a_i \al_i$ with $a_i\in\N$.
In particular, since the Weyl group $W$ acts on the set of weights of $V(\la)$, 
we can write $w(\la)- w_0(\la) = \sum_{i=1}^n a_i \al_i$ with $a_i\in\N$ for all $w\in W$.
Now, we have
$$
\begin{array}{rcl}
\sL_\la(w)-\sL_\la(w_0) 
&
=
&
\langle \la-w_0(\la), \rho^\vee\rangle - \langle \la-w(\la), \rho^\vee\rangle 
\\
&
=
&
\ds
\left\langle w(\la)-w(\la), \rho^\vee\right\rangle 
\\
&
=
&
\ds\left\langle \sum_{i=1}^n a_i \al_i, \rho^\vee\right\rangle 
\\
&
=
&
\sum_{i=1}^n a_i \geq  0.
\end{array}
$$
\end{proof}

\subsection{Generalised inversion sets and $\la$-atomic length}
\label{subsec_AL_inv_generalised}

To be complete, we will now give a formula for the
$\la$-atomic length resembling Definition \ref{atomic_length}.
This requires to introduce a notion of \textit{$\la$-inversion set}.
For this, we use the characterisation of the inversion set $N(w)$ given in \Cref{Decomposition_N}.
In the rest of this section, let us denote $\Red(w)$ the set of all reduced expressions of a given $w\in W$.

\begin{Not}\label{Def_w_jk}
Let $w \in W$, let $s_j \in S$ and let $\underline{w} = s_{i_1}\dots s_{b_1}s_j\dots s_{b_2}s_j\dots s_{b_{p_j}}s_j\dots s_{i_r} \in \Red(w)$,
where $p_j$ is the number of times that the letter $s_j$ appears in $\underline{w}$.  
We set $\underline{w}_{j,k}= s_{i_1}\dots s_{b_{1}}s_j\dots s_{b_2}s_j\dots s_{b_{k}}$. 
For example $\underline{w}_{j,1} = s_{i_1}\dots s_{b_1}$ and $\underline{w}_{j,2} = s_{i_1}\dots s_{b_1}s_j\dots s_{b_2}$.  
We also define $\sigma_j(\underline{w}) = \sum\limits_{k=1}^{p_j} \underline{w}_{j,k}(\al_j)$.
\end{Not}

\begin{Rem}\label{rem_inv_set_formula}
Recal that we gave in \Cref{Decomposition_N} the formula, for $w\in W$,
$$N(w) = \{\al_{i_1},  s_{i_1}(\al_{i_2}),s_{i_1}s_{i_2}(\al_{i_3}),\dots,s_{i_1} \dots s_{i_{r-1}}(\al_{i_r})\}.$$
Using the definition of $\underline{w}_{j,k}$ in \Cref{Def_w_jk}, we see that to compute $N(w)$,
it suffices to choose a reduced expression $\underline{w}\in\Red(w)$, then we have
$N(w) = \bigsqcup_{j=1}^n  \{\underline{w}_{j,k}(\al_j)\mid k=1,\dots,p_j\}$.
We see in particular that the inversion set does not depend on the choice of a reduced expression.
\end{Rem}

\begin{Exa}\label{example_N}
Let $w \in A_4 $ with reduced expression given by $\underline{w} = s_1s_2s_1s_3s_4s_3$.  Then we have 
\begin{align*}
N(w) =&  \{\underline{w}_{1,k}(\al_1)\mid k=1,2\} \sqcup \{\underline{w}_{2,1}(\al_2)\} \sqcup \{\underline{w}_{3,k}(\al_3)\mid k=1,2\}\sqcup \{\underline{w}_{4,1}(\al_4)\} \\
= & \{\al_1,s_1s_2(\al_1)\} \sqcup \{s_1(\al_2)\} \sqcup \{s_1s_2s_1(\al_3), s_1s_2s_1s_3s_4(\al_3)\} \sqcup \{s_1s_2s_1s_3(\al_4)\} \\
=& \{e_{12}, e_{23}\} \sqcup \{e_{13}\} \sqcup \{e_{14}, e_{45}\} \sqcup \{e_{15}\} \\
 = & \{e_{12}, e_{23}, e_{13}, e_{14}, e_{45}, e_{15}\}
\end{align*}
\end{Exa}

\begin{Def}\label{Def lambda inverse}
Let $w\in W$ and $\underline{w}\in\Red(w)$.  
Let $\lambda \in P$ be a weight with $\lambda = \sum\limits_{j=1}^nm_j\omega_j$. 
The set
$$
N_{\lambda}(\underline{w}) = \bigsqcup\limits_{j=1}^n  \{m_j\underline{w}_{j,k}(\al_j)\mid k=1,\dots,p_j\}.
$$
is called a \textit{$\lambda$-inversion set} of $w$.
\end{Def}

\begin{Rem}\
\begin{enumerate}
\item By definition, the notion of $\lambda$-inversion set depends on the choice of a reduced expression
(which was not the case for the usual inversion set $N(w)$, see \Cref{rem_inv_set_formula}).
We illustrate in \Cref{exa_pseudo_inv} how different reduced expressions yield different $\lambda$-inversion sets in general.
\item However, in the case $\lambda = \rho$, we have $N_\rho(\underline{w}) = N(w)$ for all $\underline{w}\in\Red(w)$,
so the $\rho$-inversion set of $w$ does not depend on the choice of a reduced expression.
\end{enumerate}
\end{Rem}

\begin{Exa}\label{exa_pseudo_inv}
We continue Example \ref{example_N}.  Let $\lambda = \sum_{j=1}^5m_j\omega_j$. Then by definition we have
\begin{align*}
N_{\lambda}(\underline{w}) =&  \{m_1\underline{w}_{1,k}(\al_1)\mid k=1,2\} \sqcup \{m_2\underline{w}_{2,1}(\al_2)\} \sqcup \{m_3\underline{w}_{3,k}(\al_3)\mid k=1,2\}\sqcup \{m_4\underline{w}_{4,1}(\al_4)\} \\
=& \{m_1e_{12}, m_1e_{23}\} \sqcup \{m_2e_{13}\} \sqcup \{m_3e_{14}, m_3e_{45}\} \sqcup \{m_4e_{15}\} \\
 = & \{m_1e_{12}, m_1e_{23}, m_2e_{13}, m_3e_{14}, m_3e_{45}, m_4e_{15}\}.
\end{align*}
Let us take now another reduced expression of $w$, 
say $\underline{w}' =s_2s_1s_4s_2s_3s_4$. The same kind of computations yields
\begin{align*}
N_{\lambda}(\underline{w}') = \{m_1e_{13}, m_2e_{23}, m_2e_{12}, m_3e_{15}, m_4e_{14}, m_4e_{45}\}.
\end{align*}
We can easily find coefficients $m_j$ such that $N_{\lambda}(\underline{w}') \neq N_{\lambda}(\underline{w})$. 
\end{Exa}

The following lemma is key, as the dependency on the choice of a reduced expression vanishes.

\begin{Lem}\label{Lemma omega_j}
Let $w\in W$ and $\underline{w}\in\Red(w)$. Fix $1\leq j \leq n$.
We have $\omega_j - w(\omega_j) = \sigma_j(\underline{w})$. In particular $\sigma_j(\underline{w})$ 
does not depend on the choice of a reduced expression of $w$.
\end{Lem}

\begin{proof}
We know that for any generator $s_i \in S$ and for any fundamental weight $\omega_j$ one has $s_i(\omega_j) = \omega_j$ if $j\neq i$ and $s_i(\omega_i) = \omega_i - \al_i$.  Therefore it follows that
 \begin{align*}
 \underline{w}(\omega_j) & = s_{i_1}\dots s_{b_1}s_j\dots s_{b_2}s_j\dots s_{b_{p_j}}s_j\dots s_{i_r}(\omega_j)  \\
 					 & = s_{i_1}\dots s_{b_1}s_j\dots s_{b_2}s_j\dots s_{b_{p_j}}s_j(\omega_j) \\
 					 & = s_{i_1}\dots s_{b_1}s_j\dots s_{b_2}s_j\dots s_{b_{p_j}}(\omega_j - \al_j) \\
 					 & = s_{i_1}\dots s_{b_1}s_j\dots s_{b_2}s_j\dots s_{b_{p_j}}(\omega_j) - \underline{w}_{j,p_j}(\al_j)\\
 					 & = \cdots \\
 					 & = \omega_j - \sum\limits_{k=1}^{p_j} \underline{w}_{j,k}(\al_j) \\
 					 & = \omega_j - \sigma_j(\underline{w}).
 \end{align*} 
 However, it is obvious that $\underline{w}(\omega_j) = w(\omega_j)$. The result follows.  
\end{proof}

\begin{Th}\label{th_lambda_atomic_length}
Let $w \in W$, $\underline{w}\in\Red(w)$ and $\lambda \in P^+$.
We have 
$$\lambda-w(\lambda) = \sum\limits_{\beta \in N_{\lambda}(\underline{w})}\beta
\quad\text{ and }\quad
\sL_{\lambda}(w) = \sum\limits_{\beta \in N_{\lambda}(\underline{w})} \h(\beta).
$$
In particular, in both identities, the right hand side does not depend on the choice
of the reduced expression.
\end{Th}

\begin{proof}
 Let us write $\lambda = \sum\limits_{j=1}^nm_j\omega_j$. By Lemma \ref{Lemma omega_j} we know that $\omega_j - w(\omega_j) = \sigma_j(\underline{w})$.  It follows then that
\begin{align*}
\lambda - w(\lambda)  & = \lambda - w\left(\sum\limits_{j=1}^n m_j\omega_j\right) 
                                       = \lambda - \sum\limits_{j=1}^n m_jw(\omega_j) = \lambda - \sum\limits_{j=1}^nm_j(\omega_j -  \sigma_j(\underline{w})) \\
                                       &=\sum\limits_{j=1}^n m_j\sigma_j(\underline{w}) = \sum\limits_{j=1}^n m_j\sum\limits_{k=1}^{p_j} \underline{w}_{j,k}(\al_j) = \sum\limits_{j=1}^n\sum\limits_{k=1}^{p_j} m_j\underline{w}_{j,k}(\al_j)\\
                                       & =\sum\limits_{\beta \in N_{\lambda}(\underline{w})}\beta,
\end{align*}
which proves the first identity.
By definition, $\sL_\la(w) = \langle \la-w(\la), \rho^\vee\rangle$,
so combining the first identity and (\ref{inner_prod_height}), we get the second identity.
\end{proof}

Finally, we give analogues of \Cref{inversions_weak_order_1} and \Cref{inversions_weak_order_2} for $\la$-inversion sets.
Recall that $\leq$ denotes the weak order on $W$,
that is, $w\leq w'$ if and only if there exist reduced expressions $\underline{w}, \underline{w'}$
such that $\underline{w}$ is a prefix of $\underline{w'}$.

\begin{Prop}\label{la_inv_weak_order}
Let $\la=\sum_{i\in I} m_i\om_i\in P$.
\begin{enumerate}
\item Let $\underline{w}, \underline{w'}\in \Red(W)$. Assume that $m_j\neq 0$ for all $j$.
If $N_\la(\underline{w}) \subseteq N_\la(\underline{w'})$,
then $w\leq w'$.
\item 
Let $w, w'\in W$. If $w\leq w'$, then $N_\la(\underline{w}) \subseteq N_\la(\underline{w'})$ where $\underline{w}$, $\underline{w'}$
are chosen such that $\underline{w}$ is a prefix of $\underline{w'}$.
\end{enumerate}
\end{Prop}

\begin{proof}\
\begin{enumerate}
\item 
Assume that $N_\la(\underline{w}) \subseteq N_\la(\underline{w'})$
Let us show that $\Leftrightarrow N(w) \subseteq N(w')$.
Take $\underline{w}_{j,k}(\al_j) \in N(w)$.
By definition, $m_j \underline{w}_{j,k}(\al_j) \in N_\la(\underline{w})\subseteq N_\la(\underline{w'})$.
So we can write $m_j \underline{w}_{j,k}(\al_j)
=m_{j'} \underline{w'}_{j',k'}(\al_{j'})$ for some $j', k'$.
Since $\underline{w}_{j,k}(\al_{j}), \underline{w'}_{j',k'}(\al_{j'})\in \Phi^+$, we must have $m_j=m_{j'}$ since
these are nonzero by assumption.
This yields 
$\underline{w}_{j,k}(\al_j) = \underline{w'}_{j',k'}(\al_{j'}) \in N(w')$. So we have proved that $N(w) \subseteq N(w')$.
By \cite[Corollary 2.10]{SRLE}, $N(w)\subseteq N(w') \Leftrightarrow w\leq w'$, which concludes the proof.
\item 
Let $w, w'\in W$ such that $w\leq w'$,
and choose reduced expressions $\underline{w}, \underline{w'}$ such that $\underline{w}$ is a prefix of $\underline{w'}$.
This implies that for all $1\leq j\leq n$, $1\leq k\leq p_j$,
$\underline{w}_{j,k}(\al_j) = \underline{w'}_{j,k}(\al_j)$.
Now, let $m_j \underline{w}_{j,k}(\al_j)\in N_\la(\underline{w})$.
This gives $m_j \underline{w}_{j,k}(\al_j) 
= m_j \underline{w'}_{j,k}(\al_j) 
\in N_\la(\underline{w'})$ by definition.
\end{enumerate}
\end{proof}

\begin{Cor}\label{cor_inv_weak_order}
Let  $\la\in P^+$
and $w, w'\in W$ such that $w\leq w'$. 
Then $\sL_\la (w) \leq \sL_\la(w')$.
\end{Cor}

\begin{proof}
We know by Theorem \ref{cor_inv_weak_order} 
that $\sL_{\lambda}(w) = \sum\limits_{\beta \in N_{\lambda}(\underline{w})} \h(\beta)$.  
Now, using Proposition \ref{la_inv_weak_order} (2), 
$w \leq w'$ implies $N_\la(\underline{w}) \subseteq N_\la(\underline{w'})$.  
It follows that $\sL_{\la}(w) \leq \sL_{\la}(w')$.
\end{proof}

\begin{Rem}
In particular, \Cref{la_inv_weak_order} enables us to recover directly \Cref{atomic_longest_element},
since $w_0$ is the largest element with respect to the weak order.
\end{Rem}


\section{Susanfe reflections}\label{sec_susanfe}

In this whole section, we use 
the notation of \Cref{sec_Cox}.
In particular, $W$ is a finite Weyl group,
and $S$ denotes a set of simple reflections.
This part is devoted to the study of some particular reflections of $W$ which we call ``Susanfe'',
and their properties with respect to the atomic length $\sL$.
Eventually, this will enable us to give in \Cref{sec_surj_AL} a proof of \Cref{thm_surj} which generalises
that of \cite{SackUlfarsson2011} for type $A_n$. 
In the following, if $W_A$ is a reflection subgroup of $W$, we will simply denote by $\sL_A$ the restriction to $W_A$
of the map $\sL :W\to \N$.

\subsection{General Susanfe theory}

\begin{Prop}\label{Atomic length parabolic}
Let $W_I$ be a standard parabolic subgroup of $W$.  Let $w \in W_I$. Then $\sL(w) = \sL_{I}(w)$.  
This fails in general for reflection subgroups that are not standard parabolic subgroups.
\end{Prop}

\begin{proof}
By definition we know that $\sL(w) = \sum\limits_{\al \in N(w)}\h(\al)  $ and $\sL_{I}(w) = \sum\limits_{\al \in N_{I}(w)}\h(\al)$.  Since $w \in W_I$, we know  by the second point of Lemma \ref{lemma Shi coeff dans J} that $N_{I}(w) = N(w)$.  Then the result follows.  The reason that we don't necessary have this equality for a reflection subgroup $W_A$ comes from the fact that the equality $N_{A}(w) = N(w)$ isn't true any more for $w \in W_A$.  The deep reason for 
this phenomenon is something well known: a reduced expression of $w$ in $W_A$ is a priori not a reduced expression of $w$ in $W$, and then the first concept used in the proof of Lemma \ref{lemma Shi coeff dans J} (1) cannot be used in $W_A$.
\end{proof}

\begin{Exa}
Let us illustrate the fact that,  in general, $\sL(w)\neq \sL_A(w)$ for $w\in W_A$.
Take $W$ of type $A_3$,  $A = \{s_1s_2s_1, s_3\}$ and $w = s_1s_2s_1s_3$.  
By Proposition \ref{Decomposition_N} we get easily that $N(w) = \{e_{12},e_{13},e_{23},e_{14}\}$. 
Therefore 
$$\sL(w) = \mathrm{ht}(e_{12}) +  \mathrm{ht}(e_{13}) + \mathrm{ht}(e_{23}) + \mathrm{ht}(e_{14})= 1 + 2 +1 +3 = 7.$$ 
Now, the group $W_A = \langle s_1s_2s_1, s_3 \rangle$ is a Coxeter group of type $A_2$ with simple system $S_A = \{s_1s_2s_1, s_3\}$, with set of simple roots $\Delta_A = \{e_{13}, e_{34}\}$ and set of positive roots $\Phi_A^+ = \{e_{13}, e_{34},e_{14}\}$.  The expression $s_1s_2s_1 s_3$ is a reduced expression of $w$ in $W_A$ of length $2$, i.e., $\ell_A(w) = 2$,  but it is a reduced expression of length $4$ in $W$.  Moreover, since $w \in W_A$,  it follows by Proposition \ref{propDy} that $N_A(w) = N(w) \cap \Phi_A = \{e_{13},e_{14}\}$. Therefore 
$$\sL_A(w) = \mathrm{ht}(e_{13}) + \mathrm{ht}(e_{14}) = 2+3 = 5,$$
and we see that $\sL(w) \neq \sL_A(w)$.
\end{Exa}
 
For $w \in W$, we denote $\fix_w  = \{\al \in \Phi^+\mid w(\al) = \al\}$, the set of positive roots fixed by $w$,
and $\overline{\fix_w}=\Phi^+ \backslash \fix_w$ its complement.

\begin{Def} An element $w\in W$ is called \textit{Susanfe} if $N(w) = \overline{\fix_w }$.
\end{Def}

\begin{Prop}\label{highest_reflection}
There exists a Susanfe reflection in $W$. 
In particular, the reflection associated to the highest root $\tal$ is a Susanfe reflection.
\end{Prop}

\begin{proof}
Let $t$ be the reflection associated to the highest root, so that
$t(\al) = \al - 2\displaystyle\frac{(\al \mid \tal)}{(\tal \mid \tal)}\tal$ ~for all $\al\in \Phi^+$.
Since $\Phi$ is crystallographic, we have ~$2\displaystyle\frac{(\al \mid \tal)}{(\tal\mid \tal)} \in \mathbb{Z}$. 
It is known (see \cite[Ch. VI, \S 1, Sec. 8 ]{BOURB}) that $\tal$ belongs to the fundamental chamber 
$C_0 = \{x \in V \mid  \left( x\mid\beta\right) \geq 0 ~ \forall \beta \in \Delta\}$. 
Therefore it follows that $2\displaystyle\frac{(\al \mid \tal)}{(\tal \mid \tal)} \in \N$.  
We must show that $N(t) = \overline{\mathrm{Fix}_t }$.  
The inclusion $N(t) \subseteq \overline{\fix_t}$ is obvious,
because if there were a root $\alpha \in N(t)$  such that $\al \in \mathrm{Fix}_t$, 
then we would have $t(\al) = \al$ and thus $t^{-1}(\al) = t(\al) \in \Phi^+$, which condraticts $\al\in N(t)$.
Let us now show the reverse inclusion.  
Let $\al \in \overline{\fix_t}$. 
There are two cases.
\begin{enumerate}
\item If $( \al \mid \tal) = 0$ then $t(\al) = \al$, which shows that $\al \in \fix_t$.
Just as above, this is a contradiction and therefore this situation cannot occur.
\item  If  $( \al \mid \tal ) > 0$ then $2\displaystyle\frac{(\al \mid \tal)}{(\tal \mid \tal)} > 0$ 
and then $\al - 2\displaystyle\frac{(\al \mid \tal)}{(\tal \mid \tal)}\tal$
belongs to $\Phi^-$ since $\tal$ is the highest root.  
It follows that $t(\al) \in \Phi^-$ and thus $\al \in N(t)$.  
Hence $\overline{\fix_t} \subseteq N(t)$. 
\end{enumerate}
\end{proof}

\begin{Rem}
In type $B_n$, 
computer experiments in SageMath suggest that the number of Susanfe elements 
should be given by $F_n-1$, where $F_n$ is the $n$-th Fibonacci number.
\end{Rem}

\begin{Prop}\label{Proposition Fondamentale}
Let $B\subseteq T$, so that
$W_B$ is a reflection subgroup of $W$, 
let $w \in W_B$ and let $t \in T$ be a Susanfe reflection. 
Consider the reflection subgroup $W_A= t W_B t$.  
Then we have $N(tw) =  N_{A}((tw)_A) \sqcup \left(N(t) \cap \overline{\Phi_A}\right)$.
\end{Prop}

\begin{proof}
One has $N(tw) = [N(tw) \cap \Phi_A] \sqcup [N(tw) \cap \overline{\Phi_A}]$.  Using  Proposition \ref{propDy} we know that $N(tw) \cap \Phi_A = N_{A}((tw)_A)$. Thus $N(tw) = N_{A}((tw)_A) \sqcup [N(tw) \cap \overline{\Phi_A}]$.  We consider now the following decomposition:
$$
N(tw) \cap \overline{\Phi_A} = N(tw) \cap \overline{\Phi_A} \cap \fix_t  \sqcup N(tw) \cap \overline{\Phi_A} \cap \overline{\fix_t } 
$$
Let us look at both terms of the right-hand side separately.
\begin{itemize}
\item We claim that $N(tw) \cap \overline{\Phi_A} \cap \fix_t =\emptyset$.  
Let us first show that
\begin{equation}\label{eq1}
\overline{\Phi_A} \cap \fix_t  \subseteq \overline{\Phi_B}.
\end{equation}  
Let $\al \in \overline{\Phi_A} \cap \fix_t $ such that $\al \in \Phi_B$. Then $t(\al) \in t(\Phi_B)$ and we  know by Proposition \ref{root sytem 
reflection} that $t(\Phi_B) = \Phi_A$. Thus $t(\al) \in \Phi_A$, 
and since $\al \in \fix_t $ we have $t(\al) = \al$, which implies that $\al \in \Phi_A$. This is a contradiction and we indeed obtain \Cref{eq1}.
Now, assume there exists $\al\in N(tw) \cap \overline{\Phi_A} \cap \fix_t$.
In particular, we have $k(tw, \al)=-1$, and 
$k(t,\al)=0 $ because $t^{-1}(\al) = t(\al) = \al \in \Phi^+$, that is $\al \notin N(t)$.
Thus, by the  formula of Proposition  \ref{formula_Nathan} we obtain that 
$-1 = k(tw,\al) = k(w,\al) + k(t,\al)= k(w,\al) + 0$,  that is $k(w,\al) =-1$.
But now by \Cref{eq1}, we have $\al\in\overline{\Phi_B}$. Thus, since $w \in W_B$, we can apply \Cref{lemma Shi coeff dans J} 
and obtain $k(w,\al)=0$, which is a contradiction.

\item Consider now $N(tw) \cap \overline{\Phi_A} \cap \overline{\fix_t }$.
Since $t$ is Susanfe, we have  $\overline{\fix_t } = N(t)$, and therefore
$N(tw) \cap \overline{\Phi_A} \cap \overline{\fix_t } = N(tw) \cap \overline{\Phi_A} \cap  N(t)$.
Now, we claim that
\begin{equation} \label{eq2}
N(t)\cap \overline{\Phi_A} \subseteq N(tw).
\end{equation}
So let $\al\in N(t)\cap \overline{\Phi_A}$, and let us show that $k(tw,\al) = -1$.
Since $k(tw,\al) = k(w,t(\al)) + k(t,\al)$ (by \Cref{formula_Nathan}) 
and $k(t,\al) = -1$ because $\al \in N(t)$,  it suffices to show that $k(w,t(\al)) = 0$.
Since  $\al \in \overline{\Phi_A}$, we have $t(\al) \in \overline{t(\Phi_A)}$ and thus
$t(\al) \in \overline{\Phi_B}$ by Proposition \ref{root sytem reflection}. 
Finally, we can use \Cref{lemma Shi coeff dans J} since $w \in W_B$, which yields $k(w,t(\al)) = 0$.
This proves \Cref{eq2}, and we therefore obtain 
 $N(tw) \cap \overline{\Phi_A} \cap \overline{\fix_t } = N(t) \cap \overline{\Phi_A}$.
\end{itemize}
In the end, we obtain the desired identity
$N(tw) =  N_{A}((tw)_A) \sqcup \left(N(t) \cap \overline{\Phi_A}\right)$.
\end{proof}

From \Cref{Proposition Fondamentale} we obtain the following crucial corollary, which will be used in \Cref{sec_surj_AL}.
 For $A \subseteq  T$ and $w\in W$, we denote
$$\sL(w, A) = \sum_{\al\in N(w)\cap \overline{\Phi_A}} \h(\al),$$
which corresponds to restraining the computation of the atomic length to roots outside $\Phi_A$.

\begin{Cor}\label{corollary principal}
We keep the assumptions of Proposition \ref{Proposition Fondamentale}.
Then  $$\sL(tw) = \sL_{A}((tw)_A) + \sL(t, A).$$
\end{Cor}

\subsection{Some particular Susanfe reflections}

In the following, we use the convenient notation
introduced in \Cref{Weyl_groups}.
We still use the conventions of \cite{BOURB} for the labelling of the Dynkin diagrams.

\begin{Lem}[Type $A_n$]\label{Lemma A}
Let $t$ be the reflection corresponding to the highest root $\tal$ 
and let $I = \{s_2,\dots,s_n\}$.  Then one has the following reduced expressions
and properties.
\begin{enumerate}
\item  $t = s_1s_2\dots s_{n-1}s_ns_{n-1}\dots s_2s_1$.
\item $t_I = s_ns_{n-1}\dots s_2$ and $^It = s_1s_2\dots s_{n-1}s_n$.
\item $N(t) = \{e_{12}, e_{13}, \dots, e_{1,n+1}\} \cup \{e_{2,n+1},e_{3,n+1},\dots,e_{n,n+1}\}$.
\item $N({}^It) = \{e_{12}, e_{13}, \dots, e_{1,n+1}\}$, in particular $ N(t) \cap \overline{\Phi_I}=N({}^It)$.
\item $\sL(t,I) = \sL({}^It) = \binom{n+1}{2}.$
\end{enumerate}
\end{Lem}

\begin{proof} Recall that in type $A_n$, all roots $\al$ verify $(\al\mid \al) =2$, 
so in particular $\langle x, \al^\vee\rangle = (\al\mid x)$ for all $x\in V$.
This implies that $s_\al(x) = x-(\al\mid x)\al$. 
In particular, we have $\tal= e_{1,n+1}$, and $t(x) = x-(e_{1,n+1} \mid x) e_{1,n+1}$.
\begin{enumerate}
\item  
We know that for any $k$, the simple root corresponding to $s_k$ is $e_{k,k+1}$. 
Moreover, 
In type $A_n$, we know that $\tal = e_{1,n+1}$, so $t=s_{e_{1,n+1}}$.
On the other hand, $s_1s_2\dots s_{n-1}s_ns_{n-1}\dots s_2s_1 = s_{s_1s_2\ldots s_{n-1}(e_{n,n+1})}$,
so it suffices to show that $s_1s_2\dots s_{n-1}(e_{n,n+1}) = e_{1,n+1}$.
This equality follows directly by iterating from right to left, 
since $s_{j}(e_{j+1, n+1}) = e_{j,n+1}$.
The fact that this expression of $t$ is reduced will be proved in (2) just below.
\item 
We need to show that 
\begin{enumerate}
 \item $t = (s_ns_{n-1}\dots s_2 )( s_1s_2\dots s_{n-1}s_n),$
\item $s_ns_{n-1}\dots s_2 \in W_I$ and is a reduced expression,
\item $s_1s_2\dots s_{n-1}s_n \in {}^{I}W$  and is a reduced expression.
\end{enumerate}
Using the braid relations of type $A_n$, we have
\begin{align*}
s_1s_2\dots \textcolor{purple!80}{s_{n-1}s_ns_{n-1}}\dots s_2s_1 & = s_1s_2\dots \textcolor{purple!80}{s_{n}s_{n-1}s_{n}}\dots s_2s_1 \\
& = \textcolor{purple!80}{s_n}s_1s_2\dots \textcolor{teal}{s_{n-2}s_{n-1}s_{n-2}}\dots s_2s_1\textcolor{purple!80}{s_n} \\
& = \textcolor{purple!80}{s_n}s_1s_2\dots \textcolor{teal}{s_{n-1}s_{n-2}s_{n-1}}\dots s_2s_1\textcolor{purple!80}{s_n} \\
& = \textcolor{purple!80}{s_n}\textcolor{teal}{s_{n-1}}s_1s_2\dots \textcolor{orange!89}{s_{n-3}s_{n-2}s_{n-3}}\dots s_2s_1\textcolor{teal}{s_{n-1}}\textcolor{purple!80}{s_n} \\
& = \dots \\
& =  \textcolor{purple!80}{s_n}\textcolor{teal}{s_{n-1}}\textcolor{orange!89}{s_{n-2}}\dots s_2 \cdot s_1s_2 \dots \textcolor{orange!89}{s_{n-2}}\textcolor{teal}{s_{n-1}}\textcolor{purple!80}{s_n},
\end{align*}
which proves (a).
Moreover, by definition of $I$, $s_ns_{n-1}\dots s_2 \in W_I$ and it is also clear that this expression is reduced, so we have (b).  
Finally, it is clear that the expression $s_1s_2\dots s_{n-1}s_n$ is already reduced and has no reduced expression beginning by a letter in $I$, 
implying that $\ell(ss_1s_2\dots s_{n-1}s_n) > \ell(s_1s_2\dots s_{n-1}s_n)$ for any $s \in I$. 
This means precisely that $s_1s_2\dots s_{n-1}s_n \in~^{I}W$.
Therefore, we have found the $I$-decomposition of $t$ where the two factors are expressed in a reduced form.
By \Cref{rem_A_dec}, we are ensured that $\ell(t) =  \ell(t_I) + \ell({}^It)=2n-1$.
Therefore the expression $t = s_1s_2\dots s_{n-1}s_ns_{n-1}\dots s_2s_1$, which also uses $2n-1$ generators, is also reduced.

\item From the formula $ t(e_{ij}) =  e_{ij} - ( e_{ij} \mid e_{1,n+1}) e_{1,n+1}$,
we see that the set of positive roots that are fixed under $t$ is 
$\fix_t  = \{e_{ij}\mid 2\leq i<j\leq n\}$, so that
$\overline{\fix_t } = \{e_{1j}\mid 2\leq j \leq n+1\} \cup \{e_{i,n+1}\mid 2 \leq i \leq n\}.$ 
We conclude by using Proposition \ref{highest_reflection}.

\item Since $^It = s_1s_2\dots s_{n-1}s_n$ is a reduced expression,  we know by Proposition \ref{Decomposition_N} that 
$$
N({}^It) = \{\al_1, s_1(\al_2), s_1s_2(\al_3), \dots, s_1s_2\dots s_{n-1}(\al_n)\},
$$
 that is 
\begin{align*}
N({}^It) &= \{e_{12}, s_1(e_{23}), s_1s_2(e_{34}), \dots, s_1s_2\dots s_{n-1}(e_{n,n+1})\} \\
           & = \{e_{12},  e_{13},  e_{14}, \dots, e_{1,n+1}\}.
\end{align*}

\item The height of the root $e_{ij}$ is $j-i$. Therefore, using Point (4) above, we obtain
\begin{align*}
\sL({}^It) &=\sum\limits_{\al \in N({}^It)}\h( \al) =\sum\limits_{j=2}^{n+1}\h(e_{1j}) =  \sum\limits_{j=2}^{n+1}(j-1) =  \sum\limits_{j=1}^{n}j = \binom{n+1}{2}.
\end{align*}
\end{enumerate}
\end{proof}

\begin{Exa}
Take $n=4$. The highest root is 
$\tal = e_{15}$ and its corresponding reflection $t$ 
has the following reduced expression 
$$ t = s_1 s_2 s_3s_4s_3s_2s_1 = \textcolor{purple!80}{s_4s_3s_2} \textcolor{teal}{s_1s_2s_3s_4}$$ 
where the red part is in $W_I = \langle s_2,s_3,s_4 \rangle$ and the blue part is in $^IW$.  A convenient way to see this decomposition is via the Shi vectors.  
The Shi vector  corresponding to $t$ is given in \Cref{Shi_A}
and the $I$-decomposition of $t$ is given in \Cref{decomp_A}.

\begin{figure}[h!]
\begin{center}
\begin{tikzpicture} 
\node at (0,0) {$-1$} ;
\node at (1,0) {$0$} ;
\node at (2,0) {$0$} ;
\node at (3,0) {$-1$} ;
\node at (0.5,0.5 ) {$-1$} ;
\node at (1.5, 0.5 ) {$0$} ;
\node at (2.5, 0.5) {$-1$} ;
\node at (1.1,1) {$-1$} ;
\node at (2,1) {$-1$} ;
\node at (1.5,1.5 ) {$-1$} ;
\node at (5,0) {$e_{12}$} ;
\node at (6,0) {$e_{23}$} ;
\node at (7,0) {$e_{34}$} ;
\node at (8,0) {$e_{45}$} ;
\node at (5.5,0.5 ) {$e_{13}$} ;
\node at (6.5, 0.5 ) {$e_{24}$} ;
\node at (7.5, 0.5) {$e_{35}$} ;
\node at (6,1) {$e_{14}$} ;
\node at (7,1) {$e_{25}$} ;
\node at (6.5,1.5 ) {$e_{16}$} ;
\end{tikzpicture}
\end{center}
\caption{The Shi vector of $t$ in $A_4$. The right triangle gives the coordinates of the Shi vector.}
\label{Shi_A}
\end{figure}

\begin{figure}[h!]
\begin{center}
\begin{tikzpicture} 
\node at (0,0) {$-1$} ;
\node at (1,0) {$0$} ;
\node at (2,0) {$0$} ;
\node at (3,0) {$-1$} ;

\node at (0.5,0.5 ) {$-1$} ;
\node at (1.5, 0.5 ) {$0$} ;
\node at (2.5, 0.5) {$-1$} ;

\node at (1.1,1) {$-1$} ;
\node at (2,1) {$-1$} ;

\node at (1.5,1.5 ) {$-1$} ;

\node at (4,0.7) {$=$};

\node at (5,0) {$0$} ;
\node at (6,0) {$\textcolor{purple!80}{0}$} ;
\node at (7,0) {$\textcolor{purple!80}{0}$} ;
\node at (8,0) {$\textcolor{purple!80}{-1}$} ;

\node at (5.5, 0.5 ) {$0$} ;
\node at (6.5, 0.5) {$\textcolor{purple!80}{0}$} ;
\node at (7.5, 0.5) {$\textcolor{purple!80}{-1}$} ;

\node at (6,1) {$0$} ;
\node at (7,1) {$\textcolor{purple!80}{-1}$} ;

\node at (6.5,1.5 ) {$0$} ;

\node at (9,0.7) {$\bullet$};

\node at (10,0) {$\textcolor{teal}{-1}$} ;
\node at (1	1,0) {$0$} ;
\node at (12,0) {$0$} ;
\node at (13,0) {$0$} ;

\node at (10.5, 0.5 ) {$\textcolor{teal}{-1}$} ;
\node at (11.5, 0.5) {$0$} ;
\node at (12.5, 0.5) {$0$} ;

\node at (11,1) {$\textcolor{teal}{-1}$} ;
\node at (12,1) {$0$} ;

\node at (11.5,1.5 ) {$\textcolor{teal}{-1}$} ;
\end{tikzpicture}
\end{center}
\caption{The $I$-decomposition of $t$ in terms of Shi vectors in $A_4$.}
\label{decomp_A}
\end{figure}
\end{Exa}

\begin{Lem}[Type $B_n$]\label{lemma B}
Let t be the reflection of the highest root $\tal$ and let $I = \{s_2,\dots,s_n\}$.  Then one has the following reduced expressions
and properties.
\begin{enumerate}
 \item  $t = s_2s_3\dots s_{n-1}s_ns_{n-1}\dots s_3s_2\cdot s_1 \cdot s_2s_3\dots s_{n-1}s_ns_{n-1}\dots s_3s_2$.
\item $t_I = s_2s_3\dots s_{n-1}s_ns_{n-1}\dots s_3s_2$ and $^It = s_1s_2\dots s_{n-1}s_ns_{n-1}\dots s_3s_2$.
\item $N(t) = \{e_{13},e_{14},\dots,e_{1n},e_1,e^{1n}, e^{1,n-1},\dots, e^{12}\} \cup \{e_{23},e_{24},\dots,e_{2n},e_2,e^{2n}, e^{2,n-1},\dots, e^{23}\}$.
\item $N({}^It) = \{e_{12},  e_{13}, \dots, e_{1n}\} \sqcup \{e_{1}\} \sqcup \{e^{1n}, e^{1,n-1},\dots, e^{13}\}$.
\item Let $t' =$ $^Its_1$. Then $t'$ is a Susanfe reflection that satisfies $N(t') \cap \overline{\Phi_I} = N(t')$.
\item $\sL(t',I) = \sL(t') = 2n^2 - n$.
\end{enumerate}
\end{Lem}

\begin{proof} In type $B_n$ the highest root is given by $\tal = e^{12}$, and we have
 $t(x) = x - (e^{12} \mid x) e^{12}$.
\begin{enumerate}
\item  
As in type $A_n$, it suffices to show that $s_2s_3\dots s_{n-1}s_ns_{n-1}\dots s_3s_2(\al_1) = e^{12}$.  
First of all we have  $s_{n}(e_{1n}) = e_1 -e_n - 2\left( e_1-e_n \mid e_{n}\right) e_n = e_1 +e_n = e^{1n}$. 
Moreover, the same kind of computations as above shows that $s_k(e^{1,k+1}) = e^{1k}$ for any $1 \leq k  \leq n-1$. Therefore we obtain that
\begin{align*}
s_2s_3\dots s_{n-1}s_ns_{n-1}\dots s_3s_2(\al_1) & = s_2s_3\dots s_{n-1}s_ns_{n-1}\dots s_3s_2(e_{12}) \\
																				   & = s_2s_3\dots s_{n-1}s_ns_{n-1}\dots s_3(e_{13}) \\
																				   & = \cdots \\
																				   & =  s_2s_3\dots s_{n-2}s_{n-1}s_n(e_{1n}) \\
																				   & =  s_2s_3\dots s_{n-2}s_{n-1}(e^{1n}) \\
																				   & = s_2s_3\dots s_{n-2}(e^{1,n-1}) \\
																				   & = \cdots \\
																				   & = s_2(e^{13}) \\
																				   & = e^{12}.
\end{align*}
The fact that this expression of $t$ is reduced is again proved in (2) below.

\item It is clear that $s_2s_3\dots s_{n-1}s_ns_{n-1}\dots s_3s_2 \in W_I$.  
Further, since $s_1s_2\dots s_{n-1}s_ns_{n-1}\dots s_3s_2$ starts with $s_1$ and since $s_1$ appears only once, there is no reduced expression of this word beginning with a letter in $I$, which implies that 
$\ell(ss_1s_2\dots s_{n-1}s_ns_{n-1}\dots s_3s_2) > \ell(s_1s_2\dots s_{n-1}s_ns_{n-1}\dots s_3s_2)$ for any $s \in I$,
so $s_1s_2\dots s_{n-1}s_ns_{n-1}\dots s_3s_2 \in {}^IW$.  
Since we obviously have $$t = (s_2s_3\dots s_{n-1}s_ns_{n-1}\dots s_3s_2 ) ( s_1s_2\dots s_{n-1}s_ns_{n-1}\dots s_3s_2),$$  
by uniqueness of the $I$-decomposition it follows that we have  $t_I = s_2s_3\dots s_{n-1}s_ns_{n-1}\dots s_3s_2$ and 
$^It = s_1s_2\dots s_{n-1}s_ns_{n-1}\dots s_3s_2$.  
Finally, it is also clear that these two expressions are reduced and therefore the expression given for $t$ is also reduced
(by \Cref{rem_A_dec} again).

\item  Since $t$ is Susanfe by \Cref{highest_reflection}, $N(t) = \overline{\fix_t}$, so it suffices to determine the set of 
positive roots fixed by $t$.
For any $\al \in \Phi^+$, we have $t(\al) = \al - (e^{12} \mid \al) e^{12}$,
so $\al\in\fix_t$ if and only if $(e^{12}\mid \al)=0$,
which happends exactly when there are no $1$ or $2$ appearing in the labels of the roots, 
except for $\al = e_{12}$ (since  $(e_{12} \mid e^{12} ) = 0$). 
Therefore  $\overline{\fix_t }=N(t)$ is the set of positive roots that have a label $1$ or a label $2$ with the exception of $e_{12}$,
which is precisely the desired description.

\item Let us fix $u = s_1s_2\dots s_{n-1}s_n$ and $v = s_ns_{n-1}\dots s_3s_2$.  Since $ s_1s_2\dots s_{n-1}s_ns_{n-1}\dots s_3s_2$ is a reduced expression of ${}^It$ by Point (2), it follows that $u$ and $v$ are also reduced. 
Therefore by \Cref{Decomposition_N}, we have $N({}^It) = N(u) \sqcup uN(v)$.  
Let us compute at first $N(u)$.  Since the expression of $u$ is reduced we have again by \Cref{Decomposition_N}
$$
N(u) =  \{e_{12}, s_1(e_{23}), \dots, s_1s_2\dots s_{n-1}(e_{n})\}.
$$
Moreover the generators $s_1,\dots, s_{n-1}$ are of the same form,  
and the computation of the first $n-1$  terms straightforward. 
The only thing to be cautious about is the last element $s_1 s_2 \ldots s_{n-1} (e_n)$. 
An easy calculation shows that $s_{k-1}(e_k) = e_{k-1}$ for any $k=2,3,\dots,n$. 
Iterating therefore yields  $s_1 s_2 \ldots s_{n-1} (e_n) = e_1$.  Hence
$$
 N(u) = \{e_{12},  e_{13}, \dots, e_{1n},e_{1}\}.
$$
We now compute $uN(v)$. Using \Cref{Decomposition_N} again, we obtain $N(v) = \{e_{n-1,n}, e_{n-2,n},\dots, e_{2n}\}.$
A direct inductive computation yields $$uN(v) = \{e^{1n}, e^{1,n-1},\dots, e^{13}\}.$$ 
Finally we get
\begin{align*}
N({}^It) =N(u) \sqcup uN(v) = \{e_{12},  e_{13}, \dots, e_{1n},e_{1}\} \sqcup \{e^{1n}, e^{1,n-1},\dots, e^{13}\}.
\end{align*}

\item   
The expression ${}^I t s_1 = s_1s_2\dots s_{n-1}s_ns_{n-1}\dots s_3s_2s_1$ is clearly reduced, since
there is no relation in the Coxeter group of type $B_n$ that can be applied.
Note that this means that this is actually the only reduced expression of ${}^I t s_1$.
Therefore, we can apply \Cref{Decomposition_N} and we get
\begin{align*}
N(t') & = N({}^Its_1) = N({}^It) ~\sqcup ~^It N(s_1) = N({}^It) ~\sqcup ~^It(\al_1) \\
\end{align*}
Moreover we have
\begin{align*}
^It(\al_1) & = s_1s_2\dots s_{n-1}s_ns_{n-1}\dots s_3s_2(e_{12}) \\
                      & = s_1s_2\dots s_{n-1}s_ns_{n-1}\dots s_3(e_{13}) \\
                      & = \cdots \\
                      & = s_1s_2\dots s_{n-1}s_n(e_{1n}) \\
                      & = s_1s_2\dots s_{n-1}(e^{1n}) \\
                      & = s_1s_2\dots s_{n-2}(e^{1,n-1}) \\
                      & = \cdots \\ 
                      & = s_1(e^{12}) \\
                      & = e^{12}.
\end{align*}
Finally we obtain 
$$
N(t') = \{e_{12},  e_{13}, \dots, e_{1n},e_{1}\} \sqcup \{e^{1n}, e^{1,n-1},\dots, e^{13}, e^{12}\}.
$$

Since $\overline{\Phi_I}$ is the set of roots with a label $1$
and since all the roots with a label $1$ are the ones given in $N(t')$, it follows that $N(t') = \overline{\Phi_I^+}$ and thus 
$N(t') \cap \overline{\Phi_I} = N(t')$.  
It remains to show that $t'$ is a Susanfe reflection.  
Clearly, $t'$ is a reflection, since $w s_n w^{-1}$ is a reflection for all $w\in W$.
Let us prove that $\overline{\fix_{t'}(\Phi^+)} = N(t')$.  
For all $k =1,2,\dots, n$ we have 
\begin{align*}
s_{k-1}(e_k) & = e_k - \left( e_k \mid e_{k-1}-e_k \right) (e_{k-1}-e_k) \\
& = e_{k-1}.
\end{align*}

This implies that $s_{k-1}s_{e_k}s_{k-1} = s_{e_{k-1}}$ for any $k =1,2,\dots, n$. Therefore, it follows by induction that
\begin{align*}
t' &= s_1s_2\dots s_{n-2}\textcolor{purple!80}{s_{n-1}s_{e_n}s_{n-1}}s_{n-2}\dots s_2s_1 \\ 
   &= s_1s_2\dots s_{n-2}\textcolor{purple!80}{s_{e_{n-1}}}s_{n-2}\dots s_2s_1 \\
   & = \dots \\
   & = s_{e_1}.
\end{align*}
Thus the roots $\al \in \Phi_I$ such that $t(\al) = \al$ are exactly the roots
satisfying $\left( e_1 \mid  \al\right) = 0$,
that is to say, the roots with no label $1$. 
We already know that this is $\Phi_I$. 
Therefore $\fix_{t'} = \Phi_I^+$ and then $\overline{\fix_{t'}} = \overline{\Phi_I^+} = N(t')$.

\item 
We have 
$$\sL(t') = \sum_{\al \in N(t')} \mathrm{ht}(\al) = \sum_{\al\in N(t')\cap \overline{\Phi_I}} \al
= \sL(t',I) .$$ 
since $N(t')\cap \overline{\Phi_I} = N(t')$ by (5).
Now, we know (see \cite{BOURB}) that $\h(e_{ij}) = j-i$, $\h(e_i) = n-i+1$ and $ \h(e^{ij}) = 2(n+1)-(i+j)$ in type $B_n$. 
Therefore it follows that
\begin{align*}
\sL(t') &=\sum\limits_{j=2}^{n}\h(e_{1j})  + \h(e_1) + \sum\limits_{j=2}^{n}\h(e^{1j})  =  \sum\limits_{j=2}^{n}(j-1) + n +  \sum\limits_{j=2}^{n}[2(n+1)-(1+j)] \\
& =  \sum\limits_{j=1}^{n-1}j + n +  \sum\limits_{j=2}^{n}[2n+1-j]  =  \sum\limits_{j=1}^{n-1}j + n + (n-1)(2n+1) - \sum\limits_{j=2}^{n}j \\
&  =  \sum\limits_{j=1}^{n-1}j + n + (n-1)(2n+1) -[( \sum\limits_{j=1}^{n-1}j) -1 +n]  =  n + (n-1)(2n+1)  + 1 -n \\
& = (n-1)(2n+1)  + 1  \\
& = 2n^2 - n.
\end{align*}
\end{enumerate}

\end{proof}

\begin{Exa} Take $n=4$.  
The highest root is $\tal = e^{12}$ and its corresponding reflection $t$ has the following reduced expression 
$$ t = \textcolor{purple!80}{s_2 s_3s_4s_3s_2} \textcolor{teal}{s_1s_2s_3s_4s_3s_2}$$ 
where the red part is in $W_I = \langle s_2,s_3,s_4 \rangle$ and the blue part is in $^IW$.  
The corresponding Shi vector of $t$ is given in \Cref{ShiB4}.
\begin{figure}[h!]
\begin{center}
\begin{tikzpicture} 
\node at (0,0) {$0$} ;
\node at (1,0) {$-1$} ;
\node at (2,0) {$0$} ;
\node at (3,0) {$0$} ;
\node at (0.5,0.5 ) {$-1$} ;
\node at (1.5, 0.5 ) {$-1$} ;
\node at (2.5, 0.5) {$0$} ;
\node at (0.5,1 ) {$-1$} ;
\node at (1.5, 1 ) {$-1$} ;
\node at (2.5, 1) {$0$} ;
\node at (1, 1.5 ) {$-1$} ;
\node at (2, 1.5) {$-1$} ;
\node at (1, 2 ) {$-1$} ;
\node at (2, 2) {$-1$} ;
\node at (1.5, 2.5) {$-1$} ;
\node at (1.5, 3) {$-1$} ;
\node at (5,0) {$e_{12}$} ;
\node at (6,0) {$e_{23}$} ;
\node at (7,0) {$e_{34}$} ;
\node at (8,0) {$e_4$} ;
\node at (5.5,0.5 ) {$e_{13}$} ;
\node at (6.5, 0.5 ) {$e_{24}$} ;
\node at (7.5, 0.5) {$e_3$} ;
\node at (5.5,1 ) {$e_{14}$} ;
\node at (6.5, 1 ) {$e_2$} ;
\node at (7.5, 1) {$e^{34}$} ;
\node at (6, 1.5 ) {$e_1$} ;
\node at (7, 1.5) {$e^{24}$} ;
\node at (6, 2 ) {$e^{14}$} ;
\node at (7, 2) {$e^{23}$} ;
\node at (6.5, 2.5) {$e^{13}$} ;
\node at (6.5, 3) {$e^{12}$} ;
\end{tikzpicture}
\end{center}
\caption{The Shi vector of $t$ in $B_4$. The right part gives the coordinates of the Shi vector.}
\label{ShiB4}
\end{figure}
and the $I$-decomposition of $t$ is given in \Cref{ShidecB4}.

\begin{figure}[h!]
\begin{center}
\begin{tikzpicture} 
\node at (0,0) {$0$} ;
\node at (1,0) {$-1$} ;
\node at (2,0) {$0$} ;
\node at (3,0) {$0$} ;
\node at (0.5,0.5 ) {$-1$} ;
\node at (1.5, 0.5 ) {$-1$} ;
\node at (2.5, 0.5) {$0$} ;
\node at (0.5,1 ) {$-1$} ;
\node at (1.5, 1 ) {$-1$} ;
\node at (2.5, 1) {$0$} ;
\node at (1, 1.5 ) {$-1$} ;
\node at (2, 1.5) {$-1$} ;
\node at (1, 2 ) {$-1$} ;
\node at (2, 2) {$-1$} ;
\node at (1.5, 2.5) {$-1$} ;
\node at (1.5, 3) {$-1$} ;
\node at (4,1) {$=$} ;
\node at (5,0) {$0$} ;
\node at (6,0) {$\textcolor{purple!80}{-1}$} ;
\node at (7,0) {$\textcolor{purple!80}{0}$} ;
\node at (8,0) {$\textcolor{purple!80}{0}$} ;
\node at (5.5,0.5 ) {$0$} ;
\node at (6.5, 0.5 ) {$\textcolor{purple!80}{-1}$} ;
\node at (7.5, 0.5) {$\textcolor{purple!80}{0}$} ;
\node at (5.5,1 ) {$0$} ;
\node at (6.5, 1 ) {$\textcolor{purple!80}{-1}$} ;
\node at (7.5, 1) {$\textcolor{purple!80}{0}$} ;
\node at (6, 1.5 ) {$0$} ;
\node at (7, 1.5) {$\textcolor{purple!80}{-1}$} ;
\node at (6, 2 ) {$0$} ;
\node at (7, 2) {$\textcolor{purple!80}{-1}$} ;
\node at (6.5, 2.5) {$0$} ;
\node at (6.5, 3) {$0$} ;
\node at (9,1) {$\bullet$} ;
\node at (10,0) {$\textcolor{teal}{-1}$} ;
\node at (11,0) {$0$} ;
\node at (12,0) {$0$} ;
\node at (13,0) {$0$} ;
\node at (10.5,0.5 ) {$\textcolor{teal}{-1}$} ;
\node at (11.5, 0.5 ) {$0$} ;
\node at (12.5, 0.5) {$0$} ;
\node at (10.5,1 ) {$\textcolor{teal}{-1}$} ;
\node at (11.5, 1 ) {$0$} ;
\node at (12.5, 1) {$0$} ;
\node at (11, 1.5 ) {$\textcolor{teal}{-1}$} ;
\node at (12, 1.5) {$0$} ;
\node at (11, 2 ) {$\textcolor{teal}{-1}$} ;
\node at (12, 2) {$0$} ;
\node at (11.5, 2.5) {$\textcolor{teal}{-1}$} ;
\node at (11.5, 3) {$0$} ;
\end{tikzpicture}
\end{center}
\caption{The $I$-decomposition of $t$ in terms of Shi vectors in $B_4$.}\label{ShidecB4}
\end{figure}

Finally, the Shi vector of the Susanfe reflection $t'$ considered in Lemma \ref{lemma B} is given in \Cref{ShiprimeB4}.

\begin{figure}[h!]
\begin{center}
\begin{tikzpicture} 
\node at (0,0) {$-1$} ;
\node at (1,0) {$0$} ;
\node at (2,0) {$0$} ;
\node at (3,0) {$0$} ;
\node at (0.5,0.5 ) {$-1$} ;
\node at (1.5, 0.5 ) {$0$} ;
\node at (2.5, 0.5) {$0$} ;
\node at (0.5,1 ) {$-1$} ;
\node at (1.5, 1 ) {$0$} ;
\node at (2.5, 1) {$0$} ;
\node at (1, 1.5 ) {$-1$} ;
\node at (2, 1.5) {$0$} ;
\node at (1, 2 ) {$-1$} ;
\node at (2, 2) {$0$} ;
\node at (1.5, 2.5) {$-1$} ;
\node at (1.5, 3) {$-1$} ;
\end{tikzpicture}
\end{center}
\caption{The Shi vector of $t'$ in $B_4$.}\label{ShiprimeB4}
\end{figure}
\end{Exa}

\newpage

\begin{Lem}[Type $C_n$]\label{Lemma C}
Let t be the reflection of the highest root $\tal$ and let $I = \{s_2,\dots,s_n\}$.  Then one has the following reduced expressions
and properties.
\begin{enumerate}
 \item $t = s_1s_2\dots s_{n-1}s_ns_{n-1}\dots s_2s_1$.
\item $t_I = e$ and $^It = t$. In particular $^It$ is a Susanfe reflection.
\item $N({}^It) = \overline{\fix_t} = \overline{\Phi_I^+}$. In particular $ N(t) \cap \overline{\Phi_I}=N({}^It)$.
\item $\sL(t,I) = \sL({}^It) = 2n^2-n.$
\end{enumerate}
\end{Lem}

\begin{proof} In type $C_n$ the highest root is  $\tal = 2e_{1}$.
\begin{enumerate}
\item  The strategy is again the same as in type $A_n$ and $B_n$, namely showing that $s_1s_2\dots s_{n-1}(\al_n) = 2e_1$ with $\al_n = 2e_n$.
The only thing to take care of in this computation is how the reflections in type $C_n$ act on the roots.  
We leave the details to the reader since there is no difficulty.
Note that $t$ has the same reduced expression as the element $t'$ in \Cref{lemma B}.
By the exact same argument as in that proof, we are ensured that the expression $t=s_1s_2\dots s_{n-1}s_ns_{n-1}\dots s_2s_1$ is reduced
(and is in fact the only reduced expression of $t$).
\item 
We have just seen in (1) that $t=s_1s_2\dots s_{n-1}s_ns_{n-1}\dots s_2s_1$ is the unique reduced expression of $t$.
Therefore, for any $s\in I$, we have $\ell(st) > \ell(t)$,  that is to say $t \in {}^IW$.
\item The first equality comes from Proposition \ref{highest_reflection} and from the equality $^It = t$.  With respect to the second equality we first need to see that since $\Phi_I = \mathbb{Z}\langle e_{23},e_{34},\dots, e_{n-1,n},2e_n \rangle\cap \Phi$, each root of $\Phi_I$ does not have a label $i=1$ or $j=1$.  Moreover, any root $\al = e_{ij}$ or $2e_i$ or $e^{ij}$ that doesn't have a label $1$ (i.e.  any root in $\Phi_I$) satisfies 
$\left( \al \mid \tal \right) = 0$ and then $t(\al) = \al$, that is $\al \in \fix_{t}$ and then $\Phi_I = \fix_{t}$. 
\item In type $C_n$ we know (see \cite{BOURB}) that $\h(e_{ij}) = j-i$ and $\h(2e_i) = 2(n-i)+1$ and $\h(e^{ij}) = 2n+1-(i+j)$. Therefore it follows that
\begin{align*}
\sL(^Jt) &=\sum\limits_{j=2}^{n}\h(e_{1j})  + \sum\limits_{j=2}^{n}\h(e^{1j})  + \h(2e_1)=  \sum\limits_{j=2}^{n}(j-1)  +  \sum\limits_{j=2}^{n}[2n+1-(1+j)] + 2n-1 \\
& =  \sum\limits_{j=1}^{n-1}j +  \sum\limits_{j=2}^{n}(2n-j)  + 2n-1 =  \sum\limits_{j=1}^{n-1}j  + 2n(n-1)- \sum\limits_{j=2}^{n}j + 2n-1 \\
& = \sum\limits_{j=1}^{n}j - n  + 2n(n-1) - \sum\limits_{j=1}^{n}j + 1 + 2n-1 = 2n(n-1) +n\\
& = 2n^2-n.
\end{align*}
\end{enumerate}
\end{proof}

\begin{Exa}
Take $n=4$.  
The highest root is $\tal = 2e_{1}$ and its corresponding reflection $t$ has the following reduced expression 
$$ t = \textcolor{purple!80}{e} \textcolor{teal}{s_1s_2s_3s_4s_3s_2s_1}$$ where the red part is in $W_I = \langle s_2,s_3,s_4 \rangle$ and the blue part is in ${}^I W$.  The corresponding Shi vector of $t$ is given by
\begin{figure}[h!]
\begin{center}
\begin{tikzpicture} 
\node at (0,0) {$-1$} ;
\node at (1,0) {$0$} ;
\node at (2,0) {$0$} ;
\node at (3,0) {$0$} ;
\node at (0.5,0.5 ) {$-1$} ;
\node at (1.5, 0.5 ) {$0$} ;
\node at (2.5, 0.5) {$0$} ;
\node at (0.5,1 ) {$-1$} ;
\node at (1.5, 1 ) {$0$} ;
\node at (2.5, 1) {$0$} ;
\node at (1, 1.5 ) {$-1$} ;
\node at (2, 1.5) {$0$} ;
\node at (1, 2 ) {$-1$} ;
\node at (2, 2) {$0$} ;
\node at (1.5, 2.5) {$-1$} ;
\node at (1.5, 3) {$-1$} ;
\node at (5,0) {$e_{12}$} ;
\node at (6,0) {$e_{23}$} ;
\node at (7,0) {$e_{34}$} ;
\node at (8,0) {$2e_4$} ;
\node at (5.5,0.5 ) {$e_{13}$} ;
\node at (6.5, 0.5 ) {$e_{24}$} ;
\node at (7.5, 0.5) {$e^{34}$} ;
\node at (5.5,1 ) {$e_{14}$} ;
\node at (6.5, 1 ) {$e^{24}$} ;
\node at (7.5, 1) {$2e_3$} ;
\node at (6, 1.5 ) {$e^{14}$} ;
\node at (7, 1.5) {$e^{23}$} ;
\node at (6, 2 ) {$e^{13}$} ;
\node at (7, 2) {$2e_2$} ;
\node at (6.5, 2.5) {$e^{12}$} ;
\node at (6.5, 3) {$2e_1$} ;
\end{tikzpicture}
\end{center}
\caption{The Shi vector of $t$ in $C_4$. The right part gives the coordinates of the Shi vector.}
\end{figure}
and the $I$-decomposition of $t$ is given by
\begin{figure}[h!]
\begin{center}
\begin{tikzpicture} 
\node at (0,0) {$-1$} ;
\node at (1,0) {$0$} ;
\node at (2,0) {$0$} ;
\node at (3,0) {$0$} ;
\node at (0.5,0.5 ) {$-1$} ;
\node at (1.5, 0.5 ) {$0$} ;
\node at (2.5, 0.5) {$0$} ;
\node at (0.5,1 ) {$-1$} ;
\node at (1.5, 1 ) {$0$} ;
\node at (2.5, 1) {$0$} ;
\node at (1, 1.5 ) {$-1$} ;
\node at (2, 1.5) {$0$} ;
\node at (1, 2 ) {$-1$} ;
\node at (2, 2) {$0$} ;
\node at (1.5, 2.5) {$-1$} ;
\node at (1.5, 3) {$-1$} ;
\node at (4,1) {$=$} ;
\node at (5,0) {$0$} ;
\node at (6,0) {$\textcolor{purple!80}{0}$} ;
\node at (7,0) {$\textcolor{purple!80}{0}$} ;
\node at (8,0) {$\textcolor{purple!80}{0}$} ;
\node at (5.5,0.5 ) {$0$} ;
\node at (6.5, 0.5 ) {$\textcolor{purple!80}{0}$} ;
\node at (7.5, 0.5) {$\textcolor{purple!80}{0}$} ;
\node at (5.5,1 ) {$0$} ;
\node at (6.5, 1 ) {$\textcolor{purple!80}{0}$} ;
\node at (7.5, 1) {$\textcolor{purple!80}{0}$} ;
\node at (6, 1.5 ) {$0$} ;
\node at (7, 1.5) {$\textcolor{purple!80}{0}$} ;
\node at (6, 2 ) {$0$} ;
\node at (7, 2) {$\textcolor{purple!80}{0}$} ;
\node at (6.5, 2.5) {$0$} ;
\node at (6.5, 3) {$0$} ;
\node at (9,1) {$\bullet$} ;
\node at (10,0) {$\textcolor{teal}{-1}$} ;
\node at (11,0) {$0$} ;
\node at (12,0) {$0$} ;
\node at (13,0) {$0$} ;
\node at (10.5,0.5 ) {$\textcolor{teal}{-1}$} ;
\node at (11.5, 0.5 ) {$0$} ;
\node at (12.5, 0.5) {$0$} ;
\node at (10.5,1 ) {$\textcolor{teal}{-1}$} ;
\node at (11.5, 1 ) {$0$} ;
\node at (12.5, 1) {$0$} ;
\node at (11, 1.5 ) {$\textcolor{teal}{-1}$} ;
\node at (12, 1.5) {$0$} ;
\node at (11, 2 ) {$\textcolor{teal}{-1}$} ;
\node at (12, 2) {$0$} ;
\node at (11.5, 2.5) {$\textcolor{teal}{-1}$} ;
\node at (11.5, 3) {$\textcolor{teal}{-1}$} ;
\end{tikzpicture}
\end{center}
\caption{The $I$-decomposition of $t$ in terms of Shi vectors in $C_4$. }
\label{ShiC4}
\end{figure}
\end{Exa}

\begin{Lem}[Type $D_n$]\label{Lemma D}
Let t be the reflection of the highest root $\tal$ and let $I = \{s_2,\dots,s_n\}$.  Then one has the following reduced expressions
and properties.
\begin{enumerate}
 \item  $t = s_2s_3\dots s_{n-2}s_ns_{n-1}s_{n-2}\dots s_3s_2 \cdot s_1 \cdot s_2s_3\dots s_{n-2}s_{n-1}s_ns_{n-2}\dots s_3s_2$.
\item $t_I = s_2s_3\dots s_{n-2}s_ns_{n-1}s_{n-2}\dots s_3s_2$ and ${}^It =  s_1t_I$.
\item $N(t) = \{e_{13},e_{14},\dots,e_{1n},e_{23},e_{24}\dots,e_{2n}\} \cup \{ e^{12}, e^{13}\dots,e^{1n}, e^{23},e^{24},\dots, e^{2n} \}$.
\item $\sL(t,I)  = 2n^2-4n+1$.
\end{enumerate}
\end{Lem}

\begin{proof} In type $D_n$ the highest root is $\tal=e^{12}$. 
\begin{enumerate}
 \item  Once again, we need to show that $s_2s_3\dots s_{n-2}s_ns_{n-1}\dots s_3s_2(\al_1) = \tal$ with $\tal  = e^{12}$ and $\al_1 = e_{12}$. Roughly,  $s_2(e_{12}) = e_{13}$ and by induction $s_{n-1}\dots s_3s_2(e_{12}) = e_{1n}$. Then, since $\al_n = e^{n-1,n}$, it follows that $s_n(e_{1n})=e^{1,n-1}$ and $s_{n-2}(e^{1,n-1}) = e^{1,n-2}$.  By induction it follows that $s_2s_3\dots s_{n-3}(e^{1,n-2}) = e^{12}$.

\item The fact that $t_I \in W_I$ is clear.  Concerning ${}^It$,  it is easy to see that $s_2s_3\dots s_{n-1}s_ns_{n-2}\dots s_3s_2 = t_I$ since $s_{n-2}s_{n-1}s_ns_{n-2}$ = $s_{n-2}s_{n}s_{n-1}s_{n-2}$ (the generators $s_{n-1}$ and $s_{n}$ commute).  Therefore we have $t = t_Is_1t_I$. The way to show that $s_1t_I$ belongs to ${}^IW$ is exactly the same as in types $A,B,C$ and the details are left to the reader.  Hence ${}^It = s_1t_I$. It is also clear that these two expressions are reduced.  Thus the given expression of $t$ is reduced.

\item Since $t = s^{12}$,  the set of its fixed roots is exactly the set of roots having no label 1 and 2 and the special root $e_{12}$ since we have already seen that $\left( e_{12} \mid e^{12} \right) = 0$ 
and thus $s^{12}(e_{12}) = e_{12}$. 
Therefore, by \Cref{highest_reflection} one has
$$
N(t) = \overline{\fix_t } = \{e_{13},e_{14},\dots,e_{1n},e_{23},e_{24}\dots,e_{2n}\} \cup \{ e^{12}, e^{13}\dots,e^{1n}, e^{23},e^{24},\dots, e^{2n} \},
$$
\item From Point (2) it follows that 
$
N(t) \cap \overline{\Phi_I} =  \{e_{13},e_{14},\dots,e_{1n}\} \cup \{ e^{12}, e^{13}\dots,e^{1n}\}.
$
In type $D_n$ the height of the roots is given by $\h(e_{ij}) = j-i$,  $\h(e^{in}) = n-i$ and $\h(e^{ij}) = 2n-(i+j)$ for any $1 \leq i < j <n$. Therefore we have
\begin{align*}
\sL_{\overline{\Phi_I}}(t) & = \sum\limits_{\al \in N(t) \cap \overline{\Phi_I} } \h(\al) = \sum\limits_{j=3}^{n}\h(e_{1j})  + \sum\limits_{j=2}^{n-1}\h(e^{1j})  + \h(e^{1n}) \\
& = \sum\limits_{j=3}^{n}(j-1)  + \sum\limits_{j=2}^{n-1}(2n-(1+j)) + n-1 \\
& =  \sum\limits_{j=2}^{n-1}j  +  (n-2)(2n-1)  -  \sum\limits_{j=2}^{n-1}j + n-1 \\
& =  (n-2)(2n-1)+n-1\\
& = 2n^2-4n+1.
\end{align*}
\end{enumerate}
\end{proof}


\section{Image of the atomic length}\label{sec_surj_AL}

In this section, we consider again only
finite Weyl groups, which we will again simply denote by $W$.
We prove that the image of the atomic length is the integer interval $\llbracket 0, \sL(w_0)\rrbracket$,
using two independent methods providing different insights.
Recall that we have proved in \Cref{atomic_longest_element} that the element of $W$ 
realising the largest atomic length is $w_0$, the longest element of $W$.

\begin{Lem} \label{lemma_bornes}
The integer $\sL(w_0)$ is given by the following formulas.
\begin{enumerate}
\item For the classical types $X_n$, we have
$$
\sL (w_0)=
\left\{
\begin{array}{ll}
\frac{1}{6}n(n+1)(n+2) = {n+2 \choose 3}
&
\quad \text{if $X=A$}
\\
\frac{1}{6}n(n+1)(4n-1)
&
\quad \text{if $X\in\{B,C\}$}
\\
\frac{1}{3}n(n-1)(2n-1)
&
\quad \text{if $X=D$.}
\end{array}
\right.
$$
\item For the exceptional types, we have the following values.
$$
\begin{array}{@{\hskip 0pt}l@{\hskip 20pt}  @{}l@{}}
\toprule
\text{Type}
& 
\sL(w_0)
\\
\midrule
E_6
& 
156
\\
E_7
&
399
\\
E_8
&
1240
\\
F_4
&
110
\\
G_2
&
16
\\
\bottomrule
\end{array}
$$
\end{enumerate}
\end{Lem}

\begin{proof}
By \Cref{lem_rho_inv}, we have
$\rho-w_0(\rho) = \sum_{\al\in N(w_0)}\al = \sum_{\al\in\Phi^+}\al=2\rho$.
Therefore, $\sL(w_0) = \langle 2\rho, \rho^\vee\rangle= 2 \langle \rho, \rho^\vee\rangle$.
We can compute this inner product by expressing $\rho$ and $\rho^\vee$
as a vector in the canonical basis, and then doing componentwise multiplication.
We give these vectors in the following table, which we have recovered from \cite{BOURB},
and the reader can check that one gets the desired values of $\sL(w_0)$.
Recall that $\rho^\vee=\rho$ in simply-laced type.
$$
\begin{array}{@{}l@{\hskip 20pt} @{}l@{\hskip 20pt} @{}l@{}}
\toprule
\text{Type}
& \rho
&\rho^\vee
\\ 
\midrule
A_n
&
(\frac{n}{2}, \frac{n}{2}-1, \ldots, -\frac{n}{2}+1, -\frac{n}{2})
&
\text{same as $\rho$}
\\ 
B_n
&
(n-\frac{1}{2}, n-\frac{3}{2}, \ldots, \frac{1}{2})
&
(n, n-1, \ldots, 1 )
\\ 
C_n
&
(n, n-1, \ldots, 1 )
&
(n-\frac{1}{2}, n-\frac{3}{2}, \ldots, \frac{1}{2})
\\ 
D_n
&
(n-1, n-2, \ldots, 0)
&
\text{same as $\rho$}
\\ 
E_6
&
( 0, 1, 2, 3, 4, -4, -4, 4 ) 
&
\text{same as $\rho$}
\\
E_7
&
( 0, 1, 2, 3, 4, 5, -\frac{17}{2}, \frac{17}{2}) 
&
\text{same as $\rho$}
\\
E_8
&
( 0, 1, 2, 3, 4, 5, 6, 23 )
&
\text{same as $\rho$}
\\
F_4
&
(\frac{11}{2}, \frac{5}{2}, \frac{3}{2}, \frac{1}{2})
&
(8,3,2,1)
\\
G_2
&
(-1, -2, 3)
&
(-\frac{1}{3}, -\frac{4}{3}, \frac{5}{3})
\\
\bottomrule
\end{array}
$$
\end{proof}

The formula of \Cref{lemma_bornes} depends on the chosen Dynkin type.
In fact, finding a uniform formula for $\sL(w_0)=2\langle \rho, \rho^\vee\rangle$ seems complicated,
as opposed to the "strange formula" of Freudenthal and de Vries \cite{FdV, Burns2000} which states that
$2(\rho | \rho) = ng(h+1)/6$, where $h$ and $g$ are the Coxeter and dual Coxeter number respectively.
Nevertheless, we are able to give in the following proposition an expression of $\sL(w_0)$
which is independent of the type.
In order to do this, recall the notation of \Cref{sec_affine_WG}, and consider the fundamental polytope 
$\mathcal{P}= \left\{\sum_{i=1}^nc_i\omega_i~|~0 \leq c_i \leq 1 \right\}$.
There is a unique alcove $A_{w_{\mathcal{P}}} \subseteq \cP$ whose labelling element
$w_{\mathcal{P}}\in W_a$ has maximal length
\cite{LP2018, NC21lattice}.

\begin{Prop}\
\begin{enumerate}
\item For all $\al\in\Phi^+$, we have $k(w_{\mathcal{P}},\al) = \h(\al)-1$.  
\item We have $\sL(w_0) = \ell(w_0) + \ell(w_{\mathcal{P}})$.
\end{enumerate}
\end{Prop}

\begin{proof}\
\begin{enumerate}
\item 
For $w\in W_a$, we have $A_w\subset\cP$ if and only if $k(w,\al_i) = 0$ for all $1\leq i\leq n$.
In particular, the vector
$(\h(\alpha)-1)_{\alpha \in \Phi^+}$, which is a Shi vector by \Cref{Shi ineq}, gives an element of $W_a$ labelling
an alcove of $\cP$.
Let us show that this element coincides with $w_\cP$.
By \cite[Theorem 4.1]{NC1}, we know that for all $w\in W_a$ and for all $\alpha \in \Phi^+$, 
there exists $P_{\alpha} \in \mathbb{\N}[X_1,\ldots, X_n]$ with no constant term 
and  $\lambda_{\alpha}(w) \in \llbracket 0, \h(\alpha)-1\rrbracket$ such that 
$k(w_{\mathcal{P}},\alpha) = P_{\alpha}(k(w,\alpha_1), \ldots, k(w,\alpha_n)) + \lambda_{\alpha}(w)$.  
In particular, for $A_w \subset \cP$, the characterisation above implies that
$P_{\alpha}(k(w,\alpha_1), \ldots, k(w,\alpha_n)) = 0$, and so $k(w,\alpha) = \lambda_{\alpha}(w)$. 
Moreover, by \cite[Proposition 4.3]{JYS1}, we have $\ell(w) = \sum_{\alpha \in \Phi^+} |k(w, \alpha)|$ for all $w\in W_a$.
Therefore, since $w_\cP$ realises the maximal length within $\cP$, 
each of its Shi coefficients must be maximal, that is,
$k(w_{\mathcal{P}},\alpha) = \lambda_{\alpha}(w_{\mathcal{P}}) = \h(\alpha)-1$ for all $\alpha \in \Phi^+$. 
\item 
In particular, we have $\ell(w_{\mathcal{P}}) = \sum_{\alpha \in \Phi^+} |k(w_{\mathcal{P}}, \alpha)| = \sum_{\alpha \in \Phi^+} (\h(\al)-1)$ by (1).
Now
$$
\sL(w_0) = \sum\limits_{\alpha \in \Phi^+}\h(\alpha) =   |\Phi^+| + \sum\limits_{\alpha \in \Phi^+}(\h(\alpha) -1)= \ell(w_0) +\ell(w_{\mathcal{P}}).
$$
\end{enumerate}
\end{proof}

We are ready to prove the central result of this section.

\begin{Th}\label{thm_surj} 
Let $W$ be a finite Weyl group of rank $n$.
The map $\sL : W \to \lbra 0, \sL(w_0) \rbra$ is surjective
if and only if $n\neq 2$.
\end{Th}

\begin{proof}
This is achieved in several steps.
\begin{enumerate}
\item Consider first the small rank cases. 
In type $A_1$, we have $\sL(W)=\{0,1\}$,
and for $n=2$ we get the following images by direct computation.
$$
\begin{array}{@{}l@{\hskip 20pt} @{}l@{}}
\toprule
\text{Type}
& \sL(W)
\\ 
\midrule
A_2
&
\{0,1,3,4\}
\\ 
B_2 \text{ or } C_2
&
\{0,1,3,4,6,7\}
\\ 
G_2
&
\{0, 1, 3, 5, 8, 11, 13, 15, 16\}
\\ 
\bottomrule
\end{array}
$$
\item
Consider now the remaining exceptional types $E_6, E_7, E_8$ and $F_4$. 
One can use the computer algebra softwares SageMath and GAP to obtain $\sL(W)$.
In all cases, we have obtained $\sL(W) = \lbra 0, \sL(w_0) \rbra$.
\item
Finally, assume that we are in classical type $X_n$, $X\in\{A, B, C, D\}$.
We will prove the result by induction on $n$.
To this end, it will be convenient to make the dependence on $n$ appear clearly by using the notation
$W=W(X_n)$, $b_n = \sL(w_0)$
and $\sL_n = \sL : W(X_n) \to \lbra 0, b_n \rbra$.

One checks by a direct computation that $\sL_n$ surjects onto $\lbra 0, b_n \rbra$ for
$n=3$ if $X=A$, for $n\in\{3,4\}$ if $X=B,C$ and for $n\in\{4,5\}$ for $X=D$.

Assume now that  $\sL_n$  surjects onto $\lbra 0, b_n \rbra$ for a fixed $n$, where $n\geq3$ if $X=A$, respectively $n\geq4$ if $X=B,C$, respectively $n\geq 5$ for $X=D$.
By \Cref{AL_symmetry_w0}, it suffices to show that
for all $N\in \lbra 0, b_{n+1}/2\rbra$, there exists $w\in W(X_{n+1})$ such that $\sL_{n+1}(w) = N$
to ensure that $\sL_{n+1} : W(X_{n+1}) \to \lbra 0, b_{n+1}\rbra$ is surjective.
In order to do that, let $I=\{s_2,\ldots, s_{n+1}\}$ and
consider the standard parabolic subgroup $W_I \leq W(X_{n+1})$.
By \Cref{Atomic length parabolic}, we are ensured that
$$\sL (W_I)= \sL_I (W_I) = \lbra 0, b_n \rbra$$ 
where the second identity is the induction hypothesis.
Thus, to conclude, it suffices to prove that 
$$\frac{b_{n+1}}{2} \leq b_n.$$
We show this by considering the different cases of \Cref{lemma_bornes}.
\begin{itemize}
\item If $X=A$, then $b_n=\frac{1}{6}n(n+1)(n+2)$. Therefore 
$$
\begin{array}{rcl}
b_n - \frac{b_{n+1}}{2}
& = &
\frac{1}{12} \left( 2n(n+1)(n+2) - (n+1)(n+2)(n+3) \right)
\\
& = &
\frac{1}{12} (n+1)(n+2)(n-3)
\\
& \geq &
0 
\end{array}
$$
since we have assumed $n\geq 3$.
\item If $X\in\{B,C\}$, then $b_n=\frac{1}{6}n(n+1)(4n-1)$. Therefore 
$$
\begin{array}{rcl}
b_n - \frac{b_{n+1}}{2}
& = &
\frac{1}{12} \left( 2n(n+1)(4n-1) - (n+1)(n+2)(4n+3) \right)
\\
& = &
\frac{1}{12} (n+1)(4n^2-13n-6)
\\
& \geq &
0 
\end{array}
$$
since the largest root of this polynomial is 
$\frac{1}{8}(13+\sqrt{265})\in ( 3, 4)$,
and we have assumed $n\geq 4$.
\item If $X=D$, then $b_n=\frac{1}{3}n(n-1)(2n-1)$. Therefore 
$$
\begin{array}{rcl}
b_n - \frac{b_{n+1}}{2}
& = &
\frac{1}{6} \left( 2n(n-1)(2n-1) - (n+1)n(2n+1) \right)
\\
& = &
\frac{1}{6}  n(2n^2-9n+1) 
\\
& \geq &
0 
\end{array}
$$
since the largest root of this polynomial
is $\frac{1}{4}(9+\sqrt{73})\in ( 4, 5)$,
and we have assumed $n\geq 5$.
\end{itemize}
\end{enumerate}
\end{proof}

We now give an alternative proof of \Cref{thm_surj}
based on Susanfe theory.

\begin{proof}[Proof of \Cref{thm_surj}]
We have already explained in the previous proof how to treat the exceptional cases using computer algebra.
For the classical cases, we will again proceed by induction on the rank $n$.

\medskip

Fix $X\in\{ A, B, C, D\}$.
The small rank cases $n\leq 3$ if $X=A$,
$n\leq 4$ if $X\in\{B,C\}$ and
$n\leq 5$ if $X =D$
can be checked individually.
Assume the following induction hypothesis:  
$\sL(W(X_n))=\lbra 0, b_n\rbra$
where $n+1$ is fixed, 
different from the values above.
Here, $b_n$ denotes again the image of $w_0$ (the longest element of $W(X_n)$) 
by the function $\sL:W(X_n)\to \N$.
Recall that we know that $w_0$ indeed realises the maximum by \Cref{atomic_longest_element},
and that we gave the formula for $b_n$ in \Cref{lemma_bornes}.

Take now $W = W(X_{n+1})$.
We will use the same notation as in \Cref{sec_susanfe}, namely 
$I = \{s_2,s_3,\dots, s_{n+1}\}$.
Moreover, we denote by $t$ 
the reflection corresponding the highest root $\tal$
(and $t'= {}^I t s_1$ if $X=B$).
Finally, we let $J=tIt$ so that $tW_J=W_It$.
The idea of the inductive step is to consider the union
$\sL(W_I) \cup \sL(W_It)$
(respectively $\sL(W_I) \cup \sL(W_It')$ if $X=B$).
Let us look at these two image sets separately first.
\begin{itemize}
 \item By definition of $I$, the subgroup $W_I$ verifies
$W_I\simeq W(X_n)$, so we have 
$\sL_{I}(W_I)=\lbra 0, b_n\rbra$ by induction hypothesis.  
Moreover, since $W_I$ is a standard parabolic subgroup of $W$,  by Proposition \ref{Atomic length parabolic} we have  $\sL_{I}(W_I) = \sL(W_I)$ and thus 
\begin{equation}\label{LWI}
\sL(W_I) = \lbra 0, b_n\rbra.
\end{equation}
\item For the second set $\sL(W_I t)$, let us use Corollary \ref{corollary principal} with $B=J$.
This yields, for all $w \in W_J$, 
\begin{equation}\label{Ltw}
\sL(tw) = \sL_{I}((tw)_I) + \sL(t,I).
\end{equation}
Let us analyse the right-hand side.
First, by \Cref{Proposition_Utopic}, we know that $\{(tw)_I \mid w \in W_J\} = W_I$,
and by induction hypothesis, $\sL_I(W_I)=\lbra 0, b_n\rbra$.
Second, the term $\sL(t,I)$ is a non-negative integer depending only on $X_{n+1}$, which we denote by $K_{n+1}$.
The formulas for $K_{n+1}$ for the different types are given in the lemmas of \Cref{sec_susanfe}.
Therefore, \Cref{Ltw} implies
$\sL(t W_J) = \sL_I(W_I) + \sL(t,I) = \lbra 0, b_n\rbra + K_{n+1} = \lbra K_{n+1}, b_n+K_{n+1} \rbra$,
which we can rewrite, since $tW_J=W_It$, as
\begin{equation}\label{LWIt}
\sL(W_I t) =  \lbra K_{n+1}, b_n+K_{n+1} \rbra.
\end{equation}
\end{itemize}
Therefore, combining \Cref{LWI} and \Cref{LWIt}, we have proved that
\begin{equation}\label{unionL}
\sL(W_I) \cup \sL(W_It) = \lbra 0, b_n\rbra \cup \lbra K_{n+1}, b_n+K_{n+1}\rbra.
\end{equation}
Obviously, $\sL(W_I) \cup \sL(W_It)\subseteq \sL(W)$.
In order to show that $\sL(W)=\lbra 0, b_{n+1}\rbra$, we will prove that
\begin{enumerate}
\item $\sL(W_I) \cup \sL(W_It)$ is in fact an interval, by showing for each Dynkin type that $K_{n+1}\leq b_n$,
\item the maximum of this inverval, that is $b_n+K_{n+1}$, equals $b_{n+1}$ 
(respectively $b_n+K_{n+1} = b_{n+1}-1$ if $X=D$, which will actually suffice).
\end{enumerate}
Note that in type $B_n$, we need to replace $t$ by $t'$ in the above reasoning.

\begin{itemize}
\item Assume $X=A$. We have $b_n = \binom{n+1}{3}=\frac{(n+1)n(n-1)}{2}$ and $K_{n+1}=\binom{n+1}{2}$.
Since $n+1\geq 4$ by assumption, we have
$ \binom{n+1}{2} < \binom{n+1}{3}$, hence $K_{n+1}<b_n$, proving (1).

Finally, we have the well-known formula 
$\binom{n+1}{3} + \binom{n+1}{2} = \binom{n+1}{3}$, that is, $b_n+K_{n+1} = b_{n+1}$, proving (2).

\item Assume $X\in \{B,C\}$.
We have $b_n = \frac{n(n+1)(4n-1)}{6}$ and $K_{n+1} = 2(n+1)^2-(n+1) = 2n^2+3n+1$.

We compute
\begin{align*}
\frac{n(n+1)(4n-1)}{6} - (2n^2+3n+1) 
&= 
\frac{1}{6}(n+1)(4n^2-13n-6).
\end{align*}
We recover twice the polynomial of the previous proof, and we have seen that its largest root $y$ verifies
$3<y<4$. Since we have assumed that $n\geq 4$, we are ensured that $K_{n+1}<b_n$, proving (1).

Finally, we compute 
$$b_n+K_{n+1} =  \frac{n(n+1)(4n-1)}{6} + 2(n+1)^2-(n+1) = \frac{(n+1)(n+2)(4n+3)}{6}=b_{n+1},$$
proving (2).

\item Assume $X=D$.
We have $b_n = \frac{n(n-1)(2n-1)}{3}$ and $K_{n+1} = 2(n+1)^2-4(n+1)+1 = 2n^2-1$.

We compute
\begin{align*}
\frac{n(n-1)(2n-1)}{3} -(2n^2-1)
&= 
\frac{1}{3}(2n+1)(n^2-5n+3).
\end{align*}
The largest root of this polynomial is
$\frac{1}{2}(5+\sqrt{13})\in (4,5)$.
Since we have assumed that $n\geq 5$, we are ensured that $K_{n+1}<b_n$, proving (1).

Finally, we compute 
$$b_n+K_{n+1} =  \frac{n(n-1)(2n-1)}{3} + 2n^2-1 = \frac{(n + 1)n(2 n + 1)}{3} - 1 = b_{n+1}-1.$$
We are missing the value $b_{n+1}$ by this method, but fortunately, we know that $b_{n+1}=\sL(w_0)\in\sL(W)$,
proving (2).
\end{itemize}
\end{proof} 

\begin{Rem}
In the case of type $A_n$, 
recall that the atomic length coincides with the statistic
$\mathsf{invsum}$, see \Cref{rem_inv_type_A}.
In fact, the previous proof is analogous to the proof of \cite[Section 2]{SackUlfarsson2011}
recalled in \Cref{sec_entropy}.
We see that the results on
the entropy of permutations are just a particular case of a more general phenomenon.
Moreover, in \cite[Sections 3 and 4]{SackUlfarsson2011}, the authors are able to obtain several interesting properties of the
generating function of $\mathsf{ninvsum}$ (including product formulas), involving the $q$-analogues of binomial coefficients.
On the other hand, it is well-known that the generating function of the Coxeter length can be described by certain $q$-binomial
coefficients.
It would be interesting to study the generating function for the atomic length and look for product formulas.
\end{Rem}

It is natural to look for a generalisation of \Cref{thm_surj} by asking,
for a fixed $\la\in P^+$, whether $\sL_\la:W\to \lbra 0, \sL_\la(w_0)\rbra$ is surjective.
Even more ambitious would be the classification of weights $\la\in P^+$ such that $\sL_\la$ surjects onto $\lbra 0, \sL_\la(w_0)\rbra$.
This motivates the following definition.

\begin{Def}\label{def_ideal}
Let $W$ be an affine (respectively finite) Weyl group.
An element $\la\in P^+$ is called \textit{$W$-ideal}
if $\sL_\la:W\to \N$ (respectively $\sL_\la:W\to\lbra0, \sL_\la(w_0)\rbra$)
is surjective.
\end{Def}

If there is no ambiguity about the group $W$, we will simply say that $\la\in P^+$ is ideal.
We can reformulate \Cref{thm_surj} by saying that $\overline{\rho}$ is always ideal except in rank $2$.
We will see in the next section two further examples of ideal weights:
\begin{itemize}
\item in finite types, the minuscule weights, this will be \Cref{minusc},
\item in untwisted affine type $A$, the weight $\la=\La_0$, this will be \Cref{GO}, 
a (very much non trivial) theorem by Granville and Ono.
\end{itemize}
In fact, in finite types, it is not hard to see that we have the following necessary condition.

\begin{Prop}\label{ideal_necessary}
Let $W$ be a finite Weyl group and let $\la=\sum_{i=1}^n m_i\om_i \in P^+$.
If $\la$ is ideal, then there exists $1\leq i\leq n$ such that $m_i=1$.
\end{Prop}

\begin{proof}
By contradiction, assume that $m_i\geq 2$ for all $1\leq i\leq n$.
We can use \Cref{th_lambda_atomic_length} to express the $\la$-atomic length in terms
of the $\la$-inversion set.
Let $w\in W$ and $\underline{w}\in\Red(w)$.
Since the elements appearing in $N_\la(\underline{w})$ are of the form 
$\beta = m_i\underline{w}_{i,k}(\al_i)$
for some $\underline{w}_{i,k}(\al_i)\in Q$, the condition $m_i\geq2$ implies that $\h(\beta)\geq2$.
In particular, it is impossible to find $w\in W$ verifying $\sL_\la(w)=1$,
so $\sL_\la : W \to \lbra 0,\sL_\la(w_0)\rbra$ is not surjective.
\end{proof}

In the search for ideal weights, the previous proposition rules out many possibilities.
However, among the remaining weights, it seems complicated to
classify those that are ideal, as illustrated in the following example.

\begin{Exa}
Let $W$ be the Weyl group of type $C_3$.
Consider dominants weights $\la= m_1\om_1 + m_2\om_2+m_3\om_3$.
\begin{itemize}
\item If $(m_1, m_2, m_3) = (2,1,1)$ then $\sL_\la(W)=\lbra 0, 27\rbra$ and $\la$ is ideal.
\item If $(m_1, m_2, m_3) = (1,2,1)$ then 
$$\sL_\la(W)= \{0, 1, 2, 4, 5, 6, 7, 8, 9, 10, 11, 13, 14, 15, 16, 17, 19, 20, 21, 22, 23, 24, 25, 26, 28, 29, 30\}$$
and $\la$ is not ideal.
\item If $(m_1, m_2, m_3) = (1,1,2)$ then 
$$\sL_\la(W)= \{
0, 1, 2, 3, 4, 6, 7, 8, 9, 10, 11, 13, 14, 15, 16, 17, 18, 20, 21, 22, 23, 24, 25, 27, 28, 29, 30, 31
\}$$
and $\la$ is not ideal.
\end{itemize}
\end{Exa}


\section{Links with crystal theory}\label{sec_cryst}

\subsection{Atomic length and crystal depth}

Consider the affine Kac-Moody algebra $\fg$ introduced in \Cref{sec_affine_WG}.
For each $\la\in P^+$, one can construct the corresponding irreducible highest weight module $V(\la)$ of $\fg$
already introduced in \Cref{AL_integer}.
Since this representation is \textit{integrable}, it comes equipped with a  \textit{crystal}, denoted $B(\la)$, see \cite{HongKang2002}.
This is an oriented colored graph 
whose structure mirrors the algebraic structure of $V(\la)$.
For finite classical types, there are explicit constructions of $B(\la)$ relying on tableau combinatorics, see \cite[Chapters 7 and 8]{HongKang2002} for details.
For affine type $A$, there are explicit constructions of $B(\la)$ given in terms of multipartitions/abaci combinatorics, see \cite[Chapter 6]{GeckJacon2011}.
The crystal $B(\la)$ is connected and has a unique source vertex $b_\la$ called its \textit{highest weight vertex}.
Moreover, there is a weight function $\wt : B(\la) \to P$ defined on the vertices of $B(\la)$ and determined by the properties
\begin{align}
& \wt(b_\la)  =\la, 
\label{w1}
\\
& \wt(b) = \wt(b')-\al_i \text{ if there is an arrow } b' \overset{i}{\lra} b \text{ in }B(\la).
\label{w2}
\end{align}
Indeed, each vertex $b\in B(\la)$ can be obtained from $b_\la$ by following a sequence of arrows,
so the above formulas suffice to determine $\wt$ on the entire $B(\la)$.
Moreover, one can check that for all $1 \leq i \leq n$ and for all $b\in B(\la)$,
\begin{equation}\label{w3}
\langle \wt(b), \al_i^\vee\rangle = \varphi_i(b) - \eps_i(b),
\end{equation}
where $\varphi_i(b)$ (respectively $\eps_i(b)$) is the length of the string of $i$-arrows outgoing (respectively incoming) at vertex $b$.
One recovers for instance the dimensions of the weight spaces of $V(\la)$ by counting
the number of vertices with the same weight in $B(\la)$. Clearly, these all appear at the same depth in $B(\la)$.

\medskip

Now, the Weyl group $W$ acts on $B(\la)$ in a particularly simple way.
For $0\leq i\leq n$ fixed, removing all $j$-arrows with $j\neq i$ as well as all vertices without incoming or outgoing $i$-arrows
yields a disjoint union of $i$-strings.
Then the generator $s_i\in W$
acts on vertices by reflecting with respect to the middle of the $i$-string \cite[Theorem 11.14]{BumpSchilling2017}.
In the rest of this section, we will be particularly interested in the orbit of the highest weight vertex $b_\la$,
which we will denote $\cO(\la)$.

\begin{Exa}\label{exa_crystal_S3_rho}
The crystal $B(\la)$ in type $A_2$ for $\la=\rho=\om_1+\om_2$
can be constructed using all semistandard tableaux of shape $(2,1)$.
This gives the following graph.

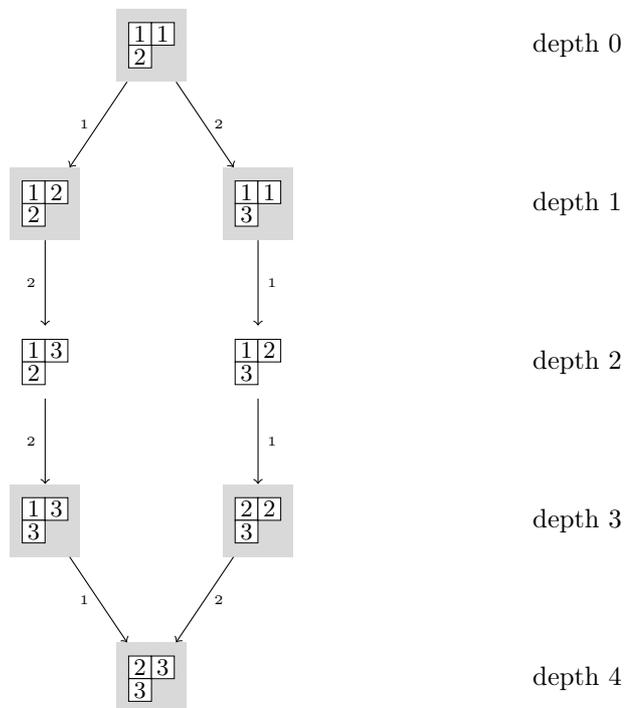
\begin{figure}[H] \centering
\begin{tikzpicture}[scale=0.7]
\Yboxdim{3mm}

\node (b) at (0,8) [fill=gray!30] {
{${ \small \young(11,2)  }$}
};

\node (c1) at (-2,5)[fill=gray!30] {
{${ \small  \young(12,2)  }$}
}; 
\node (c2) at (2,5)[fill=gray!30] {
{${ \small  \young(11,3)  }$}
};

\node (d1) at (-2,2) {
{${ \small  \young(13,2)  }$}
}; 
\node (d2) at (2,2) {
{${ \small  \young(12,3)  }$}
}; 

\node (e1) at (-2,-1)[fill=gray!30] {
{${ \small  \young(13,3)  }$}
}; 
\node (e2) at (2,-1)[fill=gray!30] {
{${ \small  \young(22,3)  }$}
};

\node (f1) at (0,-4)[fill=gray!30] {
{${ \small  \young(23,3)  }$}
};

\node (x0) at (8,8) {
{depth $0$}
}; 
\node (x1) at (8,5) {
{depth $1$}
}; 
\node (x1) at (8,2) {
{depth $2$}
}; 
\node (x1) at (8,-1) {
{depth $3$}
}; 
\node (x1) at (8,-4) {
{depth $4$}
};

\node (phantom) at (-8,8) {
{\phantom{0}}
}; 
\draw[->] (b) -- node [font=\tiny,left] {1} (c1) ;
\draw[->] (b) -- node [font=\tiny, right] {2} (c2) ;
\draw[->] (c1) -- node [font=\tiny, left] {2} (d1) ;
\draw[->] (c2) -- node [font=\tiny, right] {1} (d2) ;
\draw[->] (d1) -- node [font=\tiny, left] {2} (e1) ;
\draw[->] (d2) -- node [font=\tiny, right] {1} (e2) ;
\draw[->] (e1) -- node [font=\tiny, left] {1} (f1) ;
\draw[->] (e2) -- node [font=\tiny, right] {2} (f1) ;
\end{tikzpicture}
\end{figure}

In this picture, we have highlighted in gray
the elements of $\cO(\la)$.
They are obtained by starting from the highest weight vertex
and reflecting along $i$-strings.
\end{Exa}

\begin{Rem}\label{rem_keys}
In type $A_n$ and for any $\la\in P^+$,
there is a simple rule for computing $\cO(\la)$ recursively:
$s_i$ acts on a tableau $b\in \cO(\la)$ by
changing all possible entries $i$ to $i+1$
(that is, so that the resulting tableau is still semistandard).
From there, one can deduce an explicit description of $\cO(\la)$, namely $\cO(\la)$ consists of those tableaux such that each column contains the one to its right, see  \cite{LS1990}.
For results in other classical types, see \cite{Santos2021} and \cite{JaconLecouvey2020}.
\end{Rem}

This enables us to interpret the atomic length $\sL_\la$ in the context of crystals.

\begin{Prop}\label{crystal_interp}
Let $b \in\cO(\la)$, so that $b=w(b_\la)$ for some $w\in W$.
Then $\sL_\la(w)$ is the depth of $b$ in $B(\la)$.
\end{Prop}

\begin{proof}
By \Cref{w2}, the depth of $b$ in $B(\la)$ is the number of simple roots
substracted to $\wt(b_\la)$ to get $\wt(b)$,
that is, it is the number of simple roots that appear in the decomposition of $\wt(b_\la) - \wt(b)$.
On the other hand, $\sL_\la(w) = \langle \la-w(\la), \rho^\vee\rangle$,
that is, $\sL_\la(w)$ is the number of simple roots
that appear in the decomposition of $\la-w(\la)$.
To conclude, we first use \Cref{w1} which ensures that $\la = \wt(b_\la)$.
Finally, for all $i\in I$ and for all $a\in B(\la)$,
$$
\begin{array}{rcll}
\wt(s_i(a)) 
& = 
& \wt(a) - (\varphi_i(a) - \eps_i(a)) \al_i
& \text{ by definition of the action of $s_i$}
\\
& =
& \wt(a) - \langle \wt(a), \al_i^\vee\rangle \al_i
& \text{ by \Cref{w3}}
\\
& = 
& s_i(\wt(a)).
\end{array}
$$
Therefore $w(\la) = w (\wt(b_\la)) = \wt( w(b_\la)) = \wt(b)$.
\end{proof}

\begin{Exa}
\Cref{crystal_interp} enables us to compare \Cref{exa_atomic_length_S3} and \Cref{exa_crystal_S3_rho}: we recover the values $0,1,3,4$ as depths in the crystal $B(\rho)$.
\end{Exa}

\begin{Rem}
In finite types, the action of $-w_0$ on the simple roots $\al_i$ induces an involution of the crystal
$B(\la)$, known as the \textit{Sch\"utzenberger-Lusztig involution} and denoted by $\eta_\la$.
More precisely, recall the involution
of the Dynkin diagram $\zeta$ induced from $-w_0$ in \Cref{sec_AL}.
Then for each path 
$$b_\la\overset{i_1}{\lra} \cdots \overset{i_r}{\lra} b$$
in $B(\la)$, there exists a path
$$\eta_\la(b_\la)\overset{\zeta(i_1)}{\lra} \cdots \overset{\zeta(i_r)}{\lra} \eta_\la(b)$$
in $B(\la)$,
where $\eta_\la(b_\la)$ is the lowest weight vertex of $B(\la)$, see \cite{Lenart2005} for more details.
In particular, combining this with \Cref{crystal_interp} enables us to recover \Cref{AL_symmetry_w0}.
\end{Rem}

With this intepretation at hand, we are now able to
show that $\sL_\la$  surjects onto $\lbra 0, \sL_\la(w_0) \rbra$ for a whole family of
dominant weights in finite types.
Recall that a weight $\la\in P^+$ is called \textit{minuscule} if the Weyl group $W$ acts transitively 
on the set of weights of the module $V(\la)$.
Equivalently, $\la$ is minuscule if and only if,
for all $\al\in \Phi$, $\langle\la,\al^\vee\rangle\in\{-1,0,1\}$,
see for instance \cite[Chapter VI, Exercices, \S 1, 24]{BOURB}.
We recall the classification of minuscule weights in the table below, found in \cite[Chapter VI, Exercices, \S 4, 15]{BOURB}.
There is no minuscule weight for types $E_8, F_4$, and $G_2$.
$$
\begin{array}{@{}l@{\hskip 20pt} @{}l@{}}
\toprule
\text{Type} 
&
\text{Minuscule weights} 
\\ 
\midrule
A_n
&
\om_i \, , \, 1\leq i\leq n
\\
B_n
&
\om_n
\\
C_n
&
\om_1
\\
D_n
&
\om_1, \om_{n-1}, \om_n
\\
E_6
&
\om_1, \om_6
\\
E_7
&
\om_7
\\
\bottomrule
\end{array}
$$

For the next result, recall \Cref{def_ideal} introducing the notion of ideal weights.

\begin{Th}\label{minusc}
Let $\la\in P^+$ be minuscule. Then 
$\la$ is ideal.
\end{Th}

\begin{proof}
Since $\la$ is minuscule, the set of weights of $V(\la)$
equals $W\la$, the orbit of $\la$ under the action of $W$.
By general crystal theory, this implies that 
$B(\la) = \cO(\la)$,
that is, every vertex in the crystal graph is in $\cO(\la)$.
In particular, there is an element of
$\cO(\la)$ at every given depth in $B(\la)$.
By \Cref{crystal_interp}, this means that
$\sL_\la:W\to\lbra 0, \sL_\la(w_0)\rbra$ is surjective,
that is, $\la$ is ideal.
\end{proof}

\subsection{Affine crystals of type $A$ and the Granville-Ono theorem}

We now mention a fundamental particular case, namely
that of $\la=\La_0$ in type $A_{n}^{(1)}$.
There is a classical realisation of the corresponding
crystal $B(\La_0)$ that dates back to \cite{MisraMiwa1990}
and uses the \textit{$(n+1)$-regular} partitions,
that is, partitions where each part is repeated at most $n$ times. 
More precisely,
there is an arrow $b\overset{i}{\to}b'$ in $B(\La_0)$
if and only if $b'$ is obtained from $b$ by adding its \textit{good $i$-box}, see \cite{LLT1996} for a definition.
This yields the whole crystal graph $B(\La_0)$
by starting from the highest weight vertex $\emptyset$
(the only partition of $0$) and adding good boxes recursively,
thereby yielding all $(n+1)$-regular partitions.
In turn, similarly to the finite case, the elements of $\cO(\La_0)$
are obtained recursively from $\emptyset$
by adding all $i$-boxes at once (for each fixed $0\leq i\leq n$), which corresponds to the action of $s_i$.
The following result is well-known \cite{LLT1996}, \cite{Lascoux2001}.
Recall that a partition is called a \textit{$(n+1)$-core}
if it has no removable rim $(n+1)$-hook.

\begin{Prop}\label{cores_orbit}
Let $W=W(A_n^{(1)})$ and $b\in B(\La_0)$.
Then $b\in \cO(\La_0)$ if and only if $b$ is an $(n+1)$-core.
Moreover, the depth of $b$ in $B(\La_0)$ is the size of $b$ 
(that is, the number of boxes of $b$).
\end{Prop}

Note that the second part of the statement is obvious
since each arrow in $B(\La_0)$ corresponds to adding some box.

\begin{Exa}\label{exa_cores}
Let $n=2$. The crystal $B(\La_0)$ is realised by $3$-regular
partitions as shown in \Cref{crystal}, 
and the orbit $\cO(\La_0)$ consists precisely of the $3$-cores, which we have highlighted in gray.
Looking at depth (or counting boxes), we see that the first values of 
$\sL_{\La_0}$ are $0,1,2,4,5$ and that there is no $3$-core of size $3$.

\begin{figure}\centering
\begin{tikzpicture}[-,scale=0.5, every node/.style={scale=1}]
\node (a) at (0,-2) [fill=gray!30] {
${ \scriptstyle\emptyset }$
};

\node (b1) at (0,-5)[fill=gray!30] {
{${  \Yboxdim{10pt} 
\yng(1)
}$}
}; 

\node (c1) at (-3,-8)[fill=gray!30] {
{${  \Yboxdim{10pt} 
\yng(2)
}$}
}; 
\node (c2) at (3,-8)[fill=gray!30] {
{${  \Yboxdim{10pt} 
\yng(1,1)
}$}
};

\node (d1) at (-3,-12) {
{${  \Yboxdim{10pt} 
\yng(3)
}$}
}; 
\node (d2) at (3,-12) {
{${  \Yboxdim{10pt} 
\yng(2,1)
}$}
};

\node (e1) at (-7,-16) {
{${  \Yboxdim{10pt} 
\yng(4)
}$}
}; 
\node (e2) at (-3,-16) [fill=gray!30]{
{${  \Yboxdim{10pt} 
\yng(3,1)
}$}
}; 
\node (e3) at (3,-16) [fill=gray!30]{
{${  \Yboxdim{10pt} 
\yng(2,1,1)
}$}
}; 
\node (e4) at (7,-16) {
{${  \Yboxdim{10pt} 
\yng(2,2)
}$}
};

\node (f1) at (-11,-20) {
{${  \Yboxdim{10pt} 
\yng(5)
}$}
}; 
\node (f2) at (-5,-20) {
{${  \Yboxdim{10pt} 
\yng(4,1)
}$}
}; 
\node (f3) at (0,-20) [fill=gray!30]{
{${  \Yboxdim{10pt} 
\yng(3,1,1)
}$}
}; 
\node (f4) at (5,-20) {
{${\Yboxdim{10pt} 
\yng(2,2,1)
}$}
}; 
\node (f5) at (11,-20) {
{${  \Yboxdim{10pt} 
\yng(3,2)
}$}
}; 

\node (g) at (0,-22) {
{${\vdots}$}
}; 

\draw[->] (a) -- node [font=\tiny, left] {0} (b1) ;
\draw[->] (b1) -- node [font=\tiny, left] {1} (c1) ;
\draw[->] (b1) -- node [font=\tiny, right] {2} (c2) ;
\draw[->] (c1) -- node [font=\tiny, left] {2} (d1) ;
\draw[->] (c2) -- node [font=\tiny, right] {1} (d2) ;
\draw[->] (d1) -- node [font=\tiny, left] {0} (e1) ;
\draw[->] (d1) -- node [font=\tiny, right] {2} (e2) ;
\draw[->] (d2) -- node [font=\tiny, left] {1} (e3) ;
\draw[->] (d2) -- node [font=\tiny, right] {0} (e4) ;
\draw[->] (e1) -- node [font=\tiny, left] {1} (f1) ;
\draw[->] (e1) -- node [font=\tiny, right] {2} (f2) ;
\draw[->] (e2) -- node [font=\tiny, left] {0} (f2) ;
\draw[->] (e2) -- node [font=\tiny, right] {1} (f3) ;
\draw[->] (e3) -- node [font=\tiny, left] {2} (f3) ;
\draw[->] (e3) -- node [font=\tiny, right] {0} (f4) ;
\draw[->] (e4) -- node [font=\tiny, left] {1} (f4) ;
\draw[->] (e4) -- node [font=\tiny, right] {2} (f5) ;
\end{tikzpicture}
\caption{The beginning of the 
crystal $B(\Lambda_0)$ in type $A_2^{(1)}$. 
The shaded vertices correspond to $\cO(\La_0)$,
which consists of the $3$-core partitions.}
\label{crystal}
\end{figure}
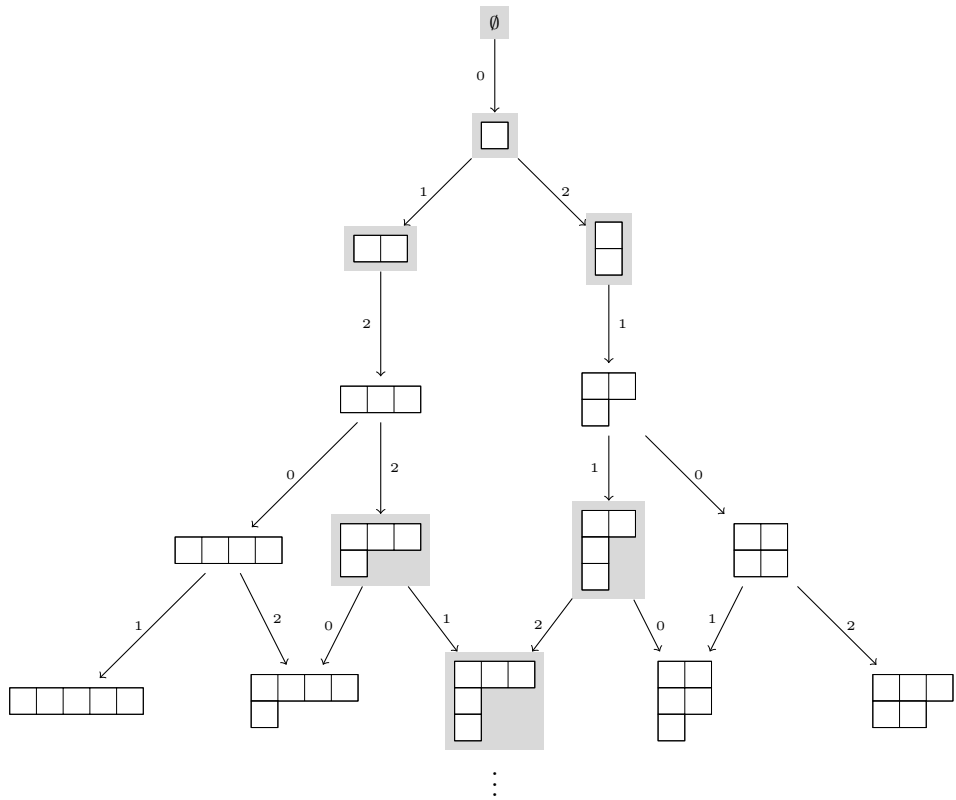
\end{Exa}

The question of the surjectivity of $\sL_\la : W\to \N$
has been solved in this case by Granville and Ono \cite[Theorem 1]{GO1996}.
Indeed, they have proved that
there exists an $(n+1)$-core 
of every given size provided $n\geq 3$.
Using the interpretation of \Cref{crystal_interp}
and using \Cref{cores_orbit},
we can rephrase their result as follows.

\newpage

\begin{Th}\label{GO}
Let $W$ be the Weyl group of type $A_n^{(1)}$.
The map $\sL_{\La_0}: W\to \N$ is surjective
if and only if $n\geq 3$.
\end{Th}

\begin{Rem}
There are generalisations of this crystal realisation
for higher level dominant weights $\la$,
achieved by certain $\ell$-partitions generalising
the $(n+1)$-regular partitions.
These are due to \cite{JMMO1991}, \cite{FLOTW1999},
see also \cite[Chapter 6]{GeckJacon2011}.
In turn, Jacon and Lecouvey gave in \cite{JaconLecouvey2020}
a description of the orbit $\cO(\la)$ resembling
that of \Cref{rem_keys} and generalising \Cref{cores_orbit}.
In particular, the depth of an $\ell$-partition $b$ in the crystal $B(\la)$ is again given by the number of boxes of $b$,
and it would be very interesting to understand for which weights we can generalise \Cref{GO}.
\end{Rem}


\section{Atomic length in affine Weyl groups}\label{sec_affine_AL}

Let $\mathfrak{g}$ be a Kac-Moody algebra of affine type, see \Cref{sec_affine_WG},
so that the Weyl group $ W $ writes  $ W  = T(M) \rtimes W_0$ where 
$W_0$ is the corresponding finite Weyl group
and $M$ is a certain lattice defined from the finite simple roots
(see \Cref{aff_weyl_gps}).
Let $\la=\overline{\la} +\ell \La_0 + z\delta \in P^+$ 
as in \Cref{dom_weight_aff}.
Recall the Coxeter number $h = \sum_{i=0}^n a_i^\vee$.

\begin{Lem}\label{lem_aff_AL}
Let $\be\in M$, $\overline{w} \in W_0$ and set $w=t_\beta \overline{w} $.
We have
$$
\sL_\la( w) 
=
\sL_{\overline{\la}}(\overline{w})
- \ell \ \h (\be)
+ h \left( ( \overline{\la}\mid \overline{w}^{-1}(\be)) + \frac{1}{2}|\beta|^2\ell\right).
$$
\end{Lem}

\begin{proof}
We write
$$\sL_\la( w) = \left\langle \la- w(\la) , \rho^\vee \right\rangle
=\left\langle \la , \rho^\vee \right\rangle 
-\left\langle  w(\la), \rho^\vee \right\rangle$$
and compute both terms.
On the one hand, we have
\begin{align*}
\left\langle \la , \rho^\vee \right\rangle
&
=
\left\langle \overline{\la} + \ell \La_0 +z\delta , \rho^\vee \right\rangle
\\
&
=
\left\langle \overline{\la}, \rho^\vee\right\rangle + \ell \left\langle \La_0  , \rho^\vee \right\rangle+ z \left\langle \delta  , \rho^\vee \right\rangle
\end{align*}

On the other hand, we have
\begin{align*}
\left\langle  w(\la) , \rho^\vee \right\rangle
&
=
\left\langle t_\beta  \overline{w} (\la) ,\rho^\vee\right\rangle
\\
&
= \left\langle t_\beta \left( \overline{w}(\overline{\la} + \ell \La_0 + z\delta ) \right) , \rho^\vee\right\rangle
\\
& 
= \left\langle t_\beta \left( \overline{w}(\overline{\la}) + \ell \La_0 + z\delta \right) , \rho^\vee\right\rangle
\\
& 
=\left\langle  \overline{w} (\overline{\la})+ \ell \La_0 + z\delta +
\ell\beta - \left(( \overline{w}(\overline{\la})\mid\beta)+\frac{1}{2}|\beta|^2\ell \right) \delta 
, \rho^\vee\right\rangle
\text{\quad\quad by (\ref{translations M})}
\\
& 
= \left\langle \overline{w} (\overline{\la} ) , \rho^\vee\right\rangle 
+ \ell \left\langle \La_0  , \rho^\vee\right\rangle 
+ z \left\langle \delta  , \rho^\vee\right\rangle 
+ \ell \left\langle \be  , \rho^\vee\right\rangle 
- \left\langle 
\left(\left(\overline{\la}\mid\overline{w}^{-1}(\be) 
\right)+ \frac{1}{2}|\beta|^2\ell 
\right)\delta  , \rho^\vee\right\rangle 
\end{align*}
Taking the difference yields
\begin{align*}
\sL_\la( w) 
& 
=
\sL_{\overline{\la}}(\overline{w}) 
-\ell \left\langle  \be  , \rho^\vee\right\rangle 
+ \left\langle \left( (\overline{\la}\mid\overline{w}^{-1}(\be))+\frac{1}{2}|\beta|^2\ell \right)\delta   , \rho^\vee\right\rangle 
\\
& 
=
\sL_{\overline{\la}}(\overline{w})
- \ell \ \h ( \be)
+ h \left( ( \overline{\la}\mid\overline{w}^{-1}(\be)) + \frac{1}{2}|\beta|^2\ell\right)
\text{\quad\quad by \Cref{formulas}.}
\end{align*}
\end{proof}

\begin{Cor}\label{aff_AL_lambda0}
For all $ w=t_\beta \overline{w} \in W $, we have
$$\ds \sL_{\La_0} ( w) =
\frac{h}{2} |\beta|^2 - \h(\beta).$$
In particular, $\sL_{\La_0}$ only depends on $\beta$.
\end{Cor}

\begin{proof}
We have taken $\la =\La_0$, so that $\ell=1$ and $\overline{\la}=0$. 
Since $(\La_0 \mid \beta)=0$ (see \Cref{formulas}), 
\Cref{lem_aff_AL} yields precisely the expected formula.
\end{proof}

\begin{Rem}\label{beta_gamma}
One could have chosen to write $w=\overline{w}t_\ga$ instead of $w=t_\be \overline{w}$.
Clearly, $\ga$ and $\be$ are related by the formula $\ga=\overline{w}^{-1}(\be)$,
and since $\overline{w}$ is an isometry, we have $|\ga|=|\be|$.
\end{Rem}

\begin{Rem}\label{cores_paths} 
At this point, let us explain how the specialisation of the affine atomic length
at $\la=\La_0$ appears in various contexts.
\begin{enumerate}
\item Let $W=W(A_n^{(1)})$, so that $W_0=W(A_n)$.
The quotient $W/W_0$ is in bijection with the fundamental chamber $C_0$.
Following \cite{Lascoux2001}, there is a bijection between
the $(n+1)$-cores and the alcoves in $C_0$,
namely the core $w(\emptyset)$ corresponds to the alcove
$A_{w^{-1}}$.
In particular,  if $w=s_{i_r}\cdots s_{i_1}$ is the reduced expression
corresponding to the path $i_1\to\cdots \to i_r$ in $\cO(\lambda)$
starting at $\emptyset$,
then  $w^{-1}=s_{i_1}\cdots s_{i_r}$ is the corresponding reduced
path in $C_0$.
\item Assume again that $W=W(A_n^{(1)})$, in particular $h=n+1$. 
Combining \Cref{cores_orbit} with \Cref{aff_AL_lambda0} immediately gives a formula
for the size of $(n+1)$-core partitions.
Surprisingly, we recover exactly the formula given in \cite[Bijection 2]{GKS1990}, 
used to prove Ramanujan's congruences via the theory of cranks.
Very recently, cranks have been used in \cite{BrunatNath2022} for labelling the
solutions of certain Pell-Fermat equations, which has shed some light on results of Han and Ono \cite{HanOno2011}.
Observations suggest that these results can be extended to other equations
by using the atomic length in other affine types.
\item
In the case of other (untwisted) affine Weyl groups, 
the recent works \cite{ThielWilliams2017} and \cite{STW2021} introduce a generalised notion of size of a core partition,
and obtain generalisations of well-known results on expected sizes of simultaneous core partitions.
Again, this new size statistic turns out to coincide with 
$\sL_{\La_0}$: compare \cite[Section 1.4]{STW2021} with \Cref{aff_AL_lambda0} (up to the change of variables $w \mapsto w^{-1}$). 
Moreover, a formula involving the inversion set (and requiring a distinction between long and short roots) is given. We have been informed \cite{2022NWilliams} that this approach with inversion sets can be adapted for finite types,  yielding the finite atomic length (and,  alternatively,  an analogous statistic
where the squared length of the simple roots is taken into account).
This should also give interesting expectation and variance formulas.
\end{enumerate}
\end{Rem}

\begin{Exa}
Take $W=W(A_2^{(1)})$, so that $h=3$.
Let us compute the first value of  $\sL_{\La_0}(w)$
using \Cref{aff_AL_lambda0}.
These are recorded in the table in \Cref{decomp_AL_aff},
and we also indicate the decomposition $w=t_\beta \overline{w}=\overline{w}t_\ga$,
where we start with elements $w$ with reduced decomposition described in \Cref{cores_paths}.
The reader interested in working out this example might find
it helpful to use the formula $s_0=\tau_{\tal} s_{\tal}$,
where $\tal=\al_1+\al_2$ is the highest root.
\begin{figure}
$$
\begin{array}{@{}l@{\hskip 20pt} @{}l@{\hskip 20pt}  @{}l@{\hskip 20pt} @{}l@{\hskip 20pt} @{}l@{}}
\toprule
\text{reduced expression of $w$}
&
\text{reduced expression of }
\overline{w}
&
\beta
&
\ga
&
\sL_{\La_0}(w)
\\ 
\midrule
e
&
e
&
0
&
0
&
0
\\
s_0
&
s_2s_1s_2
&
\al_1+\al_2
&
-\al_1-\al_2
&
1
\\
s_1 s_0
&
s_2s_1
&
\al_2
&
-\al_1-\al_2
&
2
\\
s_2 s_0
&
s_1s_2
&
\al_1
&
-\al_1-\al_2
&
2
\\
s_2s_1 s_0
&
s_1
&
-\al_2
&
-\al_1-\al_2
&
4
\\
s_1s_2 s_0
&
s_2
&
-\al_1
&
-\al_1-\al_2
&
4
\\
s_2s_1s_2 s_0
&
e
&
-\al_1-\al_2
&
-\al_1-\al_2
&
5
\\
s_0s_2s_1 s_0
&
s_1s_2
&
2\al_1+\al_2
&
-\al_1-2\al_2
&
6
\\
s_0s_1s_2 s_0
&
s_2s_1
&
2\al_2+\al_1
&
-\al_2-2\al_1
&
6
\\
s_0s_2s_1s_2s_0
&
s_1s_2s_1
&
2\al_1+2\al_2
&
-2\al_1-2\al_2
&
8
\\
s_1s_0s_2s_1 s_0
&
s_2
&
\al_2-\al_1
&
-\al_1-2\al_2
&
9
\\
s_2s_0s_1s_2 s_0
&
s_1
&
\al_1-\al_2
&
-2\al_1-\al_2
&
9
\\
\bottomrule
\end{array}
$$
\caption{The first values of $\sL_{\La_0}$ in type $A_2^{(1)}$, computed by determining $\be\in M$
and using \Cref{aff_AL_lambda0}.}
\label{decomp_AL_aff}
\end{figure}
We can compare this table with 
\Cref{exa_cores}, where we computed the first values of $\sL_{\La_0}$ using crystals.
\end{Exa}

\medskip

Combining \Cref{lem_aff_AL} and \Cref{aff_AL_lambda0} yields the following theorem.
Recall the relationship $\be=\overline{w}(\ga)$ explained in \Cref{beta_gamma}.

\begin{Th}\label{thm_aff_AL}
For all $ w=\overline{w}t_\ga\in W $, we have
\begin{align*}
\sL_\la( w) 
&
=
\sL_{\overline{\la}}(\overline{w})
+ 
\ell \sL_{\La_0}( w) + h (\overline{\la}\mid \ga).
\end{align*}
\end{Th}

\Cref{thm_aff_AL} is interesting because
it expresses the affine atomic length in terms of
\begin{itemize}
\item its finite counterpart $\sL_{\overline{\la}}$, which depends only on $\overline{w}$, 
and which we can control in some cases by \Cref{thm_surj}, \Cref{minusc} and \Cref{ideal_necessary},
\item the $\La_0$-atomic length, which depends only on $\be$, and  which we understand by \Cref{GO},
\item the map $w \to h(\overline{\la}\mid\ga )$,
which depends only on $\ga$.
Write $\ga= \sum_{i=1}^n c_i\al_i$
and $\overline{\la} = \sum_{i=1}^n m_i\la_i$.
Then using the formulas of \Cref{sec_affine_WG}, one checks that the linear form
$$\varphi : (c_1,\ldots, c_n) \mapsto \sum_{i=1}^n \frac{a_i^\vee}{a_i} m_i c_i$$
in the variables $c_i$ verifies $\varphi(c_1,\ldots, c_n) = (\overline{\la}\mid \ga)$.
Understanding which integers are representable by the linear form $\varphi$
seems reasonable. For instance,
solutions to the ``coin problem'' give some control over $\varphi$
in the case where the $c_i$'s are nonnegative.
\end{itemize}
We believe that the decomposition of \Cref{thm_aff_AL}
will help in the search for affine ideal weights.


\section*{Acknowledgements}

We thank Cédric Lecouvey and Emily Norton for
stimulating conversations.
We are grateful to Meinolf Geck, Gerhard Hiss and Frank Lübeck for useful discussions regarding the exceptional cases of \Cref{thm_surj} and the use of computer algebra softwares.
Finally, we thank Max Alekseyev for pointing out the reference  \cite{SackUlfarsson2011}
and Nathan Williams for interesting remarks.


\bibliographystyle{alphaurl}
\bibliography{biblio}

\end{document}